%% file: main.tex
\documentclass[12pt,a4paper]{article}
\input{preamble.tex}

\begin{document}
\sloppy
\input{fig/data.tex}
\input{en/frontmatter.tex}
\input{en/r1_intro_path.tex}
\input{en/r2_complexity.tex}
\input{en/r3_fgm.tex}
\input{en/r4_ad.tex}
\input{en/r5_ipm_pi2.tex}
\input{en/r6_traffic.tex}
\input{en/r7_smooth_pd.tex}
\input{en/r8_cubic_gn.tex}
\input{en/r9_cd_zo.tex}
\input{en/r10_univ.tex}
\input{en/r11_tensor_superfast.tex}
\input{en/r12_concl.tex}
\input{en/r13_bib.tex}
\end{document}

%% file: preamble.tex
\usepackage{cmap}
\usepackage[T1]{fontenc}
\usepackage[utf8]{inputenc}
\usepackage[english]{babel}

\usepackage{amsmath,amssymb,amsthm,mathtools}
\usepackage{geometry}
\usepackage{indentfirst}
\usepackage{ifthen}
\usepackage{graphicx}
\usepackage{needspace}
\usepackage{xcolor}
\usepackage{tikz}
\usepackage{pgfplots}
\pgfplotsset{compat=1.18}
\usetikzlibrary{arrows.meta,positioning,calc,patterns,decorations.pathmorphing,
                shapes.geometric,fit,backgrounds}
\usepackage{booktabs}
\usepackage{array}
\usepackage{enumitem}
\usepackage{caption}
\usepackage{hyperref}
\definecolor{lnkcol}{rgb}{0.00,0.00,0.55}
\definecolor{citcol}{rgb}{0.00,0.35,0.10}
\definecolor{urlcol}{rgb}{0.45,0.00,0.10}
\hypersetup{colorlinks=true,linkcolor=lnkcol,citecolor=citcol,urlcolor=urlcol,
            pdftitle={How to Look Inside the Black Box: A Survey of Yu. E. Nesterov's Contributions (On the Award of the 2026 Gauss Prize)},
            pdfauthor={A. V. Gasnikov, A. Yu. Gornov, D. M. Dvinskikh, Yu. V. Dorn, D. I. Kamzolov, Yu. V. Maximov, V. V. Matyukhin, R. F. Hildebrand},bookmarksnumbered=true}

\pgfplotsset{
  umnstyle/.style={
    axis line style={black!70,line width=0.5pt},
    tick align=outside,
    tick label style={font=\footnotesize},
    label style={font=\small},
    title style={font=\small},
    legend style={font=\footnotesize,draw=black!30,fill=white,
                  /tikz/every even column/.append style={column sep=1.1em}},
    legend cell align=left,
    grid=none,
    max space between ticks=40pt,
  },
  legendbelow/.style={
    legend style={at={(0.5,-0.30)},anchor=north,legend columns=#1},
  },
  legendbelow/.default=2,
}

\theoremstyle{plain}
\newtheorem{theorem}{Theorem}[section]
\newtheorem{lemma}[theorem]{Lemma}
\newtheorem{proposition}[theorem]{Proposition}

\theoremstyle{definition}
\newtheorem{definition}[theorem]{Definition}
\newtheorem{example}[theorem]{Example}
\newtheorem{algo}[theorem]{Scheme}
\theoremstyle{remark}
\newtheorem{remark}[theorem]{Remark}

\DeclareMathOperator*{\argmin}{arg\,min}

\DeclareMathOperator*{\Argmin}{Arg\,min}
\DeclareMathOperator{\diag}{diag}
\DeclareMathOperator{\Tr}{tr}
\DeclareMathOperator{\sgn}{sign}

\DeclareMathOperator{\vol}{vol}
\DeclareMathOperator{\dist}{dist}
\DeclareMathOperator{\Lin}{Lin}

\DeclareMathOperator{\rank}{rank}
\DeclareMathOperator{\range}{range}
\newcommand{\R}{\mathbb{R}}
\newcommand{\Exp}{\mathbb{E}}
\newcommand{\Prb}{\mathbb{P}}
\newcommand{\Prob}{\mathop{\mathsf{P}}\nolimits}
\newcommand{\eps}{\varepsilon}
\newcommand{\grad}{\nabla}
\newcommand{\dotp}[2]{\langle #1,#2\rangle}
\newcommand{\nrm}[1]{\|#1\|}
\newcommand{\Cl}{\mathcal{C}}
\newcommand{\Mfr}{\mathfrak{C}}
\newcommand{\bt}{\bar t}
\newcommand{\bff}{\bar f}
\newcommand{\Ss}{\mathbb{S}}
\newcommand{\Oh}{\mathcal{O}}
\newcommand{\Ot}{\widetilde{\mathcal{O}}}
\newcommand{\Om}{\Omega}
\newcommand{\Ndata}{\mathcal{N}}

\newcommand{\practice}{\smallskip\noindent\textit{Practical interpretation.}\ }
\newcommand{\hypoth}{\smallskip\noindent\textit{Research hypothesis.}\ }

\allowdisplaybreaks
\numberwithin{equation}{section}

%% file: fig/data.tex
\def\GDpath{(-2,2.4) (-2.47742,1.31727) (-2.40309,1.27775) (-2.331,1.23942) (-2.26107,1.20223) (-2.19324,1.16617) (-2.12744,1.13118) (-2.06362,1.09725) (-2.00171,1.06433) (-1.94166,1.0324) (-1.88341,1.00143) (-1.82691,0.971384) (-1.7721,0.942242) (-1.71894,0.913975) (-1.66737,0.886556) (-1.61735,0.859959) (-1.56883,0.83416) (-1.52176,0.809135) (-1.47611,0.784861) (-1.43183,0.761315) (-1.38887,0.738476) (-1.34721,0.716322) (-1.30679,0.694832) (-1.26759,0.673987) (-1.22956,0.653768) (-1.19267,0.634154) (-1.15689,0.61513) (-1.12218,0.596676) (-1.08852,0.578776) (-1.05586,0.561412) (-1.02419,0.54457) (-0.993462,0.528233) (-0.963658,0.512386) (-0.934748,0.497014) (-0.906706,0.482104) (-0.879504,0.467641) (-0.853119,0.453612) (-0.827526,0.440003) (-0.8027,0.426803) (-0.778619,0.413999) (-0.75526,0.401579) (-0.732603,0.389532) (-0.710625,0.377846) (-0.689306,0.36651) (-0.668627,0.355515) (-0.648568,0.34485) (-0.629111,0.334504) (-0.610237,0.324469) (-0.59193,0.314735) (-0.574172,0.305293) (-0.556947,0.296134) (-0.540239,0.28725) (-0.524032,0.278633) (-0.508311,0.270274) (-0.493061,0.262165) (-0.47827,0.2543) (-0.463921,0.246671) (-0.450004,0.239271) (-0.436504,0.232093) (-0.423409,0.22513)}
\def\FGMpath{(-2,2.4) (-2.47742,1.31727) (-2.35072,1.2499) (-2.19358,1.16635) (-2.02036,1.07425) (-1.84134,0.979058) (-1.66372,0.884616) (-1.49239,0.793518) (-1.3305,0.707439) (-1.17992,0.627373) (-1.04158,0.55382) (-0.915772,0.486925) (-0.802296,0.426588) (-0.700656,0.372545) (-0.610156,0.324426) (-0.529986,0.281799) (-0.459284,0.244205) (-0.397174,0.211181) (-0.342801,0.18227) (-0.295348,0.157039) (-0.254049,0.13508) (-0.218196,0.116017) (-0.187142,0.0995049) (-0.160299,0.0852323) (-0.13714,0.0729188) (-0.117195,0.0623138) (-0.100045,0.0531949) (-0.0853201,0.0453655) (-0.0726947,0.0386524) (-0.0618832,0.0329039) (-0.0526361,0.0279871) (-0.0447358,0.0237864) (-0.0379931,0.0202013) (-0.0322441,0.0171445) (-0.0273468,0.0145406) (-0.0231787,0.0123243) (-0.019634,0.0104396) (-0.0166219,0.00883804) (-0.0140642,0.00747808) (-0.0118939,0.00632408) (-0.0100534,0.0053455) (-0.00849371,0.00451619) (-0.0071727,0.00381379) (-0.00605449,0.00321923) (-0.00510845,0.00271621) (-0.0043085,0.00229087) (-0.00363241,0.00193138) (-0.00306126,0.0016277) (-0.002579,0.00137128) (-0.00217195,0.00115485) (-0.00182855,0.000972255) (-0.00153894,0.000818268) (-0.0012948,0.000688457) (-0.00108906,0.000579065) (-0.000915752,0.000486914) (-0.000769806,0.000409313) (-0.000646944,0.000343986) (-0.000543549,0.00028901) (-0.000456563,0.000242759) (-0.000383403,0.000203859)}
\def\GDcurve{(0,1) (1,0.143688) (2,0.135196) (3,0.127206) (4,0.119688) (5,0.112614) (6,0.105959) (7,0.0996966) (8,0.0938045) (9,0.0882607) (10,0.0830445) (11,0.0781366) (12,0.0735187) (13,0.0691737) (14,0.0650856) (15,0.061239) (16,0.0576198) (17,0.0542145) (18,0.0510104) (19,0.0479957) (20,0.0451591) (21,0.0424902) (22,0.039979) (23,0.0376163) (24,0.0353932) (25,0.0333014) (26,0.0313333) (27,0.0294815) (28,0.0277392) (29,0.0260998) (30,0.0245573) (31,0.0231059) (32,0.0217404) (33,0.0204555) (34,0.0192466) (35,0.0181091) (36,0.0170389) (37,0.0160319) (38,0.0150844) (39,0.0141929) (40,0.0133541) (41,0.0125649) (42,0.0118223) (43,0.0111236) (44,0.0104662) (45,0.00984764) (46,0.00926565) (47,0.00871805) (48,0.00820281) (49,0.00771802) (50,0.00726189) (51,0.00683271) (52,0.0064289) (53,0.00604895) (54,0.00569146) (55,0.00535509) (56,0.00503861) (57,0.00474082) (58,0.00446064) (59,0.00419702) (60,0.00394897) (61,0.00371559) (62,0.003496) (63,0.00328938) (64,0.00309498) (65,0.00291207) (66,0.00273997) (67,0.00257803) (68,0.00242567) (69,0.00228231) (70,0.00214743) (71,0.00202052) (72,0.0019011) (73,0.00178875) (74,0.00168303) (75,0.00158357) (76,0.00148998) (77,0.00140192) (78,0.00131907) (79,0.00124111) (80,0.00116776) (81,0.00109875) (82,0.00103381) (83,0.000972711) (84,0.000915224) (85,0.000861134) (86,0.000810241) (87,0.000762356) (88,0.000717301) (89,0.000674908) (90,0.000635021) (91,0.000597492) (92,0.00056218) (93,0.000528955) (94,0.000497694) (95,0.00046828) (96,0.000440605) (97,0.000414565) (98,0.000390064) (99,0.000367011) (100,0.000345321) (101,0.000324913) (102,0.00030571) (103,0.000287643) (104,0.000270643) (105,0.000254648) (106,0.000239598) (107,0.000225438) (108,0.000212115) (109,0.000199579) (110,0.000187784) (111,0.000176686) (112,0.000166243) (113,0.000156418) (114,0.000147174) (115,0.000138476) (116,0.000130292) (117,0.000122592) (118,0.000115347) (119,0.00010853) (120,0.000102116) (121,9.60806e-05) (122,9.04023e-05) (123,8.50595e-05) (124,8.00325e-05) (125,7.53026e-05) (126,7.08522e-05) (127,6.66648e-05) (128,6.27249e-05) (129,5.90179e-05) (130,5.55299e-05) (131,5.22481e-05) (132,4.91602e-05) (133,4.62549e-05) (134,4.35212e-05) (135,4.09491e-05) (136,3.8529e-05) (137,3.62519e-05) (138,3.41095e-05) (139,3.20936e-05) (140,3.01969e-05) (141,2.84122e-05) (142,2.67331e-05) (143,2.51531e-05) (144,2.36666e-05) (145,2.22679e-05) (146,2.09519e-05) (147,1.97136e-05) (148,1.85485e-05) (149,1.74523e-05) (150,1.64209e-05) (151,1.54504e-05) (152,1.45373e-05) (153,1.36781e-05) (154,1.28698e-05) (155,1.21092e-05) (156,1.13935e-05) (157,1.07201e-05) (158,1.00866e-05) (159,9.49047e-06) (160,8.92958e-06) (161,8.40184e-06) (162,7.90529e-06) (163,7.43809e-06) (164,6.9985e-06) (165,6.58489e-06) (166,6.19572e-06) (167,5.82955e-06) (168,5.48503e-06) (169,5.16086e-06) (170,4.85586e-06) (171,4.56887e-06) (172,4.29885e-06) (173,4.04479e-06) (174,3.80574e-06) (175,3.58083e-06) (176,3.3692e-06) (177,3.17008e-06) (178,2.98273e-06) (179,2.80645e-06) (180,2.64059e-06) (181,2.48453e-06) (182,2.33769e-06) (183,2.19953e-06) (184,2.06954e-06) (185,1.94723e-06) (186,1.83215e-06) (187,1.72387e-06) (188,1.62199e-06) (189,1.52613e-06) (190,1.43594e-06) (191,1.35107e-06) (192,1.27122e-06) (193,1.19609e-06) (194,1.12541e-06) (195,1.05889e-06) (196,9.96313e-07) (197,9.37431e-07) (198,8.82029e-07) (199,8.29901e-07) (200,7.80854e-07) (201,7.34706e-07) (202,6.91284e-07) (203,6.5043e-07) (204,6.11989e-07) (205,5.75821e-07) (206,5.4179e-07) (207,5.0977e-07) (208,4.79642e-07) (209,4.51296e-07) (210,4.24624e-07) (211,3.99529e-07) (212,3.75917e-07) (213,3.537e-07) (214,3.32796e-07) (215,3.13128e-07) (216,2.94622e-07) (217,2.7721e-07) (218,2.60827e-07) (219,2.45412e-07) (220,2.30908e-07) (221,2.17261e-07) (222,2.04421e-07) (223,1.9234e-07) (224,1.80973e-07) (225,1.70277e-07) (226,1.60214e-07) (227,1.50745e-07) (228,1.41836e-07) (229,1.33454e-07) (230,1.25567e-07) (231,1.18146e-07) (232,1.11163e-07) (233,1.04593e-07) (234,9.84119e-08) (235,9.25958e-08) (236,8.71234e-08) (237,8.19744e-08) (238,7.71297e-08) (239,7.25713e-08) (240,6.82824e-08) (241,6.42469e-08) (242,6.04499e-08) (243,5.68773e-08) (244,5.35159e-08) (245,5.03531e-08) (246,4.73772e-08) (247,4.45772e-08) (248,4.19427e-08) (249,3.94639e-08) (250,3.71316e-08) (251,3.49371e-08) (252,3.28723e-08) (253,3.09296e-08) (254,2.91016e-08) (255,2.73817e-08) (256,2.57635e-08) (257,2.42408e-08) (258,2.28082e-08) (259,2.14602e-08) (260,2.01919e-08) (261,1.89986e-08) (262,1.78758e-08) (263,1.68193e-08) (264,1.58253e-08) (265,1.489e-08) (266,1.401e-08) (267,1.3182e-08) (268,1.2403e-08) (269,1.167e-08) (270,1.09803e-08) (271,1.03313e-08) (272,9.72075e-09)}
\def\FGMcurve{(0,1) (1,0.143688) (2,0.129367) (3,0.11265) (4,0.0955609) (5,0.079376) (6,0.0648011) (7,0.0521419) (8,0.041443) (9,0.032593) (10,0.0253986) (11,0.0196334) (12,0.0150692) (13,0.0114929) (14,0.00871571) (15,0.00657583) (16,0.00493837) (17,0.00369302) (18,0.00275109) (19,0.00204215) (20,0.00151097) (21,0.00111459) (22,0.000819902) (23,0.000601564) (24,0.000440303) (25,0.000321545) (26,0.000234322) (27,0.000170421) (28,0.000123716) (29,8.96534e-05) (30,6.48617e-05) (31,4.68523e-05) (32,3.37933e-05) (33,2.43401e-05) (34,1.7508e-05) (35,1.25776e-05) (36,9.02486e-06) (37,6.46821e-06) (38,4.63076e-06) (39,3.31183e-06) (40,2.36619e-06) (41,1.68895e-06) (42,1.20444e-06) (43,8.58175e-07) (44,6.10942e-07) (45,4.34584e-07) (46,3.08894e-07) (47,2.19393e-07) (48,1.55712e-07) (49,1.10439e-07) (50,7.82768e-08) (51,5.54452e-08) (52,3.92488e-08) (53,2.77669e-08) (54,1.96326e-08) (55,1.38734e-08) (56,9.7984e-09)}
\def\GDiters{272}
\def\FGMiters{56}
\def\Kappa{33}
\def\SqKappa{5.8}
\def\LevAa{1.2821}
\def\LevAb{7.4021}
\def\LevBa{0.8600}
\def\LevBb{4.9655}
\def\LevCa{0.5439}
\def\LevCb{3.1404}
\def\LevDa{0.3007}
\def\LevDb{1.7359}
\def\LevEa{0.1517}
\def\LevEb{0.8758}
\def\LevRot{62.0}
\def\ScaleGD{(10,77) (30,222) (100,687) (300,1904) (1000,5754) (3000,15619) (10000,46049) (30000,121674)}
\def\ScaleFGM{(10,28) (30,49) (100,88) (300,145) (1000,247) (3000,398) (10000,661) (30000,1037)}
\def\RefLin{(10,63.4667) (30,190.4) (100,634.667) (300,1904) (1000,6346.67) (3000,19040) (10000,63466.7) (30000,190400)}
\def\RefSqrt{(10,26.4733) (30,45.853) (100,83.7158) (300,145) (1000,264.733) (3000,458.53) (10000,837.158) (30000,1450)}
\def\CentralPath{(0.0982,-0.4338) (0.0942,-0.4399) (0.0901,-0.4460) (0.0860,-0.4519) (0.0819,-0.4578) (0.0778,-0.4636) (0.0737,-0.4693) (0.0695,-0.4749) (0.0654,-0.4805) (0.0612,-0.4860) (0.0570,-0.4913) (0.0529,-0.4966) (0.0487,-0.5019) (0.0445,-0.5070) (0.0403,-0.5121) (0.0361,-0.5171) (0.0318,-0.5220) (0.0276,-0.5268) (0.0234,-0.5316) (0.0191,-0.5363) (0.0149,-0.5409) (0.0106,-0.5455) (0.0064,-0.5499) (0.0021,-0.5543) (-0.0021,-0.5587) (-0.0064,-0.5630) (-0.0106,-0.5672) (-0.0149,-0.5713) (-0.0192,-0.5754) (-0.0234,-0.5794) (-0.0277,-0.5834) (-0.0319,-0.5873) (-0.0362,-0.5911) (-0.0405,-0.5949) (-0.0447,-0.5986) (-0.0490,-0.6023) (-0.0532,-0.6059) (-0.0574,-0.6094) (-0.0617,-0.6129) (-0.0659,-0.6164) (-0.0659,-0.6164) (-0.0789,-0.6267) (-0.0928,-0.6372) (-0.1078,-0.6481) (-0.1237,-0.6592) (-0.1407,-0.6704) (-0.1589,-0.6819) (-0.1781,-0.6934) (-0.1985,-0.7051) (-0.2200,-0.7167) (-0.2426,-0.7284) (-0.2662,-0.7400) (-0.2908,-0.7516) (-0.3163,-0.7630) (-0.3426,-0.7742) (-0.3695,-0.7852) (-0.3969,-0.7960) (-0.4246,-0.8065) (-0.4526,-0.8167) (-0.4805,-0.8266) (-0.5082,-0.8361) (-0.5356,-0.8453) (-0.5624,-0.8542) (-0.5886,-0.8626) (-0.6140,-0.8707) (-0.6385,-0.8785) (-0.6621,-0.8858) (-0.6846,-0.8928) (-0.7060,-0.8995) (-0.7262,-0.9058) (-0.7454,-0.9118) (-0.7634,-0.9174) (-0.7804,-0.9227) (-0.7963,-0.9277) (-0.8111,-0.9325) (-0.8249,-0.9369) (-0.8378,-0.9411) (-0.8498,-0.9450) (-0.8609,-0.9487) (-0.8712,-0.9522) (-0.8808,-0.9554) (-0.8897,-0.9585) (-0.8979,-0.9613) (-0.9055,-0.9640) (-0.9126,-0.9665) (-0.9191,-0.9688) (-0.9251,-0.9710) (-0.9307,-0.9730) (-0.9359,-0.9749) (-0.9406,-0.9766) (-0.9451,-0.9783) (-0.9491,-0.9798) (-0.9529,-0.9812) (-0.9564,-0.9826) (-0.9597,-0.9838) (-0.9626,-0.9850) (-0.9654,-0.9860) (-0.9680,-0.9870) (-0.9703,-0.9880) (-0.9725,-0.9888) (-0.9746,-0.9896) (-0.9764,-0.9904) (-0.9782,-0.9911) (-0.9798,-0.9917) (-0.9813,-0.9923) (-0.9826,-0.9929) (-0.9839,-0.9934) (-0.9851,-0.9939) (-0.9862,-0.9943) (-0.9872,-0.9947) (-0.9882,-0.9951) (-0.9890,-0.9954) (-0.9898,-0.9958) (-0.9906,-0.9961) (-0.9913,-0.9964) (-0.9919,-0.9966) (-0.9925,-0.9969) (-0.9930,-0.9971) (-0.9936,-0.9973) (-0.9940,-0.9975)}
\def\PolyVerts{(-1.0000,-1.0000) (1.0000,-1.0000) (1.0000,0.5500) (0.4500,1.0000) (-0.3500,1.0000) (-1.0000,0.1000)}
\def\DikAx{0.0982}
\def\DikAy{-0.4338}
\def\DikAa{0.6044}
\def\DikAb{0.4890}
\def\DikAr{-177.93}
\def\DikBx{-0.2908}
\def\DikBy{-0.7516}
\def\DikBa{0.5527}
\def\DikBb{0.2425}
\def\DikBr{-178.90}
\def\DikCx{-0.6609}
\def\DikCy{-0.8887}
\def\DikCa{0.1370}
\def\DikCb{0.0512}
\def\DikCr{-179.86}
\def\Fzero{-0.087975}
\def\Gzero{-0.5730}
\def\Hzero{-1.7300}
\def\NewtStep{-0.3312}
\def\NewtPt{-0.0312}
\def\NewtVal{-0.000974}
\def\CubStep{0.4852}
\def\CubPt{0.7852}
\def\CubVal{-0.521478}

%% file: en/frontmatter.tex
\thispagestyle{empty}
\begin{center}
{\small UDC 519.853.62+519.856}\\[14pt]
{\Large\bfseries How to Look Inside the Black Box:\\[4pt]
A Survey of Yu.\,E.~Nesterov's Contributions}\\[8pt]
{\large\bfseries (On the Award of the 2026 Gauss Prize)}\\[16pt]
{\large A.\,V.~Gasnikov, A.\,Yu.~Gornov, D.\,M.~Dvinskikh, Yu.\,V.~Dorn,
D.\,I.~Kamzolov, Yu.\,V.~Maximov, V.\,V.~Matyukhin, R.\,F.~Hildebrand}%
\begin{NoHyper}\let\thefootnote\relax\footnotetext{This research was supported by the Ministry of Science and Higher Education of the Russian Federation (state assignment), project no.~FSMG-2024-0011. Modern AI models were used in part in preparing the text of this article.}\end{NoHyper}
\end{center}

\vspace{6pt}
{\small
\hspace*{1.5em}\emph{On July 23, 2026, at the opening of the International Congress of Mathematicians in Philadelphia, the International Mathematical Union announced the award of the Gauss Prize to Yurii Evgenievich Nesterov---an award for mathematical work with significant applications outside mathematics. This article explains the mathematical content of the recognized work of Yu.\,E.~Nesterov and traces its influence on modern computational practice. The exposition is based on a concept that Yu.\,E.~Nesterov himself regards as one of his main results: the classical complexity theory of convex optimization is oriented toward a black-box oracle, but its pessimistic bounds can be overcome by using a priori information about the internal structure of the objective function. We trace Yu.\,E.~Nesterov's scientific career, then introduce the language of oracle complexity together with a classification of objective function models (geometric, algebraic, and structural), and subsequently discuss his main results in chronological order: the fast gradient method, fast automatic differentiation, the theory of self-concordant barriers and polynomial-time interior-point methods, the $\pi/2$ theorem, the stable dynamics model of traffic flows, the smoothing technique, primal-dual methods, cubic regularization of Newton's method, the modified Gauss--Newton scheme, coordinate descent methods and methods for huge-scale problems, gradient-free methods, universal gradient methods, implementable tensor methods, and ``superfast'' methods. Each topic is accompanied by historical background, precise statements specifying the problem class, oracle, and type of guarantee, examples, and a discussion of the current state of the art.}

\smallskip
{\small\textbf{Bibliography:} 194 references.}

}

\vspace{10pt}
{\small\tableofcontents}
\clearpage

%% file: en/r1_intro_path.tex
\section{Introduction}
\label{sec:intro}

\subsection{The Gauss Prize}
The Carl Friedrich Gauss Prize is an award of the International Mathematical Union, presented every four years for outstanding contributions to the application of mathematics in other scientific and applied disciplines. Until recently, the list of recipients remained very short: K.~It\^o (2006), Y.~Meyer (2010), S.~Osher (2014), D.~Donoho (2018), and E.~Lieb (2022). The sixth recipient was Yurii Evgenievich Nesterov: the decision was announced on July 23, 2026, at the opening ceremony of the International Congress of Mathematicians in Philadelphia.\footnote{The official award page and the prize committee's citation:
\url{https://www.mathunion.org/imu-awards/carl-friedrich-gauss-prize}.
Here and below, biographical information and award details are taken from the official profiles and announcements of the respective organizations.}
The prize committee particularly noted his fundamental work in mathematical optimization, which laid the theoretical and algorithmic foundations for applied computational disciplines and data analysis methods, making it possible to overcome previous computational barriers and solve problems on a new scale.

It is not easy to cover the mathematical content of such a broad and diverse range of results in a single survey. The work of Yu.\,E.~Nesterov spans many different areas: from complexity theory and interior-point methods to semidefinite relaxations, traffic equilibria, automatic differentiation, and high-order algorithms. At the same time, the best-known tools he created are now so deeply integrated into computational practice that they are often perceived as basic techniques. Nevertheless, behind all this diversity lies a coherent and fundamental research program, with which we begin.

\subsection{The main idea: looking inside the black box}
\label{sec:lookinside}

By the mid-1980s, the theory of numerical methods for convex optimization appeared to be complete, thanks to the work of A.\,S.~Nemirovski and D.\,B.~Yudin~\cite{NemYud1979}. Their framework rigorously fixed an information model in the form of a black-box oracle that, in response to a query, returned the function value and its (sub)gradient at a given point. For the principal problem classes, the authors mathematically established complexity lower bounds---the minimum number of oracle calls needed to attain a prescribed accuracy. Finally, for each class, optimal methods were proposed that were guaranteed to attain these theoretical limits. Moreover, this information model precisely matched the style of writing optimization packages accepted at the time: the user was asked to write a subroutine for computing the function value and gradient, whose internal workings were hidden from the package itself~\cite{Evtushenko1982}. It seemed that all the important questions had been asked and almost all had been answered~\cite{NesDiss2013}.

By the early 2000s, as noted in~\cite[Introduction]{NesDiss2013}, the limitations of this unified classification based on the black-box model had become apparent. The reason was unexpected. The range of real applications of polynomial-time interior-point methods virtually coincided with that of methods operating in the black-box model, and wherever interior-point methods were applicable, they generally outperformed, both theoretically and in practice, the unimprovable optimal methods. This apparent contradiction is resolved by one key point: constructing an effective (self-concordant) barrier requires access to the internal structure of the problem, which effectively takes interior-point methods outside the classical black-box model.
The lower bound is not violated; it simply applies to a different, broader class of methods.

Hence the program: some of the most pessimistic limits of complexity theory can be pushed back, and it is clear how---one must exploit the structure of the problem. In the terminology of Yu.\,E.~Nesterov himself~\cite{NesDiss2013}, this is structural optimization, based on refining the concept of an objective function model: alongside an oracle that supplies numbers, the method is given an analytical description of how these numbers are structured. We will call this program looking into the black box. It runs throughout the article and explains why the section headings below include both ``optimal method'' and ``faster than the lower bound'': in the former case, the information model remains unchanged and its limit is attained; in the latter, the model changes and the limit shifts.

Below (\S\,\ref{sec:models}), we present the classification of objective function models from~\cite{NesDiss2013} (geometric, algebraic, and structural) and indicate which section of the article corresponds to each. We note at once that the program did not amount merely to an analysis of individual examples: by the time of the doctoral dissertation defense, the outcome was formulated quite sharply---the theory of convex optimization was no longer covered by a unified system of lower bounds for all the structural models in use, and developing new principles for classifying optimization problems remains an open direction~\cite{NesDiss2013}.

For a broad international audience, this program was formulated in Yu.\,E.~Nesterov's book of the early 2000s~\cite{NesBook2004}, which set the mathematical standard of exposition in this field; its role as a model that shaped the style of modern work on the complexity of convex optimization is explicitly discussed by S.~Bubeck in his widely known survey~\cite{Bubeck2015}. The 2018 monograph~\cite{NesBook2018} is a substantially expanded version of the same book; it overlaps considerably with the 2013 doctoral dissertation~\cite{NesDiss2013}, which was subsequently published unchanged under the same title~\cite{NesBook2024RU}.

\subsection{Three layers of exposition}

Our exposition will be organized around three interrelated aspects.

\smallskip
\emph{Foundations.} Before constructing a method, we must agree on what ``better'' means. The answer is provided by oracle complexity theory (\emph{black-box complexity} in the English-language literature): a problem class is fixed, a class of methods using the same type of available information is fixed, and a lower bound is proved for this pair of classes. Only then does it make sense to speak of optimality. This language, together with the classification of objective function models, is presented in \S\,\ref{sec:complexity}.

\smallskip
\emph{Methods.} The work of Yu.\,E.~Nesterov constitutes a systematic implementation of this methodology and a profound investigation of the computational possibilities it opens up. Sometimes one can construct a method that attains the lower bound {in the same} information model: this is how the fast gradient method of 1983 arose (\S\,\ref{sec:fgm}). Much more often---and this is precisely what looking into the black box means---one can identify an objective function model under which the pessimistic bound no longer applies: this is how self-concordant barriers (\S\,\ref{sec:ipm}), the smoothing technique (\S\,\ref{sec:smooth}), primal-dual methods (\S\,\ref{sec:pd}), coordinate descent methods (\S\,\ref{sec:cd}), tensor methods (\S\,\ref{sec:tensor}), and ``superfast'' methods (\S\,\ref{sec:superfast}) arose. A separate line of work is the search for {the right notion of regularity} under which a problem suddenly becomes solvable: self-concordance (\S\,\ref{sec:ipm}), the H\"older scale (\S\,\ref{sec:univ}), gradient domination (\S\,\ref{sec:cubic}), and relative smoothness (\S\,\ref{sec:superfast}).

\smallskip
\emph{Applications.} Whenever a method has been obtained, the question of its applicability arises, and here Yu.\,E.~Nesterov generally pursued the matter through to a specific applied model: semidefinite relaxations of nonconvex problems (\S\,\ref{sec:pi2}), traffic equilibria (\S\,\ref{sec:traffic}), and huge-scale problems (\S\,\ref{sec:cd}). The fast gradient method and the dual averaging method, in turn, have become components of almost every modern scheme for training machine learning models.

\subsection{Organization}

In each section, we try to follow the same pattern: the background of the problem; a precise statement explicitly specifying the \emph{function class}, the \emph{available oracle}, the \emph{normalization of the initial error}, and the \emph{type of guarantee} (deterministic or in expectation); the idea of the proof or a complete short proof; an illustrative example; a discussion of what has changed; and the current state of the art. Statements that go beyond what has been proved---about neural networks, software packages, or empirical observations---are explicitly marked as practical interpretations or research hypotheses.

Sections \S\S\,\ref{sec:fgm}\,--\,\ref{sec:superfast} follow the chronology of the first publications on each topic; the year is indicated in the heading and, where there is a substantial discrepancy, so is the year of the CORE preprint. The only exception is \S\,\ref{sec:pd}: primal-dual methods (journal publication in 2009, preprint in 2005) are presented immediately after the smoothing technique, since these two topics are technically inseparable. Two sections---on self-concordant barriers (\S\,\ref{sec:ipm}) and on traffic models (\S\,\ref{sec:traffic})---are treated in greater detail than the others: they provide two extended illustrations of how optimization theory enters applications.

We make no claim to completeness. In particular, Yu.\,E.~Nesterov's work on quasi-variational inequalities, conditional gradient methods, ellipsoidal approximation of convex bodies, and the analysis of ranking algorithms lies outside the scope of this survey; some of it is mentioned along the way without detailed treatment.

\clearpage
\section{Scientific career}
\label{sec:path}

\noindent
\begin{minipage}[t]{0.66\textwidth}
\vspace{0pt}
Yurii Evgenievich Nesterov was born on January 25, 1956, in Moscow. In 1977, he graduated from the Faculty of Computational Mathematics and Cybernetics at Moscow State University with a degree in applied mathematics, and in the same year joined the Central Economics and Mathematics Institute of the USSR Academy of Sciences, where he worked until 1992 in a laboratory devoted to optimization theory and numerical methods. He defended his Candidate of Sciences dissertation in 1984 at the Institute of Control Sciences of the USSR Academy of Sciences; his supervisor was Boris Teodorovich Polyak.
\end{minipage}\hfill
\begin{minipage}[t]{0.30\textwidth}
\vspace{0pt}
\centering
\includegraphics[width=\linewidth]{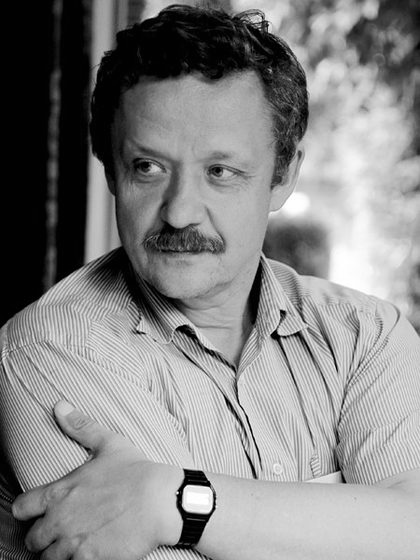}
\captionof*{figure}{Yu.\,E.~Nesterov.}
\label{fig:photo}
\end{minipage}
\par\medskip

This early experience of computational practice had a marked influence on the style of his subsequent work. The first papers of the young CEMI researcher were devoted not to building theories but to \emph{computation}: interactive systems for analyzing algorithms, comparing implementations, and understanding why a method that looks good on paper loses on real examples. Many of Yu.\,E.~Nesterov's works are characterized by simultaneous attention to iteration bounds and to the implementability of a single iteration.

B.\,T.~Polyak's influence can be traced in at least two respects, both of which we will encounter below: the \emph{heavy-ball method} (1964)~\cite{Polyak1964}---the first inertial scheme, whose analysis gave rise to the fast gradient method of 1983 (\S\,\ref{sec:fgm})---and the \emph{gradient domination condition} (1963)~\cite{Polyak1963}, now known as the Polyak--\L{}ojasiewicz condition; it will appear in \S\,\ref{sec:cubic} as the assumption under which the cubically regularized Newton method converges superlinearly. The joint paper by teacher and student on cubic regularization~\cite{NesPolyak2006} appeared more than twenty years after the dissertation defense.

A second, and apparently even more decisive, collaboration was with Arkadi Semenovich Nemirovski. It produced two major groups of results: \emph{optimal methods for smooth convex minimization}~\cite{NemNes1985}, which already considered intermediate (H\"older) smoothness classes and the composite structure of the problem (\S\,\ref{sec:univ}), and the theory of \emph{self-concordant barriers}, developed in 1988\,--\,1994 and presented in the monograph~\cite{NesNem1994} (\S\,\ref{sec:ipm})---it was this book, as noted in the prize committee's citation, that established the author's reputation in optimization.

Since 1993, Yu.\,E.~Nesterov has worked at the Center for Operations Research and Econometrics (CORE) at the Catholic University of Louvain; he was a professor at CORE in 1993\,--\,2023 and currently holds the status of professor emeritus. He is also a professor in the School of Data Science at the Chinese University of Hong Kong, Shenzhen. The vast majority of the results discussed below were published during this period, typically following a pattern in which a CORE Discussion Paper appeared two to five years before the journal publication, so that both dates occur in the literature.

\Needspace{4\baselineskip}
Yu.\,E.~Nesterov maintains friendly and scientific ties with colleagues in Russia. In 2013, he defended his Doctor of Sciences dissertation, ``Algorithmic Convex Optimization,'' at MIPT---367~pages, essentially a coherent exposition of the entire program described below~\cite{NesDiss2013}; the work was carried out at MIPT's Laboratory of Structural Methods of Data Analysis in Predictive Modeling (PreMoLab). The dissertation was subsequently published unchanged under the same title~\cite{NesBook2024RU}. Several generations of specialists have studied from his books published in Russian~\cite{NesRus1989,NesRus2010}; the 1989 book~\cite{NesRus1989} was difficult to obtain for a long time and was reissued as a second, unchanged edition, so it is appropriate to cite both. A substantial proportion of the works cited below in describing the current state of a topic were carried out by his students, coauthors, and coauthors' students; in the textbook~\cite{GasnBook}, a substantial part of the material presented (the universal method, primal-dual properties, and coordinate descent methods) is developed to the level of a university course.

In addition to the Gauss Prize, Yu.\,E.~Nesterov is a recipient of the Dantzig Prize (2000), awarded jointly by the Mathematical Optimization Society and SIAM: the citation recognizes his contributions to the theory of interior-point methods, including the concept of self-concordant functions, and his earlier work on convex optimization and automatic differentiation. Further awards include the John von Neumann Theory Prize of the Institute for Operations Research and the Management Sciences (INFORMS) (2009, shared with Yinyu Ye), the Gold Medal of the Association of European Operational Research Societies (EURO) (2016, shared with M.~Queyranne), the INFORMS Lanchester Prize (2022) for the monograph~\cite{NesBook2018}, and the World Laureates Association (WLA) Prize (2023, jointly with A.\,S.~Nemirovski) for work on the theory of convex optimization. He was elected a member of Academia Europaea (2021) and a foreign member of the US National Academy of Sciences (2022).

\begin{figure}[htbp]
\centering
\begin{tikzpicture}[x=1cm,y=1cm,
  result/.style={anchor=north west,font={\small\hyphenpenalty=10000\exhyphenpenalty=10000\relax},
    text width=5.15cm,align=left,inner xsep=0pt,inner ysep=3pt},
  meaning/.style={anchor=north west,font={\small\hyphenpenalty=10000\exhyphenpenalty=10000\relax},
    text width=6.05cm,align=left,inner xsep=0pt,inner ysep=3pt}]
  \node[anchor=south east,font=\footnotesize\bfseries,black!65]
    at (1.83,0.18) {Years};
  \node[anchor=south west,font=\footnotesize\bfseries,black!65]
    at (2.38,0.18) {Result};
  \node[anchor=south west,font=\footnotesize\bfseries,black!65]
    at (8.08,0.18) {What it enables};
  \foreach \i/\yr/\name/\effect in {%
    0/1983/{Fast gradient method \mbox{(\S\,\ref{sec:fgm})}}/{Attain the optimal order of convergence with the same oracle},
    1/1984/{Automatic differentiation \mbox{(\S\,\ref{sec:ad})}}/{Compute the gradient at a cost of the same order as the value},
    2/{1986--1987}/{Lexicographic differentiation \mbox{(\S\,\ref{sec:lexdiff})}}/{Make derivatives consistent throughout a nonsmooth computational graph},
    3/{1988--1997}/{Self-concordant barriers \mbox{(\S\,\ref{sec:ipm})}}/{Choose a step in the local geometry of the constraints},
    4/1998/{$\pi/2$ theorem \mbox{(\S\,\ref{sec:pi2})}}/{Bound the quality of rounding a semidefinite relaxation},
    5/2003/{Stable dynamics \mbox{(\S\,\ref{sec:traffic})}}/{Describe queues through dual variables},
    6/{2003--2009}/{Smoothing and primal-dual methods \mbox{(\S\S\,\ref{sec:smooth}--\ref{sec:pd})}}/{Exploit the max structure and obtain an accuracy certificate},
    7/{2006--2007}/{Cubic regularization; Gauss--Newton scheme \mbox{(\S\S\,\ref{sec:cubic}--\ref{sec:gn})}}/{Control the step using an upper model},
    8/{2010--2014}/{Coordinate descent methods; huge-scale \mbox{(\S\,\ref{sec:cd})}}/{Reduce the cost of updates through sparsity},
    9/{2011--2017}/{Gradient-free methods \mbox{(\S\,\ref{sec:zo})}}/{Estimate the gradient of a smoothed function from two values},
    10/{2013--2015}/{Universal methods \mbox{(\S\,\ref{sec:univ})}}/{Adapt to an unknown smoothness class},
    11/{2018--2021}/{Tensor and ``superfast'' methods \mbox{(\S\S\,\ref{sec:tensor}--\ref{sec:superfast})}}/{Exploit higher smoothness and convex high-order models},
    12/2026/{Gauss Prize}/{Recognition of the contribution of mathematical optimization to applications}}
  {
    \ifnum\i=0
      \coordinate (top\i) at (0,0);
    \else
      \pgfmathtruncatemacro{\prev}{\i-1}
      \coordinate (top\i) at ($(row\prev.south west)+(0,-0.13)$);
    \fi
    \node[result] (result\i) at (2.38,0 |- top\i) {\name};
    \node[meaning] (meaning\i) at (8.08,0 |- top\i) {\effect};
    \node[fit=(result\i)(meaning\i),inner sep=0pt] (row\i) {};
    \coordinate (mark\i) at (2.03,0 |- row\i.center);
    \node[anchor=east,font=\footnotesize\bfseries,blue!50!black]
      at (1.83,0 |- row\i.center) {\yr};
    \fill[blue!55!black] (mark\i) circle (0.060);
    \draw[->,>=Latex,black!35,line width=0.45pt]
      (7.67,0 |- row\i.center)--(7.94,0 |- row\i.center);
  }
  \begin{scope}[on background layer]
    \draw[blue!20!black,line width=0.6pt]
      (mark0)--(mark12);
  \end{scope}
\end{tikzpicture}
\caption{The main strands of Yu.\,E.~Nesterov's research program: from attaining the limits of the oracle model to exploiting problem structure and applications. The intervals group related works discussed in the corresponding sections; vertical distances are schematic.}
\label{fig:timeline}
\end{figure}

%% file: en/r2_complexity.tex

\section{Foundations: the black box and the objective-function model}
\label{sec:complexity}

\subsection{Why an information model is needed}

Consider the problem
\begin{equation}\label{eq:prob}
   \min_{x\in Q}\ f(x),
\end{equation}
where $Q\subseteq\R^n$ is a closed convex set and $f$ is a convex
function, and ask: how many operations are needed to find a solution with
accuracy $\eps$? In this form, the question has no answer. If the function
is given by the formula $f(x)=\nrm{x-a}^2$ and the algorithm used to solve
the problem returns the fixed point $a$ for every input, the answer is zero
iterations. If, however, the function is given by a program of ten thousand
lines, or the algorithm searches for solutions by enumerating feasible
points, the answer will be different. Thus one must specify not only the
problem class, but also the \emph{mode of access} to the problem and the
corresponding class of methods.

The model proposed by A.\,S.~Nemirovski and
D.\,B.~Yudin~\cite{NemYud1979} works as follows. The method has no access to
the analytic expression for the function. It can only submit points $x\in Q$
to an \emph{oracle} and receive local information in return: the value $f(x)$
(a zeroth-order oracle); the value and a (sub)gradient (first order);
the value, gradient, and Hessian (second order); and so on. The function
itself remains a \emph{black box}. We call this scheme the \emph{black-box
model}, and say that a method using nothing beyond the oracle's responses
\emph{operates in the black-box model}. The measure of computational effort
is the \emph{number of oracle calls} needed to attain accuracy $\eps$ on the
worst-case problem in the class.

Two points should be emphasized at once. First, oracle complexity is
\emph{not} arithmetic complexity: it does not account for the cost of an
iteration. For gradient methods this separation is justified, since
computing the gradient is generally the most expensive part of an iteration
(\S\,\ref{sec:ad}); in \S\S\,\ref{sec:cd},~\ref{sec:superfast}, however,
we will have to count arithmetic operations separately. Second, the problem
class must include a \emph{normalization of the initial error}: otherwise
the complexity is infinite, since the solution can be moved arbitrarily far
away.

\subsection{Problem classes and formulation of the question}

Let $\Cl^{1,1}_L$ denote the class of differentiable convex functions on
$\R^n$ with an $L$-Lipschitz gradient,
$\nrm{\grad f(x)-\grad f(y)}_2\le L\nrm{x-y}_2$, and let
$\Cl^{0,0}_M$ denote the class of convex $M$-Lipschitz functions. A problem
class consists of \emph{pairs} ``function $+$ initial point'':
\begin{equation}\label{eq:class}
   \Mfr_{L,R}=\bigl\{(f,x_0):\ f\in\Cl^{1,1}_L,\
   X^{*}\ne\varnothing,\
   \dist(x_0,X^{*})\le R\bigr\},\qquad
   X^{*}=\Argmin_{x\in\R^n}f(x),
\end{equation}
and similarly for the other classes. For a method $\mathcal M$ generating a
sequence $x_0,x_1,\dots$ through oracle calls, define
\begin{equation}\label{eq:complexity}
   N(\mathcal M,\Mfr_{L,R},\eps)
   =\sup_{(f,x_0)\in\Mfr_{L,R}}\
   \inf\bigl\{k:\ f(x_k)-f^{*}\le\eps\bigr\},
\end{equation}
where $f^{*}=\min_{x\in\R^n}f(x)$. The \emph{oracle complexity of the class}
is $\inf_{\mathcal M}N(\mathcal M,\Mfr_{L,R},\eps)$, where the infimum is
taken over methods of the specified type. For a \emph{randomized} method,
the criterion must be stated separately: either
$\Exp[f(x_k)-f^{*}]\le\eps$ or
$\Prb\{f(x_k)-f^{*}\le\eps\}\ge1-\beta$; these formulations are not
equivalent, and in bounds of the form~\eqref{eq:complexity} we always
specify which one is intended.

We also adopt three conventions used throughout.
\emph{All logarithms have dimensionless arguments}: instead of $\ln(1/\eps)$
we write $\ln\bigl(LR^2/\eps\bigr)$ or
$\ln\bigl((f(x_0)-f^{*})/\eps\bigr)$. The symbol $\Om(\cdot)$ always
denotes a \emph{lower} bound, $\Oh(\cdot)$ an upper bound, and
$\Ot(\cdot)$ an upper bound up to logarithmic factors.
Finally, $\range B$ is the range of the linear operator $B$,
$\lambda_{\min}(\cdot)$ and $\lambda_{\max}(\cdot)$ are the extreme
eigenvalues of a symmetric matrix, and $\lambda^{+}_{\min}(\cdot)$ is its
smallest positive eigenvalue.

\subsection{The simplest upper bound}

Let $Q=\R^n$ and $f\in\Cl^{1,1}_L$. Gradient descent
\begin{equation}\label{eq:gd}
   x_{k+1}=x_k-\tfrac1L\,\grad f(x_k)
\end{equation}
gives~\cite[Ch.~2]{NesBook2018}
\[
   f(x_k)-f^{*}\ \le\ \frac{2L\nrm{x_0-x^{*}}^2}{k+4},
   \qquad\text{that is,}\quad
   N(\eps)=\Oh\!\left(\frac{LR^2}{\eps}\right),
\]
and for $\mu$-strongly convex functions,
$N(\eps)=\Oh\bigl(\tfrac L\mu\ln\tfrac{LR^2}{\eps}\bigr)$. As in
classical algorithmic complexity theory, the bound is expressed in
$\Oh(\cdot)$ notation: it gives not a specific number, but the
characteristic dependence on the desired accuracy $\eps$ and the problem
parameters. A natural question is: can one do better?

\subsection{A lower bound and a ``universally bad function''}

The answer depends on how broad a class of methods we allow.
We distinguish two levels of generality.

\smallskip
\textbf{(A) Linear-span methods (span methods).} Suppose the method
constructs points in the linear span of the gradients obtained:
\begin{equation}\label{eq:span}
   x_k\in x_0+\Lin\bigl\{\grad f(x_0),\dots,\grad f(x_{k-1})\bigr\}.
\end{equation}
Condition~\eqref{eq:span} is satisfied by gradient descent, the conjugate
gradient method, and all the accelerated schemes in \S\,\ref{sec:fgm}.
For this class, the lower bound is proved using a single explicitly given
function.

Fix the number $N\ge1$ of oracle calls to be performed and the \emph{chain
length} $m$, setting $m=2N+1\le n$ (the latter is precisely the dimension
condition $N\le\frac12(n-1)$). Set
\begin{equation}\label{eq:worst}
   f_m(x)=\frac{L}{4}\left\{\frac12\Bigl[(x^{(1)})^2
   +\sum_{i=1}^{m-1}\bigl(x^{(i)}-x^{(i+1)}\bigr)^2+(x^{(m)})^2\Bigr]
   -x^{(1)}\right\},
\end{equation}
where $x^{(i)}$ are the coordinates of the vector $x$.
Function~\eqref{eq:worst} is convex and quadratic, and its gradient is
Lipschitz continuous with constant at most $L$. Its key property is its
\emph{chain structure}: the matrix of the quadratic form is tridiagonal,
and hence
\begin{equation}\label{eq:chain}
   x^{(j+1)}=\dots=x^{(n)}=0
   \ \Longrightarrow\
   \bigl(\grad f_m(x)\bigr)^{(j+2)}=\dots=\bigl(\grad f_m(x)\bigr)^{(n)}=0 .
\end{equation}
In other words, \emph{one oracle call advances the support by exactly one
coordinate}: starting from $x_0=0$, after $N$ calls a span method has a point
$x_N$ supported on the first $N$ coordinates.

Meanwhile, the minimizer of $f_m$ can be found explicitly: writing down
the stationarity conditions gives
\begin{equation}\label{eq:worstmin}
   (x_m^{*})^{(i)}=\frac{m+1-i}{m+1},\quad i=1,\dots,m;
   \qquad (x_m^{*})^{(i)}=0,\quad i>m .
\end{equation}
It is essential that $m=2N+1>N$: the solution also has nonzero coordinates
numbered $N+1,\dots,m$, so \emph{more than half} of the minimizer's
coordinates are inaccessible to the method. This yields the lower bound:
the residual $f_m(x_N)-f_m^{*}$ is no smaller than the contribution of these
``invisible'' coordinates. We emphasize that it is precisely the growth of
the chain length with the number of iterations that makes the construction
work: for a fixed chain of length $N$, the conjugate gradient method would
terminate in exactly $N$ steps.

\begin{theorem}[{\cite{NemYud1979};~\cite[Theorem 2.1.7]{NesBook2018}}]
\label{th:lower}
For every $k$, $1\le k\le\frac12(n-1)$, and every method
satisfying~\eqref{eq:span}, there exist a function $f\in\Cl^{1,1}_L$ and
an initial point $x_0$ such that
\begin{equation}\label{eq:lb1}
   f(x_k)-f^{*}\ \ge\ \frac{3\,L\,\nrm{x_0-x^{*}}^2}{32\,(k+1)^2} .
\end{equation}
Consequently, for class~\eqref{eq:class},
$N(\eps)=\Om\bigl(\sqrt{LR^2/\eps}\bigr)$. For the class of
$\mu$-strongly convex functions with an $L$-Lipschitz gradient (with
$n=\infty$), similarly,
$N(\eps)=\Om\bigl(\sqrt{L/\mu}\,\ln\frac{LR^2}{\eps}\bigr)$.
\end{theorem}

\smallskip
\textbf{(B) Arbitrary deterministic methods.} Condition~\eqref{eq:span}
does not hold for these methods, and instead of a single fixed function
one uses a \emph{resisting oracle}: the function is constructed during the
interaction so as to be consistent with all responses already given and
remain ``bad'' for the next step. The order of the bound remains the same,
but a dimension condition is required (the number of iterations is small
compared with $n$), and the constants change. It is specifically case~(A)
to which the constant $3/32$ and the restriction $k\le\frac12(n-1)$ in
Theorem~\ref{th:lower} apply. Randomized methods require a third
formulation, with a bound in expectation.

\subsection{Construction~\eqref{eq:worst} as a tool for obtaining lower bounds}

Function~\eqref{eq:worst} was devised by Yu.\,E.~Nesterov and has long
ceased to be a technical device for a single theorem. The reasoning
``exhibit a problem with a chain structure in which information propagates
by one coordinate per oracle call'' has become a \emph{standard tool} for
obtaining lower bounds: it is used for high-order methods, coordinate and
gradient-free schemes, distributed optimization, saddle-point problems,
and variational inequalities. The type of chain changes (one matrix,
two matrices, a block-bipartite graph), as does the type of oracle calls
being counted, but the idea remains the same.

Here is a recent example of interest because oracle calls of
\emph{different} costs are counted. The paper~\cite{BorGasKov2025}
considers a saddle-point problem with bilinear coupling,
\begin{equation}\label{eq:spp}
   \min_{x\in\R^{d_x}}\max_{y\in\R^{d_y}}\
   \bigl[F(x,y)=f(x)+\dotp{y}{Bx}-g(y)\bigr],
\end{equation}
where $f$ is $\mu_x$-strongly convex and $L_x$-smooth, $g$ is
$\mu_y$-strongly convex and $L_y$-smooth (the strong convexity constants
may be zero), and $B$ is the coupling matrix. Here the method accesses
three distinct sources of information whose access costs are generally
different: evaluation of $\grad f$, evaluation of $\grad g$, and
multiplication by $B$ and $B^{\!\top}$. Denote the time for one such call
by $\tau_f$, $\tau_g$, and $\tau_B$. It turns out that lower bounds on all
three counts can be obtained \emph{separately}, each expressed through
its own condition number.

We give a simplified statement; the full set of assumptions is given
in~\cite{BorGasKov2025}. Let $L_{xy}^2\ge\lambda_{\max}(B^{\!\top}B)$,
and let the numbers $\mu_{xy},\mu_{yx}\ge0$ satisfy
\begin{equation}\label{eq:rangecond}
   \mu_{xy}^2\le
   \begin{cases}
     \lambda^{+}_{\min}(B^{\!\top}B), & \text{if }
       \grad f(x)\in\range B^{\!\top}\ \ \forall x,\\
     \lambda_{\min}(B^{\!\top}B) & \text{otherwise;}
   \end{cases}
   \end{equation}
\begin{equation}\label{eq:rangecond2}
   \mu_{yx}^2\le
   \begin{cases}
     \lambda^{+}_{\min}(BB^{\!\top}), & \text{if }
       \grad g(y)\in\range B\ \ \forall y,\\
     \lambda_{\min}(BB^{\!\top}) & \text{otherwise}
   \end{cases}
\end{equation}
(thus $\mu_{xy}$ and $\mu_{yx}$ are generally distinct and do not reduce to
the same smallest positive singular value of $B$). Set
\begin{equation}\label{eq:kappas}
   \delta_x=\mu_x+\frac{\mu_{xy}^2}{L_y},\quad
   \delta_y=\mu_y+\frac{\mu_{yx}^2}{L_x},\qquad
   \kappa_x=\frac{L_x}{\delta_x},\quad
   \kappa_y=\frac{L_y}{\delta_y},\quad
   \kappa_{xy}=\frac{L_{xy}^2}{\delta_x\delta_y}.
\end{equation}
In the model of methods under consideration and under the full set of
structural assumptions in~\cite{BorGasKov2025}, the condition
$\min\{\delta_x,\delta_y\}>0$ is necessary and sufficient for linear
convergence. The quality of an approximate solution is measured by the
weighted distance to the solution sets $S_x\subseteq\R^{d_x}$ and
$S_y\subseteq\R^{d_y}$ of problem~\eqref{eq:spp} in the corresponding
groups of variables,
\begin{equation}\label{eq:weighteddist}
   R^2_{\delta_x\delta_y}(x,y)
   =\delta_x\dist^2(x;S_x)+\delta_y\dist^2(y;S_y),
\end{equation}
and the total running time of any linear-span method then satisfies the
lower bound
\begin{equation}\label{eq:lbspp}
   \tau\ \ge\ \Om\!\left(
   \tau_f\sqrt{\kappa_x}\,\ln\frac{cR^2}{\eps}
   +\tau_g\sqrt{\kappa_y}\,\ln\frac{cR^2}{\eps}
   +\tau_B\sqrt{\kappa_{xy}}\,\ln\frac{cR^2}{\eps}\right),
\end{equation}
where $c>0$ is a constant depending on $\kappa_x$ and $\kappa_y$.
The worst-case examples in~\cite{BorGasKov2025} are built precisely from
blocks of the form~\eqref{eq:worst}: a chain-structured quadratic form
for the ``functional'' part and a special sparse matrix $B$ for the
``coupling'' part. Bound~\eqref{eq:lbspp} unifies all previously known
special cases (strongly convex--strongly concave problems, bilinear games,
and minimization under affine constraints) and, importantly, is
\emph{attained} (\S\,\ref{sec:fgm}). This formula should be read primarily
as an illustration of \emph{complexity separation}: if information sources
differ in cost, they should also be assessed separately.

\subsection{The objective-function model}
\label{sec:models}

We now return to the main theme of \S\,\ref{sec:lookinside}. The lower
bound in Theorem~\ref{th:lower} is proved for methods operating
\emph{in the black-box model}. If a method knows not only the numbers
returned by the oracle but also how those numbers are structured, the bound
does not apply to it. The question arises: what kind of additional
knowledge can be useful?

The answer given in~\cite{NesDiss2013} consists of a classification of
\emph{objective-function models}. It is useful to distinguish three types.

\smallskip
\textbf{1. Geometric models.} One establishes that a simple convex set
containing the current point is included in a more complicated object:
the feasible set or the epigraph of the objective function. This makes it
possible to move in the desired direction while preserving feasibility
and/or guaranteeing a decrease in the function value. Examples include
the Dikin ellipsoid contained in the feasible region
(\S\,\ref{sec:ipm}); the quadratic upper model
$f(y)\le f(x)+\dotp{\grad f(x)}{y-x}+\frac L2\nrm{y-x}^2$
(\S\,\ref{sec:fgm}); a cubic upper model yielding Newton's method with
global guarantees (\S\,\ref{sec:cubic}); an upper model for the residual
norm of a system of equations (\S\,\ref{sec:gn}); and a regularized
Taylor model of arbitrary order (\S\,\ref{sec:tensor}).

\smallskip
\textbf{2. Algebraic models.} A particular representation of the function
is specified. For convex functions, a natural form is a maximum:
\begin{equation}\label{eq:algmodel}
   f(x)=\max_{u\in U}\ \bigl\{\dotp{Ax}{u}-\hat\phi(u)\bigr\}.
\end{equation}
Such a model makes it possible to construct approximate dual solutions
(\S\,\ref{sec:pd}) and, under favorable circumstances, to apply smoothing
techniques, improving the order of the convergence rate
(\S\,\ref{sec:smooth}).

\smallskip
\textbf{3. Structural models.} The structure of the function is described
by identifying two interacting objects, a ``simple'' one and a
``complicated'' one. The simple object may have poor functional
properties (it may be nonsmooth or even discontinuous), but it must allow
auxiliary problems involving it to be solved explicitly. An example is a
\emph{composite} function $f(x)+\Psi(x)$, where $f$ is smooth and $\Psi$
is a simple general convex function, such as the indicator of a set or
$\nrm{x}_1$ (\S\,\ref{sec:composite}). Relative smoothness
(\S\,\ref{sec:superfast}) and generalized $(L_0,L_1)$-smoothness
(\S\,\ref{sec:l0l1}) also belong here.

\smallskip
Table~\ref{tab:models} brings this map together. It should be read as a
table of contents for the second half of the article: almost every
subsequent section answers the question ``what does knowledge of a model
of this type give a method?''

\begin{table}[htbp]
\centering
\caption{Objective-function models according to~\cite{NesDiss2013} and the
sections of the article in which they are used.}
\label{tab:models}
\small
\renewcommand{\arraystretch}{1.25}
\begin{tabular}{@{}p{2.8cm}p{3.85cm}p{4.9cm}p{1.8cm}@{}}
\toprule
Model type & What is described & Examples & \S\\
\midrule
geometric
 & inclusion of a simple set in a more complicated object; an explicit
   upper model
 & Dikin ellipsoid; quadratic, cubic, and tensor upper models;
   an upper model for the residual norm
 & \S\S\,\ref{sec:ipm}, \ref{sec:fgm},
   \ref{sec:cubic}, \ref{sec:gn},
   \ref{sec:tensor}\\
algebraic
 & a representation of the function, most often~\eqref{eq:algmodel}
 & smoothing techniques; primal-dual and barrier subgradient methods
 & \S\S\,\ref{sec:smooth}, \ref{sec:pd}\\
structural
 & interaction between ``simple'' and ``complicated'' objects
 & composite functions $f+\Psi$; relative smoothness;
   $(L_0,L_1)$-smoothness
 & \S\S\,\ref{sec:composite},
   \ref{sec:l0l1},
   \ref{sec:superfast}\\
\bottomrule
\end{tabular}
\end{table}

\subsection{A second scale: cutting-plane methods}
\label{sec:cutting}

Oracle complexity depends not only on the function class but also on
\emph{which} local information is available. If a \emph{separation oracle}
is available instead of a gradient (given a point $x$, it either reports
that $x\in Q$ and returns a subgradient, or returns a hyperplane separating
$x$ from $Q$), then for the class of convex $M$-Lipschitz functions on a set
of diameter $R$ there exist methods whose complexity is \emph{polynomial
in the dimension and logarithmic in the accuracy}. Several levels must be
distinguished:

\begin{itemize}[leftmargin=1.6em,itemsep=3pt]
\item the \emph{center-of-gravity method}: $\Oh\bigl(n\ln(MR/\eps)\bigr)$
      oracle calls, which is optimal in order; the method, however,
      is not implementable, since it requires computing the center of
      gravity of a polyhedron;
\item the \emph{ellipsoid method} of Yudin--Nemirovski and Shor:
      $\Oh\bigl(n^2\ln(MR/\eps)\bigr)$ calls, since in one step the
      volume of the enclosing ellipsoid decreases only by a factor of
      $e^{1/(2n)}$; the method, however, is fully implementable. It was this
      method that enabled L.\,G.~Khachiyan~\cite{Khachiyan1979} to prove
      the polynomial-time solvability of linear programming in 1979;
\item \emph{Vaidya's method}~\cite{Vaidya1996}, using the volumetric
      barrier of a polyhedron: $\Oh\bigl(n\ln(\cdot)\bigr)$ oracle
      calls at polynomial iteration cost---an implementable method with
      an order-optimal number of calls;
\item modern schemes of Lee--Sidford--Wong~\cite{LeeSidfordWong2015}
      with improved arithmetic cost.
\end{itemize}

\begin{remark}[primal-dual certificates for a broad class of methods]
\label{rem:pdeverywhere}
Cutting-plane methods have another feature of interest in the context of
\S\,\ref{sec:pd}. First, they are robust to subgradient inexactness:
the error does not accumulate, as established for Vaidya's method
in~\cite{GladinEtAl2022}. Second, viewed appropriately, they turn out to
be \emph{primal-dual}: during the computation, the accumulated cuts allow
one to recover an approximate solution of the dual problem and a computable
gap bound. For Vaidya's method and related schemes with an inexact oracle,
the corresponding accuracy certificates were constructed
in~\cite{GladinGasnikovDvurechensky2025}. The observation is general:
for a broad class of convex optimization methods, one can construct a
primal-dual certificate by choosing the averaged objects and weights
appropriately; for the general mechanism, see \S\,\ref{sec:pd}
and~\cite{NemOnnRothblum2010}.
\end{remark}

Comparing the two scales, $\eps^{-1/2}$ and $\ln(1/\eps)$, is instructive.
The logarithmic scale has a substantially better dependence on accuracy,
but a worse dependence on dimension. For problems with
$n\sim10^6$--$10^9$, typical of data analysis, the choice is clear; for
problems with $n\sim10^2$ requiring high accuracy, it is equally clear,
but the opposite. Much of the story that follows is the story of attempts
to combine the advantages of both scales; interior-point methods
(\S\,\ref{sec:ipm}) are the most successful example.

\subsection{A summary map of bounds}

Table~\ref{tab:map} collects the lower bounds to be discussed. It should
be read carefully: each row is a \emph{triple} ``function class $+$ oracle
$+$ type of guarantee'', and bounds in different rows are not directly
comparable.

\begin{table}[p]
\centering
\caption{Bounds on the number of oracle calls for convex problems.
Block A shows limitations for all methods in the given class;
block B gives guarantees for the specified algorithms.}
\label{tab:map}
\small
\setlength{\tabcolsep}{4pt}
\renewcommand{\arraystretch}{1.15}
\begin{tabular}{@{}>{\raggedright\arraybackslash}p{4.8cm}
  >{\raggedright\arraybackslash}p{3.7cm}
  >{\raggedright\arraybackslash}p{5.85cm}@{}}
\toprule
Function class and oracle & Bound on the number of calls & Method and section\\
\midrule
\multicolumn{3}{@{}l@{}}{\textbf{A. Lower bounds; the methods attain the stated order}}\\[5pt]
$M$-Lipschitz.
\newline Oracle: subgradient.
 & $\displaystyle\Om\!\left(\frac{M^2R^2}{\eps^2}\right)$
 & Subgradient descent.
   \newline \S\,\ref{sec:complexity}.\\[4pt]
$L$-smooth.
\newline Oracle: gradient.
 & $\displaystyle\Om\!\left(\sqrt{\frac{LR^2}{\eps}}\right)$
 & Fast gradient method.
   \newline \S\,\ref{sec:fgm}.\\[4pt]
$\mu$-strongly convex,
$L$-smooth; $0<\mu<L$.
\newline Oracle: gradient.
 & $\displaystyle\Om\!\left(\sqrt{\frac L\mu}\,\ln\frac{LR^2}{\eps}\right)$
 & Fast gradient method
   with restarts.
   \newline \S\,\ref{sec:fgm}.\\[4pt]
Lipschitz Hessian;
constant $L_2$.
\newline Oracle: gradient
and Hessian.
 & $\displaystyle\Om\!\left(\left(\frac{L_2R^3}{\eps}\right)^{2/7}\right)$
 & Optimal second-order methods.
   \newline \S\,\ref{sec:tensor}.\\[4pt]
Lipschitz $p$th derivative;
constant $L_p$.
\newline Oracle: derivatives
up to order $p$.
 & $\displaystyle\Om\!\left(\left(\frac{L_pR^{p+1}}{\eps}\right)^{\frac{2}{3p+1}}\right)$
 & Optimal tensor methods.
   \newline \S\,\ref{sec:tensor}.\\[4pt]
\midrule
\multicolumn{3}{@{}l@{}}{\textbf{B. Upper bounds for specific methods}}\\[5pt]
Lipschitz Hessian;
constant $L_2$.
\newline Oracle: gradient
and Hessian.
 & $\displaystyle\Oh\!\left(\left(\frac{L_2R^3}{\eps}\right)^{1/3}\right)$
 & Accelerated cubic regularization.
   \newline \S\,\ref{sec:cubic}.\\[4pt]
$L$-smooth.
\newline Oracle: exact function
values at two points.
 & $\displaystyle\Oh\!\left(n\sqrt{\frac{LR^2}{\eps}}\right)$
   \newline Guarantee in expectation.
 & Gaussian smoothing
   and acceleration.
   \newline \S\,\ref{sec:zo}, Theorem~\ref{th:zo}.\\[4pt]
$M$-Lipschitz on $Q$,
$D=\operatorname{diam}Q$.
\newline Oracle: subgradient
and separation from $Q$.
 & $\displaystyle\Oh\!\left(n\ln\frac{MD}{\eps}\right)$
 & Center-of-gravity method;
   implementable schemes based on
   the volumetric barrier are discussed in \S\,\ref{sec:cutting}.\\[3pt]
\bottomrule
\end{tabular}
\par\medskip
\begin{minipage}{\textwidth}
\footnotesize\raggedright
\textbf{Notation.}
$\eps$ is the required error $f(x)-f^*$;
$R$ is an upper bound on the distance from the initial point to a solution;
$n$ is the dimension. In the last row, $D$ denotes the diameter of $Q$.
$\Om$ and $\Oh$ denote lower and upper bounds, respectively,
up to constant factors.

\smallskip
\textbf{Conditions for interpretation.}
The lower bounds in block A apply in sufficiently large dimension;
the strongly convex row considers asymptotics in $\eps$ for fixed $L$ and
$\mu$. The function-value row requires $\Exp[f(x)-f^*]\le\eps$;
in the other rows the guarantee is deterministic. Precise conditions are
given in the indicated sections. The cost of a call and the algorithm's
internal computations are not included here.
\end{minipage}
\end{table}

Two rows of the table deserve special comment. The row for second-order
methods appears twice: the lower bound has exponent $2/7$, but
\emph{accelerated cubic regularization} in \S\,\ref{sec:cubic}, despite
its name, gives exponent $1/3$; the optimal exponent is attained by another,
substantially more complicated construction (\S\,\ref{sec:tensor}). The
row for the gradient-free oracle refers to a specific two-point minimax
model and a guarantee in expectation; over a moderate range of accuracies,
the stated order cannot be improved, whereas under different assumptions
on noise, the number of points per query, and the geometry of the set,
the dimensional dependence changes (\S\,\ref{sec:zo}).

Finally, oracle complexity is not the only measure. The total
\emph{arithmetic} complexity equals the number of oracle calls times the
cost of a call, and the latter may depend dramatically on the
implementation: see \S\,\ref{sec:ad} (why the gradient is cheap) and
\S\,\ref{sec:cd} (how to make an iteration $n$ times cheaper).

\subsection{Summary of notation}
\label{sec:notation}

The discussion that follows introduces quite a few constants and radii,
with different normalizations natural in different sections.
Table~\ref{tab:notation} gathers them in one place; it is convenient to
refer back to it when reading \S\S\,\ref{sec:fgm}--\ref{sec:superfast}.

\begin{table}[htbp]
\centering
\caption{Summary of the main notation. The same letter may denote different
quantities in different sections if this is customary in the relevant
literature; such cases are noted separately.}
\label{tab:notation}
\footnotesize
\renewcommand{\arraystretch}{1.2}
\begin{tabular}{@{}p{2.15cm}p{9.5cm}p{1.9cm}@{}}
\toprule
Notation & Meaning & Introduced in\\
\midrule
$\Cl^{1,1}_L$, $\Cl^{0,0}_M$
 & class of convex functions with an $L$-Lipschitz gradient;
   class of convex $M$-Lipschitz functions
 & \S\,\ref{sec:complexity}\\

$\Mfr_{L,R}$
 & class of \emph{pairs} ``function $+$ initial point'' with normalization
   $\dist(x_0,X^{*})\le R$, where $X^{*}=\Argmin_{x\in\R^n}f(x)$
 & \eqref{eq:class}\\

$R$
 & upper bound on the distance from the initial point to the solution set:
   $\dist(x_0,X^{*})\le R$; in \S\,\ref{sec:smooth},
   the quantity $\max_{x\in Q}d_Q(x)$
 & \eqref{eq:class}\\

$R_{\mathrm{lev}}$, $D$
 & radius of the initial sublevel set
   $\sup\{\dist(x,X^{*}):f(x)\le f(x_0)\}$; needed where the method
   is monotone but the initial distance does not control the trajectory
 & \S\S\,\ref{sec:cubic}, \ref{sec:cd}\\

$\Delta_0$
 & initial function-value gap $f(x_0)-\inf_{x\in\R^n}f(x)$
 & \S\,\ref{sec:cubic}\\

$L$, $\mu$
 & Lipschitz constant of the gradient; strong convexity parameter;
   $\kappa=L/\mu$ is the condition number
 & \S\S\,\ref{sec:complexity}, \ref{sec:fgm}\\

$L_p$
 & Lipschitz constant of the $p$th derivative; $L_2$ for the Hessian,
   $L_3$ for the third derivative
 & \S\S\,\ref{sec:cubic}, \ref{sec:tensor}\\

$L_\nu$, $\nu$
 & H\"older constant and exponent of the gradient, $\nu\in[0,1]$
 & \eqref{eq:holder}\\

$L_0$, $L_1$
 & parameters of generalized smoothness
   $\nrm{\grad^2f}\le L_0+L_1\nrm{\grad f}$
 & \eqref{eq:l0l1}\\

$M_p$, $M$
 & regularization parameter of the upper model
   (cubic, tensor, prox-linear)
 & \S\S\,\ref{sec:cubic}, \ref{sec:gn}, \ref{sec:tensor}\\

$d_Q$, $d_U$; $\sigma_Q$, $\sigma_U$; $D_Q$, $D_U$
 & prox-functions on the primal and dual sets;
   their strong convexity moduli; their maxima
 & \S\,\ref{sec:smooth}\\

$\delta_x$, $\delta_y$; $\kappa_x$, $\kappa_y$, $\kappa_{xy}$
 & effective strong convexity parameters and three condition numbers
   of the saddle-point problem with bilinear coupling
 & \eqref{eq:kappas}\\

$\tau_f$, $\tau_g$, $\tau_B$
 & time for one call to each of the three information sources
   in the saddle-point problem
 & \S\,\ref{sec:complexity}\\

$\Lambda$, $R_\lambda$
 & bounded set of dual variables and its radius
 & \S\,\ref{sec:pd}\\

$\eta_e$
 & waiting time in the queue before arc $e$
   (dual variable of the traffic problem)
 & \S\,\ref{sec:traffic}\\

$\tau$
 & gradient-domination parameter (\S\,\ref{sec:cubic});
   in \S\,\ref{sec:smooth}, the step size of the excessive-gap scheme
 & \eqref{eq:pl}\\
\bottomrule
\end{tabular}
\end{table}

A separate point concerns the \emph{type of guarantee}. Unless stated
otherwise, the bound is deterministic and applies to the last point
constructed. For randomized schemes (\S\S\,\ref{sec:cd},~\ref{sec:zo}),
the guarantee is in expectation, as specified in the theorem statement;
averaged points are explicitly identified (as in~\eqref{eq:primalrate}
and~\eqref{eq:dualrate}).

%% file: en/r3_fgm.tex

\section{Fast gradient method (1983)}
\label{sec:fgm}

\subsection{Background: the heavy-ball method}

By the early 1980s, the situation was as follows. The lower bound
$\Om\bigl(\sqrt{LR^2/\eps}\bigr)$ for the class $\Cl^{1,1}_L$ was known
(Theorem~\ref{th:lower}), gradient descent gave
$\Oh\bigl(LR^2/\eps\bigr)$, and the gap between them was a square root.
Did a method attaining the lower bound exist?

There was a candidate. In 1964, B.\,T.~Polyak proposed the
\emph{heavy-ball method}~\cite{Polyak1964},
\begin{equation}\label{eq:hb}
   x_{k+1}=x_k-\alpha\,\grad f(x_k)+\beta\,(x_k-x_{k-1}),
\end{equation}
motivated by a mechanical analogy: a point with inertia rolling down a
``ravine'' does not bounce off its walls as gradient descent does, but
accumulates motion along the bottom. For quadratic functions, an appropriate
choice of $\alpha,\beta$ indeed gives the rate
$\Oh\bigl(\sqrt{L/\mu}\,\ln\frac1\eps\bigr)$---exactly what is needed.
For general smooth convex functions, however, scheme~\eqref{eq:hb} may
fail to converge: there are known examples of strongly convex functions
with Lipschitz gradients on which the heavy-ball method cycles with its
``quadratic'' parameters. Mechanical intuition is misleading here, and the
remedy has two parts. First, for the general convex class, the extrapolation
coefficient must depend on the iteration number (in the strongly convex
case, a constant coefficient is still possible; see~\eqref{eq:fgmrate2}).
Second, and no less importantly, Nesterov's scheme differs from the
heavy-ball method in \emph{how the gradient step and extrapolation are
coupled}: the gradient is evaluated not at the current point, but at an
auxiliary point obtained by extrapolation.

This is precisely what was done in~\cite{Nes1983}.

\subsection{The fast gradient method}

Let $f$ be convex and differentiable on $\R^n$, with an $L$-Lipschitz
gradient. Take $x_0=y_0$, $a_0=1$, and set
\begin{equation}\label{eq:fgm}
   \left\{
   \begin{aligned}
   x_{k+1}&=y_k-\tfrac1L\,\grad f(y_k),\\
   a_{k+1}&=\tfrac12\bigl(1+\sqrt{1+4a_k^2}\bigr),\\
   y_{k+1}&=x_{k+1}+\tfrac{a_k-1}{a_{k+1}}\,(x_{k+1}-x_k).
   \end{aligned}\right.
\end{equation}
The first line is an ordinary gradient step, but taken from the
\emph{auxiliary} point $y_k$; the third is extrapolation with coefficient
$\frac{a_k-1}{a_{k+1}}$, which tends to one as
$1-3/k+\Oh(k^{-2})$.

\begin{theorem}[\cite{Nes1983,NemNes1985}]\label{th:fgm}
\emph{Class:} $f\in\Cl^{1,1}_L$, $Q=\R^n$, $\Argmin_{x\in\R^n}f(x)\ne\varnothing$.
\emph{Oracle:} first order (value and gradient), one call per iteration.
\emph{Normalization:} $R=\nrm{x_0-x^{*}}_2$. \emph{Guarantee:}
deterministic, for the last point. Then method~\eqref{eq:fgm} satisfies
\begin{equation}\label{eq:fgmrate}
   f(x_k)-f^{*}\ \le\ \frac{2LR^2}{(k+1)^2},
   \qquad\text{that is,}\qquad
   N(\eps)=\Oh\!\left(\sqrt{\frac{LR^2}{\eps}}\right).
\end{equation}
If, in addition, $f$ is $\mu$-strongly convex, then the variant with the
constant extrapolation coefficient
$\beta=\frac{\sqrt{L/\mu}-1}{\sqrt{L/\mu}+1}$ gives
\begin{equation}\label{eq:fgmrate2}
   f(x_k)-f^{*}\ \le\ \bigl(f(x_0)-f^{*}+\tfrac\mu2 R^2\bigr)
   \Bigl(1-\sqrt{\tfrac\mu L}\Bigr)^{k},
   \qquad
   N(\eps)=\Oh\!\left(\sqrt{\tfrac L\mu}\,
   \ln\frac{LR^2}{\eps}\right).
\end{equation}
\end{theorem}

A terminological remark: in the Russian literature,
scheme~\eqref{eq:fgm} is usually called the \emph{fast} gradient method,
whereas in the English-language literature it is called
\emph{accelerated}; below we use the two words synonymously.

Comparing~\eqref{eq:fgmrate} with Theorem~\ref{th:lower}, we see that the
bounds match in order. \emph{The method is optimal on its class in its
oracle model}---and this is the only section of the article devoted to a
result in which the information model \emph{does not} change. Note that
optimality applies specifically to the pair ``class $+$ oracle'': on
quadratic functions, the conjugate gradient method is no worse, while
access to a separation oracle changes the picture entirely
(\S\,\ref{sec:cutting}).

\subsection{Where acceleration comes from: primal and dual methods}
\label{sec:primaldualgrad}

Scheme~\eqref{eq:fgm} looks as though it has been guessed. Yet there is a
presentation in which it arises as a \emph{combination of two elementary
methods}; this is how the material is organized in~\cite{NesDiss2013}.
We reproduce this construction because it gives a particularly transparent
interpretation of acceleration.

Let $Q$ be a closed convex set equipped with a \emph{prox-function}
$d(x)$: continuous, strongly convex on $Q$ with parameter~1,
differentiable, with center $x_0=\argmin_{x\in Q}d(x)$ and normalization
$d(x_0)=0$. In this section, the prox-center coincides with the method's
initial point, so the notation $x_0$ causes no ambiguity. Set
\begin{equation}\label{eq:bregman}
   \rho(x,y)=d(y)-d(x)-\dotp{\grad d(x)}{y-x},
   \qquad
   \mathcal B_Q(x,g)=\argmin_{y\in Q}
   \bigl\{\dotp{g}{y-x}+\rho(x,y)\bigr\}
\end{equation}
(the Bregman distance and Bregman projection), and let
\begin{equation}\label{eq:gradmap}
   T_Q(x)=\argmin_{y\in Q}
   \Bigl\{\dotp{\grad f(x)}{y-x}+\tfrac L2\nrm{y-x}^2\Bigr\}
\end{equation}
denote the \emph{gradient mapping}---a step based on a quadratic upper
model (the geometric model in Table~\ref{tab:models}).

\smallskip
The \emph{primal gradient method} consists of successive applications of
the Bregman projection:
\begin{equation}\label{eq:primalgm}
   x_{k+1}=\mathcal B_Q\bigl(x_k,\tfrac1L\grad f(x_k)\bigr),
   \qquad k\ge0 .
\end{equation}
For this method~\cite[Theorem 1.2.1]{NesDiss2013}, with
$\bar x_k=\frac1k\sum_{i=1}^kx_i$,
\begin{equation}\label{eq:primalrate}
   f(\bar x_k)-f^{*}\ \le\
   \frac1k\sum_{i=1}^{k}\bigl[f(x_i)-f^{*}\bigr]\ \le\
   \frac{L\,d(x^{*})}{k}.
\end{equation}

\smallskip
The \emph{dual gradient method} operates not on points but on a
\emph{model} of the objective function: it accumulates linearizations.
\begin{equation}\label{eq:dualgm}
   \psi_0(x)=\tfrac12\rho(x_0,x),\qquad
   x_k=\argmin_{x\in Q}\psi_k(x),
\end{equation}
\begin{equation}\label{eq:dualgm2}
   \psi_{k+1}(x)=\psi_k(x)+\tfrac1L
   \bigl[f(x_k)+\dotp{\grad f(x_k)}{x-x_k}\bigr],
   \qquad k\ge0 .
\end{equation}
For this method~\cite[Theorem 1.2.2]{NesDiss2013}, with $y_k=T_Q(x_k)$
and $\bar y_k=\frac1{k+1}\sum_{i=0}^ky_i$,
\begin{equation}\label{eq:dualrate}
   f(\bar y_k)-f^{*}\ \le\
   \frac1{k+1}\sum_{i=0}^{k}\bigl[f(y_i)-f^{*}\bigr]\ \le\
   \frac{L\,d(x^{*})}{k+1}.
\end{equation}

Both bounds are of order $L\,d(x^{*})/k$, so neither method alone yields
acceleration. But the methods work differently: the first moves a
\emph{point}, the second refines a \emph{model}, and the quantities they
accumulate ``work in opposite phases''. The key observation is that they
can be coupled.

\begin{algo}[the fast gradient method as a combination of two methods]
\label{algo:fgmcomb}
For $k\ge0$: evaluate $f(x_k)$ and $\grad f(x_k)$; set $y_k=T_Q(x_k)$;
find
\begin{equation}\label{eq:fgmcomb}
   z_k=\argmin_{x\in Q}\Bigl\{L\,d(x)+\sum_{i=0}^{k}\tfrac{i+1}{2}
   \bigl[f(x_i)+\dotp{\grad f(x_i)}{x-x_i}\bigr]\Bigr\};
   \qquad
   x_{k+1}=\tfrac{2}{k+3}\,z_k+\tfrac{k+1}{k+3}\,y_k .
\end{equation}
\end{algo}

\begin{theorem}[{\cite[Theorem 1.2.3]{NesDiss2013}}]\label{th:fgmcomb}
For scheme~\ref{algo:fgmcomb}, for all $k\ge0$,
\begin{equation}\label{eq:fgmcombest}
   \frac{(k+1)(k+2)}{4}\,f(y_k)\ \le\
   \min_{x\in Q}\Bigl\{L\,d(x)+\sum_{i=0}^{k}\tfrac{i+1}{2}
   \bigl[f(x_i)+\dotp{\grad f(x_i)}{x-x_i}\bigr]\Bigr\},
\end{equation}
and consequently,
\begin{equation}\label{eq:fgmcombrate}
   f(y_k)-f^{*}\ \le\ \frac{4L\,d(x^{*})}{(k+1)(k+2)} .
\end{equation}
\end{theorem}

The meaning of~\eqref{eq:fgmcombest} is transparent. On the right is a
\emph{lower model} of the problem: a sum of linearizations with weights
$\frac{i+1}{2}$ plus the regularizer $L\,d(x)$. On the left is the
function value at the current point multiplied by the total weight
$\sum_{i\le k}\frac{i+1}{2}\asymp k^2/4$. The inequality says that the
method keeps the function value below the accumulated model, while the
weights grow \emph{linearly}---hence $k^2$ in the denominator
of~\eqref{eq:fgmcombrate}. For $d(x)=\frac12\nrm{x-x_0}^2$,
formula~\eqref{eq:fgmcombrate} agrees with~\eqref{eq:fgmrate} up to a
constant factor.

\begin{remark}\label{rem:linearcoupling}
The idea of ``coupling a primal step with a dual step'' was later presented
independently under the name linear coupling~\cite{AllenZhuOrecchia2017},
and became widely known in that form. It is useful to remember, however,
that construction~\eqref{eq:fgmcomb} and its interpretation as a
combination of two elementary schemes belong to the method's author and
are presented in~\cite{NesDiss2013,NesBook2018}.
\end{remark}

\smallskip
\emph{Compatibility of norms.} Bound~\eqref{eq:fgmcombrate} relies on an
easily overlooked condition: the constant $L$ in~\eqref{eq:gradmap} and
the prox-function $d$ must refer to \emph{the same} norm $\nrm{\cdot}$.
The gradient must be $L$-Lipschitz with respect to this norm,
$\nrm{\grad f(x)-\grad f(y)}_{*}\le L\nrm{x-y}$, where
$\nrm{\cdot}_{*}$ is the dual norm, and $d$ must be $1$-strongly convex
with respect to the same norm; only then does the ``distance'' $d(x^{*})$
from the initial point to the solution enter the bound as
in~\eqref{eq:fgmcombrate}. The pair ``norm $+$ prox-function'' (the
\emph{prox-structure}) is chosen according to the geometry of $Q$.
For example, on the standard simplex $\Delta_n$, the Euclidean
prox-structure gives $\Oh(L/k^2)$, whereas the pair ``$\nrm{\cdot}_1$ $+$
entropy $d(x)=\ln n+\sum_ix^{(i)}\ln x^{(i)}$'', for which
$d(x^{*})\le\ln n$, gives $\Oh(\hat L\ln n/k^2)$, where $\hat L$ is the
Lipschitz constant of the gradient as a mapping from
$(\R^n,\nrm{\cdot}_1)$ to $(\R^n,\nrm{\cdot}_\infty)$.
Since $\hat L\le L\le n\hat L$, the gain can be as large as a factor of
$n/\ln n$.

What if the norms are \emph{not compatible}? Suppose a function whose
gradient is Lipschitz in the Euclidean norm with constant $L$ must be
minimized over the $\ell_1$-ball $\{\nrm{x}_1\le R\}$. Any
prox-structure forces a choice of one of the two norms: the Euclidean
one gives $\Oh(LR^2/N^2)$ after $N$ iterations (in the Euclidean norm,
the $\ell_1$-ball is ``as large'' as a Euclidean ball of the same radius),
while the entropic one gives $\Oh(\hat LR^2\ln n/N^2)$, and it is
possible that $\hat L=L$. Meanwhile, the construction of C.~Guzm\'an and
A.\,S.~Nemirovski~\cite{GuzmanNemirovski2015} yields only a lower bound
of order $LR^2/N^3$ (when $n\gtrsim N$), and the question of the correct
exponent---$2$ or $3$?---remained open for a long
time~\cite{Guzman2015,VorontsovaGasnikovGorbunov2019}. The answer is given
in~\cite{LoginovEtAl2026}: when $n\gtrsim N$, the minimax error after $N$
first-order oracle calls is of order $LR^2/N^3$ (up to a factor
$\ln N\ln n$ in the upper bound), for both deterministic and randomized
methods. The optimal method is not a prox-method but an accelerated
variant of the level method, in which the center of the current
localization set is its \emph{Steiner point}: the Gaussian average of the
points maximizing random linear functionals. The error of the accelerated
scheme is proportional to $L/N^2$ times the sum of squared displacements
of its ``anchor'' sequence. For prox-methods, this sum is of order $R^2$,
whereas the Steiner points of nested subsets of the $\ell_1$-ball move in
steps of length $\Ot(R/N)$ along a path of Euclidean length
$\Oh\bigl(R\sqrt{\ln N\ln n}\bigr)$, so that the sum of squared
displacements is of order $\Ot(R^2/N)$---hence the extra power of $N$.
Simultaneously and independently, D.~Mart\'inez-Rubio, C.~Guzm\'an, and
coauthors obtained a general answer for $\ell_p$-balls and smoothness
measured in the $\ell_q$-norm, where $1\le p,q\le2$ are arbitrary numbers
unrelated to one another~\cite{MartinezRubioEtAl2026b}, as well as for
the nonsmooth case: for $M$-Lipschitz functions, mismatched geometry gives
the rate $\Ot(MR/N)$ instead of the classical
$\Oh(MR/\sqrt N)$~\cite{MartinezRubioGuzman2026,MartinezRubioEtAl2026a};
the role of the Steiner point there is played by a specially constructed
``stable center''. Thus the rule ``one norm for $L$ and $d$'' is a
property of prox-methods, not of the problem itself: when the norms are
mismatched, the optimal method has a different structure.

\subsection{Estimate functions}
\label{sec:estimatefun}

Formally, Theorem~\ref{th:fgm} is conveniently proved using the technique
of \emph{estimate functions} (estimate sequences), which has since become
a standard tool.

\begin{definition}\label{def:est}
A pair of sequences $\{\varphi_k(x)\}_{k\ge0}$ (functions) and
$\{\vartheta_k\}_{k\ge0}$ (numbers, $\vartheta_k\to0$) is called an
\emph{estimate sequence} for $f$ if, for all $k$ and all $x\in\R^n$,
\begin{equation}\label{eq:est}
   \varphi_k(x)\ \le\ (1-\vartheta_k)f(x)+\vartheta_k\varphi_0(x).
\end{equation}
\end{definition}

\begin{lemma}\label{lem:est}
If, for some point $x_k$,
$f(x_k)\le\min_{x\in\R^n}\varphi_k(x)$, then
$f(x_k)-f^{*}\le\vartheta_k\bigl(\varphi_0(x^{*})-f^{*}\bigr)$.
\end{lemma}

\begin{proof}
Substitute $x=x^{*}$ into~\eqref{eq:est}:
$\min_{x\in\R^n}\varphi_k(x)\le\varphi_k(x^{*})\le(1-\vartheta_k)f^{*}
+\vartheta_k\varphi_0(x^{*})$, whence
$f(x_k)\le f^{*}+\vartheta_k(\varphi_0(x^{*})-f^{*})$.
\end{proof}

The task reduces to constructing an estimate sequence with rapidly
decreasing $\vartheta_k$ for which the condition
$f(x_k)\le\min_{x\in\R^n}\varphi_k(x)$ is maintained automatically.
Take $\varphi_0(x)=f(x_0)+\frac{\gamma_0}2\nrm{x-x_0}^2$ and recursively
\begin{equation}\label{eq:estrec}
   \varphi_{k+1}(x)=(1-\alpha_k)\varphi_k(x)
   +\alpha_k\bigl[f(y_k)+\dotp{\grad f(y_k)}{x-y_k}
   +\tfrac\mu2\nrm{x-y_k}^2\bigr];
\end{equation}
all $\varphi_k$ remain quadratic,
$\vartheta_k=\prod_{j<k}(1-\alpha_j)$, and choosing $\alpha_k$ from
$L\alpha_k^2=(1-\alpha_k)\gamma_k+\alpha_k\mu$ leads precisely
to~\eqref{eq:fgm}. For $\mu=0$, one obtains
$\vartheta_k=\Oh(k^{-2})$, hence~\eqref{eq:fgmrate}.

\subsection{Acceleration as gradient descent in moving coordinates}
\label{sec:accsimplesec}

We give another, very recent and very short
explanation~\cite[\S\,4.2]{BorGasKov2025}: \emph{an accelerated method is
ordinary gradient descent applied to a suitably reparameterized function}.

Let $h$ be convex and $L$-smooth. Define the sequence
\begin{equation}\label{eq:alpha}
   \alpha_0=1,\qquad
   \alpha_{t+1}=\frac{2}{1+\sqrt{1+4/\alpha_t^2}},\qquad t\ge0,
\end{equation}
for which $\alpha_t\in(0,1]$ and
\begin{equation}\label{eq:alpharec}
   \alpha_t^{-2}=\alpha_{t-1}^{-2}+\alpha_t^{-1},
   \qquad \alpha_t^{-1}\ge\tfrac12(t+2),
   \qquad \alpha_{-1}^{-2}:=0 .
\end{equation}
(The first equality: setting $v=\alpha_t^{-1}$,
$u=\alpha_{t-1}^{-2}$, the definition gives
$2v-1=\sqrt{1+4u}$, that is, $v^2-v=u$. The second:
$\alpha_t^{-1}=\frac12+\frac12\sqrt{1+4\alpha_{t-1}^{-2}}
\ge\frac12+\alpha_{t-1}^{-1}$ and $\alpha_0^{-1}=1$.)

Fix an ``anchor'' point $\bar z^{\,t}$ and consider the
\emph{reparameterized} function
\begin{equation}\label{eq:ht}
   h_t(z)=\alpha_t^{-2}\,
   h\bigl(\alpha_t z+(1-\alpha_t)\bar z^{\,t}\bigr).
\end{equation}
Direct differentiation gives
$\grad h_t(z)=\alpha_t^{-1}\grad h(\alpha_tz+(1-\alpha_t)\bar z^{\,t})$
and $\grad^2h_t(z)=\grad^2h(\alpha_tz+(1-\alpha_t)\bar z^{\,t})$; hence
\emph{$h_t$ is convex and $L$-smooth with the same constant $L$}---this is
the key observation.

\begin{algo}[the accelerated method as gradient descent]\label{algo:acc}
$z^0=\bar z^{\,0}$; for $t=0,1,2,\dots$,
\begin{equation}\label{eq:accsimple}
   z^{t+1}=z^{t}-\tfrac1L\,\grad h_t(z^{t}),
   \qquad
   \bar z^{\,t+1}=\alpha_t z^{t+1}+(1-\alpha_t)\bar z^{\,t}.
\end{equation}
\end{algo}

Writing $y^{t}=\alpha_t z^{t}+(1-\alpha_t)\bar z^{\,t}$ (the point where
the gradient of $h$ is actually evaluated), we rewrite~\eqref{eq:accsimple}
in the familiar three-sequence form
\begin{equation}\label{eq:accthree}
   y^{t}=\alpha_t z^{t}+(1-\alpha_t)\bar z^{\,t},\qquad
   z^{t+1}=z^{t}-\frac{1}{L\alpha_t}\grad h(y^{t}),\qquad
   \bar z^{\,t+1}=y^{t}-\frac1L\grad h(y^{t}),
\end{equation}
so $\bar z^{\,t+1}$ is an ordinary gradient step from $y^t$ (the ``primal''
sequence), while $z^{t+1}$ is a step with a multiplier increased by a
factor $1/\alpha_t$ (the ``dual'' sequence). The connection with
scheme~\ref{algo:fgmcomb} in \S\,\ref{sec:primaldualgrad} is immediately
apparent on comparing the formulas.

\begin{lemma}[one gradient-descent step]\label{lem:gdstep}
Let $\phi$ be convex and $L$-smooth, and
$z^{+}=z-\frac1L\grad\phi(z)$. Then, for all $u$,
\begin{equation}\label{eq:gdstep}
   \phi(z^{+})+\tfrac L2\nrm{z^{+}-u}^2
   \ \le\ \phi(u)+\tfrac L2\nrm{z-u}^2 .
\end{equation}
\end{lemma}

\begin{proof}
By $L$-smoothness,
$\phi(z^{+})\le\phi(z)-\frac1{2L}\nrm{\grad\phi(z)}^2$.
Expanding the square gives
$\frac L2\nrm{z^{+}-u}^2=\frac L2\nrm{z-u}^2-\dotp{\grad\phi(z)}{z-u}
+\frac1{2L}\nrm{\grad\phi(z)}^2$. Adding and using convexity,
$\phi(z)-\dotp{\grad\phi(z)}{z-u}\le\phi(u)$,
yields~\eqref{eq:gdstep}.
\end{proof}

\begin{theorem}[\cite{BorGasKov2025}]\label{th:accsimple}
For scheme~\ref{algo:acc} and any $u\in\R^n$, for all $t\ge0$,
\begin{equation}\label{eq:accrec}
   \tfrac L2\nrm{z^{t+1}-u}^2
   +\alpha_t^{-2}\bigl[h(\bar z^{\,t+1})-h(u)\bigr]
   \ \le\
   \tfrac L2\nrm{z^{t}-u}^2
   +\alpha_{t-1}^{-2}\bigl[h(\bar z^{\,t})-h(u)\bigr],
\end{equation}
and, consequently,
$h(\bar z^{\,N})-h^{*}\le\dfrac{2L\nrm{z^{0}-z^{*}}^2}{(N+1)^2}$.
\end{theorem}

\begin{proof}
Apply Lemma~\ref{lem:gdstep} to $\phi=h_t$ (it is $L$-smooth!) at
$z=z^{t}$:
$h_t(z^{t+1})+\frac L2\nrm{z^{t+1}-u}^2\le h_t(u)+\frac L2\nrm{z^{t}-u}^2$.
By definitions~\eqref{eq:ht} and~\eqref{eq:accsimple},
$h_t(z^{t+1})=\alpha_t^{-2}h(\bar z^{\,t+1})$; bound the right-hand side
using convexity of $h$:
\[
   h_t(u)=\alpha_t^{-2}h\bigl(\alpha_tu+(1-\alpha_t)\bar z^{\,t}\bigr)
   \le\alpha_t^{-1}h(u)+\bigl(\alpha_t^{-2}-\alpha_t^{-1}\bigr)
   h(\bar z^{\,t}).
\]
Subtracting $\alpha_t^{-2}h(u)$ from both sides and substituting
$\alpha_t^{-2}-\alpha_t^{-1}=\alpha_{t-1}^{-2}$
from~\eqref{eq:alpharec} gives~\eqref{eq:accrec}. Summing over
$t=0,\dots,N-1$ (for $t=0$, the right-hand side equals
$\frac L2\nrm{z^0-u}^2$, since $\alpha_{-1}^{-2}=0$) and setting
$u=z^{*}$ gives
$\alpha_{N-1}^{-2}[h(\bar z^{\,N})-h^{*}]\le\frac L2\nrm{z^0-z^{*}}^2$;
it remains to use $\alpha_{N-1}^{-1}\ge\frac12(N+1)$.
\end{proof}

Acceleration here arises not from a new step, but from taking the step at
each iteration on a \emph{different} function $h_t$, obtained from $h$ by
contraction toward the current anchor point; the factor $\alpha_t^{-2}$
in~\eqref{eq:ht} ensures accumulation of order $\Oh(t^2)$. Below
(\S\,\ref{sec:sliding}) we will see that this viewpoint is precisely what
allows acceleration to be extended to structured saddle-point problems.

\subsection{Two more perspectives: the continuous limit and control theory}

If in method~\eqref{eq:fgm} we set $t=k/\sqrt L$ and let $L\to\infty$,
the trajectory converges to the solution of the equation
\begin{equation}\label{eq:ode}
   \ddot X(t)+\frac3t\,\dot X(t)+\grad f\bigl(X(t)\bigr)=0
\end{equation}
(W.~Su, S.~Boyd, E.~Cand\`es~\cite{SuBoydCandes2016}). This describes the
motion of a point in the potential $f$ with friction whose coefficient
decreases as $3/t$. The decrease is essential: with constant friction
(the heavy-ball method~\eqref{eq:hb}), exponential convergence is obtained
only in the strongly convex case, whereas the law $3/t$ gives
$f(X(t))-f^{*}=\Oh(t^{-2})$ already for convex functions.
Equation~\eqref{eq:ode} also reveals the method's \emph{nonmonotonicity}:
the point ``overshoots'' the minimum and undergoes damped oscillations.

There is also a third explanation, through control theory. An iterative
method is viewed as a linear system with nonlinear feedback, and the
proof of the convergence rate as a search for a Lyapunov function; the
latter reduces to the feasibility of a small semidefinite inequality.
This approach was developed systematically by L.~Lessard, B.~Recht, and
A.~Packard~\cite{LessardRechtPackard2016}; for the accelerated method it
yields a short, computer-verifiable proof~\cite{HuLessard2017}.

\begin{remark}\label{rem:nonmono}
Nonmonotonicity is not a defect of the presentation but a property of the
scheme: the values $f(x_k)$ in~\eqref{eq:fgm} may increase. In practice,
this is addressed by \emph{restarts}: when an increase in the function
value is detected, the method is restarted with $a_k:=1$. Under additional
assumptions (for example, quadratic growth near the minimum), adaptive
restarting recovers the optimal linear rate~\eqref{eq:fgmrate2} without
knowing $\mu$~\cite{ODonoghueCandes2015}; there is no unconditional theorem
that ``restarting always fixes everything''. There is, however, also a
\emph{monotone} accelerated method; see \S\,\ref{sec:l0l1}.
\end{remark}

\subsection{What this looks like in an example}

\begin{figure}[htbp]
\centering
\includegraphics[width=\textwidth]{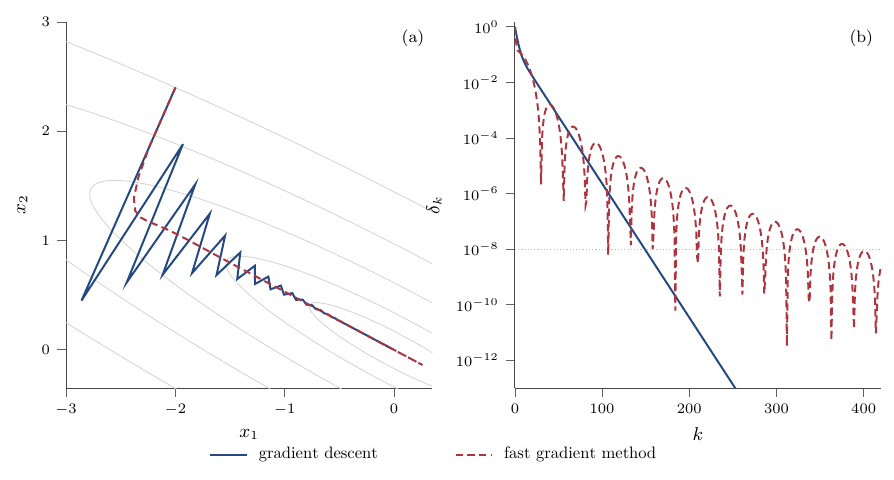}
\caption{Dynamics of two methods on a quadratic function with
$A=Q\operatorname{diag}(1,0{,}03)Q^\top$, where $Q$ specifies a rotation
by $62^\circ$, and initial point $x_0=(-2,\,2{,}4)$.
(a)~The first $100$ steps in $\R^2$; level curves are shown in gray.
(b)~The raw relative residual $\delta_k$; the gray dashed line marks
$10^{-8}$. The blue line is gradient descent with $h=1{,}8/L$;
the red dashed line is scheme~\eqref{eq:fgm} with varying $a_k$ and
$\widehat L=2L$. The spikes in the red curve are actual oscillations
of the residual, not its smoothed envelope.}
\label{fig:fgm}
\end{figure}

\begin{example}\label{ex:kappa}
Take
\[
   f(x)=\frac12\dotp{Ax}{x},\qquad x\in\R^2,
\]
where the eigenvalues of $A$ are $L=1$ and $\mu=0{,}03$,
so that $\kappa=L/\mu=100/3\approx33{,}3$. The eigenvector matrix is
rotated by $62^\circ$, and the initial point is $x_0=(-2,\,2{,}4)$.
The gray ellipses in Fig.~\ref{fig:fgm}\,(a) are level curves of the
function.

With the exact step size $1/L$, the component of gradient descent in the
eigendirection with $\lambda=L$ vanishes in one step. The trajectory
therefore immediately reaches the slow eigenaxis and looks nearly
straight; for a two-dimensional quadratic function, this is a correct but
special case for visualization. To show the characteristic transverse
motion in a narrow ravine, gradient descent in Fig.~\ref{fig:fgm}\,(a),~(b)
uses the stable step size $h=1{,}8/L$. In the stiff direction, the
one-step multiplier is $1-hL=-0{,}8$, so the blue trajectory forms a sharp,
damped ``sawtooth''.

The red trajectory is generated by scheme~\eqref{eq:fgm} with varying
$a_k$ and the valid conservative estimate $\widehat L=2L$, that is, with
gradient step size $1/\widehat L$. Accumulated inertia causes the method
to pass close to the minimum repeatedly and overshoot it. Thus the raw
residual in Fig.~\ref{fig:fgm}\,(b) is nonmonotone and has characteristic
spikes, consistent with the continuous-limit picture~\eqref{eq:ode}.
Here
\[
   \delta_k=\frac{f(x_k)-f^*}{f(x_0)-f^*}.
\]
Gradient descent falls below $10^{-8}$ and remains below it from $k=149$
onward. The fast gradient method first crosses this level at $k=107$,
but because of oscillations only stays below it permanently from $k=381$
onward. This experiment illustrates the shapes of the trajectories;
the dependence of complexity on $\kappa$ is compared separately below,
with both methods using the same step size.
\end{example}

To examine the dependence on the condition number, consider a separate
family of two-dimensional problems,
\[
   f_\kappa(x)=\frac12\left(x_1^2+\kappa^{-1}x_2^2\right),
   \qquad
   x_{0,\kappa}=\frac{(1,\sqrt\kappa)}{\sqrt{1+\kappa}}.
\]
This normalization gives $\nrm{x_{0,\kappa}}=1$ and equal initial
contributions of the two eigendirections to the function value. Thus, as
$\kappa$ varies, the slow mode does not disappear, and the problems being
compared have the same initial-error structure. Both methods use step size
$1/L=1$ and one gradient evaluation per iteration; the accelerated scheme
uses the constant coefficient
$\beta=(\sqrt\kappa-1)/(\sqrt\kappa+1)$ from~\eqref{eq:fgmrate2}.

\begin{figure}[htbp]
\centering
\includegraphics[width=0.72\textwidth]{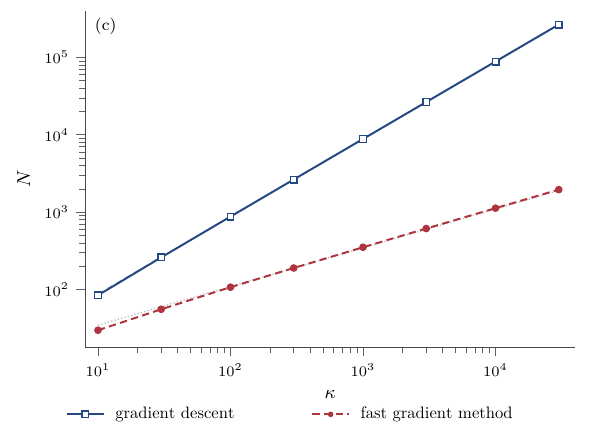}
\caption{Dependence of the number of gradient evaluations $N$ on the
condition number $\kappa$ for the family $f_\kappa$ and initial points
$x_{0,\kappa}$ specified in the text. The stopping criterion is
$\delta_k\le10^{-8}$. Both methods use step size $1/L$; the fast gradient
method uses the constant extrapolation coefficient from~\eqref{eq:fgmrate2}.
The gray dashed line is the reference $N\propto\kappa$, and the gray dotted
line the reference $N\propto\sqrt\kappa$; both are aligned with the
corresponding data series at $\kappa=300$. These indicate theoretical
slopes rather than lines fitted to the data.}
\label{fig:fgm-scaling}
\end{figure}

Figure~\ref{fig:fgm-scaling} shows the number of gradient evaluations $N$
needed to attain $\delta_k\le10^{-8}$. On logarithmic axes, the
gradient-descent data are nearly parallel to the reference
$N\propto\kappa$, and the fast-gradient data to the reference
$N\propto\sqrt\kappa$. It is the difference between these slopes that
illustrates the acceleration in bound~\eqref{eq:fgmrate2}; a single
experiment at fixed $\kappa$ cannot establish this conclusion.

\subsection{What changed}

\textbf{(i) The notion of an optimal method emerged.} Before
Yu.\,E.~Nesterov's work, accelerated methods of
A.\,S.~Nemirovski~\cite{NemYud1979} were known, but they required solving
auxiliary optimization problems in dimension 2 or 1 at each iteration.
Nesterov's fast gradient method required no auxiliary low-dimensional
optimization. The acceleration idea proposed by Yu.\,E.~Nesterov proved
very useful for incorporation into other iterative procedures. Before the
work of A.\,S.~Nemirovski and Yu.\,E.~Nesterov, optimization had largely
been an art: methods were proposed, compared on test sets, and selected.
Afterward, it became meaningful to ask, for each problem class, whether the
theoretical limit had been reached. The fast gradient method is one of the
first examples where a complete answer was obtained: a constructively
established lower bound and a method attaining it.

\smallskip
\textbf{(ii) Acceleration as an element of machine learning.} Inertial
schemes have become an integral part of neural-network training
algorithms. Precision is needed here: the broad family of ``momentum''
methods (including the momentum component of Adam) is not identical to
scheme~\eqref{eq:fgm}; it is more accurate to say that the \emph{idea of
controlled inertia} and its specific implementation known as Nesterov
momentum are part of the standard toolkit and are implemented in all major
libraries. Separate lines of development include extending acceleration to
stochastic optimization (\S\,\ref{sec:accbatch}) and, more recently, to
reinforcement-learning problems~\cite{BelomestnyEtAl2026}.

\smallskip
\textbf{(iii) Refining the constants.} Bound~\eqref{eq:fgmrate} is optimal
in order but not in its constant. The Performance Estimation Problem
approach (numerically solving for the worst-case trajectory as a
semidefinite program) has yielded exact constants for a broad family of
schemes~\cite{TaylorEtAl2017}; Y.~Drori established the exact
information-theoretic lower bound for the last point~\cite{Drori2017},
and D.~Kim and J.~Fessler constructed OGM, which is \emph{exactly}
optimal in this information model and for this criterion, improving
bound~\eqref{eq:fgmrate} by precisely a factor of two~\cite{KimFessler2016}.

\smallskip
\textbf{(iv) The nonconvex case and nonstandard step-size schedules.}
For smooth nonconvex problems, accelerated schemes are known that find a
point with $\nrm{\grad f}\le\eps$ in $\Ot(\eps^{-7/4})$ calls, compared
with $\Oh(\eps^{-2})$ for gradient descent~\cite{CarmonEtAl2017}.
Moreover, step-size schedules irregular in $k$ (``silver'' schedules),
which improve the dependence on $\kappa$ for smooth convex problems,
are being actively studied~\cite{AltschulerParrilo2023}.

\subsection{Composite functions: the first structural model}
\label{sec:composite}

Real-world learning problems have the form
$\min_{x\in Q}\{f(x)+\Psi(x)\}$, where $f$ is smooth and $\Psi$ is a
simple convex function (for example, $\nrm{x}_1$ or the indicator of a set).
This is the simplest \emph{structural} model in the sense of
Table~\ref{tab:models}. Acceleration extends to this case through the
\emph{composite gradient mapping}
\begin{equation}\label{eq:composite}
   T_L(y)=\argmin_{x\in Q}\Bigl\{f(y)+\dotp{\grad f(y)}{x-y}
   +\tfrac L2\nrm{x-y}^2+\Psi(x)\Bigr\},
   \quad
   g_L(y)=L\,B\bigl(y-T_L(y)\bigr),
\end{equation}
where $\nrm{h}=\dotp{Bh}{h}^{1/2}$ with $B=B^{\!\top}\succ0$.
The vector $g_L(y)$ plays the role of a gradient: when $Q=\R^n$ and
$\Psi\equiv0$, it coincides with $\grad f(y)$ for any $L>0$. The
assumption that $\Psi$ is simple means precisely that
operation~\eqref{eq:composite} is implementable: for
$\Psi=\nrm{\cdot}_1$, it is computed explicitly (soft thresholding),
and for $\Psi=\delta_Q$, it becomes a projection onto $Q$; here $\delta_Q$
is the \emph{indicator function} of $Q$, equal to zero on $Q$ and
$+\infty$ outside it. Replacing the gradient step by~\eqref{eq:composite}
preserves bound~\eqref{eq:fgmrate}, and in the scheme with adaptive
selection of $L$ and backtracking it also preserves the bound
$\Oh(\gamma_uL\nrm{x_0-x^{*}}^2/k^2)$, where $\gamma_u>1$ is the factor
by which the current estimate of $L$ is increased after an unsuccessful
trial (usually $\gamma_u=2$)~\cite{Nes2013composite}.

The chronology matters here: paper~\cite{Nes2013composite} was published
in a journal in 2013, but appeared as a CORE preprint in 2007, that is,
\emph{before} the widely known FISTA algorithm~\cite{BeckTeboulle2009};
moreover,~\cite{Nes2013composite} treats a more general setting, including
strong convexity and an inexact oracle.

\subsection{Monotone acceleration and generalized $(L_0,L_1)$-smoothness}
\label{sec:l0l1}

Remark~\ref{rem:nonmono} on nonmonotonicity has an unexpected
continuation. Paper~\cite{NesGasGumDvur2021} constructs the accelerated
method \emph{AGMsDR} (Accelerated Gradient Method with small-Dimensional
Relaxation), whose step is determined by exact minimization of the
function over a low-dimensional (usually one-dimensional) set. Three of
its properties are essential: it is \emph{monotone} (function values do
not increase), requires no knowledge of $L$, and is \emph{primal-dual}
in the sense of \S\,\ref{sec:pd}.

Monotonicity proved to be more than a cosmetic convenience, for the
following reason.

\begin{definition}[$(L_0,L_1)$-smoothness]\label{def:l0l1}
A twice-differentiable function $f$ is called
\emph{$(L_0,L_1)$-smooth} if, for every $x$,
\begin{equation}\label{eq:l0l1}
   \nrm{\grad^2f(x)}\ \le\ L_0+L_1\nrm{\grad f(x)}.
\end{equation}
For $L_1=0$, this is ordinary $L_0$-smoothness. For once-differentiable
functions, a local version is used:
$\nrm{\grad f(y)-\grad f(x)}\le\bigl(L_0+L_1\nrm{\grad f(x)}\bigr)
\nrm{y-x}$ for $\nrm{y-x}\le1/L_1$.
\end{definition}

We distinguish three levels of statements usually associated with this
condition.

\smallskip
\emph{(i) Theoretical examples.} The simplest example of a function
satisfying~\eqref{eq:l0l1} but not globally smooth in the classical sense
is $f(x)=e^x$ on $\R$. Indeed, $f''(x)=f'(x)=e^x$, so
condition~\eqref{eq:l0l1} holds with $L_0=0$ and $L_1=1$, as an equality.
At the same time, $\sup_{x\in\R}|f''(x)|=\infty$, so no finite global
smoothness constant $L$ exists. This class also includes
$\nrm{Ax-b}^p$ for $p>2$, log-sum-exp, and cross-entropy loss functions.
Class~\eqref{eq:l0l1} is closed under a number of operations, allowing
more complicated examples to be constructed.

\smallskip
\emph{(ii) Empirical motivation for gradient clipping.}
Condition~\eqref{eq:l0l1} was proposed by J.~Zhang, T.~He, S.~Sra, and
A.~Jadbabaie~\cite{ZhangEtAl2020} to explain why \emph{gradient clipping}
works so well in training deep networks. The reasoning is as follows:
the optimal step size for~\eqref{eq:l0l1} is of order
$\bigl(L_0+L_1\nrm{\grad f(x)}\bigr)^{-1}$, and in the large-gradient
regime it becomes of order $1/(L_1\nrm{\grad f(x)})$, so the step is
\emph{normalized} by the gradient norm---which is precisely what clipping
does. Thus the heuristic proved to be the optimal step for a suitably
chosen function class.

\practice \emph{(iii) Neural networks.} It is often said that
condition~\eqref{eq:l0l1} ``describes'' the loss functions used in training
large language models. This statement should be understood as a practical
interpretation, not a theorem: a property of individual scalar losses and
linear-composition examples does not automatically extend to the full
nonconvex function of a deep network's parameters, and the global
constants $L_0,L_1$ depend on the architecture, parameter domain, and layer
normalization. Known results for networks are generally local and
architecture-dependent.


Acceleration for class~\eqref{eq:l0l1} long resisted analysis: any
increase in the function value is dangerous here, because the local
smoothness constant grows along with the gradient. This is precisely why
\cite{VankovEtAl2025} uses the \emph{monotone} AGMsDR method
from~\cite{NesGasGumDvur2021} as its underlying construction. Its bound
(in iterations, each including a one-dimensional search) has the
characteristic \emph{additive} structure
\[
   N(\eps)\ =\ \Oh\Bigl(\sqrt{\frac{L_0R^2}{\eps}}
   +(L_1R)^{2/3}\ln\frac{f(x_0)-f^{*}}{\eps}\Bigr),
\]
with $L_1$ leaving the accelerated term unaffected and no exponential
dependence on $L_0$ or $L_1$ arising---this was the main technical
difficulty.

This bound was not final, however. Unlike the $L_0$-smooth case, the
distance to a solution bounds the function-value residual only
exponentially: at distance $R$ from the solution it can reach
$(L_0/L_1^2)\bigl(e^{L_1R}-L_1R-1\bigr)$, and the bound above then
guarantees reaching the region where the problem behaves as an
$L_0$-smooth one only after $\Oh\bigl((L_1R)^{5/3}\bigr)$ iterations.
The exact cost of the parameter $L_1$ was found very
recently~\cite{RudomanovEtAl2026}: for $L_0R^2\ge\eps$, the optimal
number of oracle calls is
\[
   N(\eps)\ =\ \widetilde\Theta\Bigl(
   \underbrace{L_1R}_{\substack{\text{reaching the region where the problem}\\
   \text{behaves as an $L_0$-smooth one}}}
   \ +\
   \underbrace{\sqrt{L_0R^2/\eps}}_{\text{accelerated phase}}\Bigr),
\]
where the tilde hides factors logarithmic in $L_1R$. The lower bound holds
for \emph{any} deterministic first-order method, while the upper bound is
given by a two-phase scheme. First, a monotone accelerated method based
on estimate sequences, close to AGMsDR but using an inexact
one-dimensional search involving only gradients, finds a point where
$f(x)-f^{*}\le L_0/(128L_1^2)$ in $\Ot(L_1R)$ oracle calls; then an
accelerated method with constant step size $1/(2L_0)$~\cite{Tyurin2026}
finishes in $\Oh\bigl(\sqrt{L_0R^2/\eps}\bigr)$ calls.
The linear dependence on $L_1R$ is intuitively natural: where the gradient
is large, the linear model of the function can be trusted only over a
distance of order $1/L_1$ (cf.\ Definition~\ref{def:l0l1}), so traversing
a distance $R$ requires on the order of $L_1R$ steps---and the lower bound
shows that no deterministic first-order method can substantially shorten
this path. The title of~\cite{RudomanovEtAl2026}, ``Descending from
Everest'', captures this picture well: the initial height may be
exponential, while the descent takes only $\Ot(L_1R)$ steps.

All these schemes are nonadaptive: they use knowledge of the constants
$L_0$ and $L_1$ (and the optimal one also uses $R$ and $\eps$).
A fully adaptive method was constructed by E.~Borodich and
D.~Kovalev~\cite{BorodichKovalev2026}, who combined Nesterov acceleration
with the step-size rule from the golden ratio algorithm (GRAAL), proposed
by Y.~Malitsky~\cite{Malitsky2020} for variational inequalities. The step
size is determined by a local estimate of inverse curvature from two
already computed points,
\[
   \lambda(x,z)=\frac{2\bigl[f(x)-f(z)-\dotp{\grad f(z)}{x-z}\bigr]}
   {\nrm{\grad f(x)-\grad f(z)}^2},
\]
and can either decrease or increase (no faster than geometrically);
neither $L$, $L_0$, $L_1$ nor backtracking is needed. For $L$-smooth
functions the method is optimal:
$N(\eps)=\Oh\bigl(\sqrt{LR^2/\eps}\bigr)$ with an additive term
$\Oh\bigl(\ln\frac{1}{\eta_0L}\bigr)$ for choosing a sufficiently small
initial step size $\eta_0$. On class~\eqref{eq:l0l1},
\[
   N(\eps)=\Oh\Bigl(\sqrt{\frac{L_0R^2}{\eps}}+(L_1R)^3
   +\bigl(1+(L_1R)^2\bigr)\ln\frac{1}{\eta_0L_0}\Bigr),
\]
so the same additive structure is preserved (although the
$\eps$-independent term $(L_1R)^3$ is noticeably larger here than the
optimal $L_1R$), without any prior information about the constants.
Interestingly, monotonicity is not used here: a decreasing Lyapunov
function of the method keeps the trajectory within a ball of radius
$\Oh(R)$ around the solution.

\subsection{Inexact oracles and batching}
\label{sec:accbatch}

Extending acceleration to stochastic problems
$\min_{x\in Q}\Exp_\xi f(x,\xi)$, where only an unbiased stochastic
gradient with variance $\sigma^2$ is available, encounters a conflict
between acceleration and inexactness: the more ``aggressive'' a scheme,
the more strongly it accumulates oracle errors. The language for
describing this conflict is supplied by the notion of an \emph{inexact
oracle}, introduced by O.~Devolder, Yu.\,E.~Nesterov, and
F.~Glineur~\cite{DevolderNesterovGlineur2014}. In the stochastic setting,
the following modification, proposed in the survey~\cite{AvetEtAl2026},
is convenient.

\begin{definition}[$(\delta_1,\delta_2,L)$-oracle; \cite{AvetEtAl2026}]
\label{def:inexact}
A $(\delta_1,\delta_2,L)$-oracle is said to be available at $x$ if one can
compute a vector $\grad_\delta f(x)$ such that, for all $x,y\in Q$,
\begin{equation}\label{eq:inexactoracle}
   -\delta_1(x,y)\ \le\ f(y)-f(x)-\dotp{\grad_\delta f(x)}{y-x}
   \ \le\ \tfrac L2\nrm{y-x}_2^2+\delta_2(x,y),
\end{equation}
with $\Exp[\delta_1(x,y)\mid x]=0$ and
$\Exp\,\delta_2(x,y)\le\delta$.
\end{definition}

The left inequality in~\eqref{eq:inexactoracle} is the
$\delta_1$-subgradient condition (B.\,T.~Polyak), reducing to ordinary
convexity for $\delta_1=0$; the right inequality with $\delta_2=0$ follows
from $L$-smoothness. The difference from~\cite{DevolderNesterovGlineur2014}
is that the errors are separated into ``zero in expectation'' and
``positive in expectation'', matching the stochastic gradient setting.

Suppose a method using such an oracle has the bound
\begin{equation}\label{eq:inexactrate}
   \Exp f(x^N)-f^{*}=\Oh\!\left(\frac{LR^2}{N^{\alpha}}
   +N^{\zeta}\delta\right).
\end{equation}
For gradient descent in model~\eqref{eq:inexactoracle},
$\alpha=1$, $\zeta=0$: the inexactness \emph{does not} accumulate.
For the standard accelerated schemes under consideration in the same
model, $\alpha=2$, $\zeta=1$: the price of acceleration is linear
accumulation of error~\cite{DevolderNesterovGlineur2014}. We do not claim
that \emph{every} possible accelerated scheme behaves this way under
\emph{every} error model; the statement concerns a specific pair
``scheme $+$ model''.

This yields the technique of \emph{batching}. Take
$\grad_\delta f(x)=\frac1r\sum_{j=1}^r\grad_xf(x,\xi^j)$; then
$\delta=\sigma^2/(2Lr)$, and balancing the two terms
in~\eqref{eq:inexactrate} gives, for the accelerated method,
\begin{equation}\label{eq:batch}
   N=\Oh\!\left(\sqrt{\frac{LR^2}{\eps}}\right),\qquad
   r=\max\left\{1,\ 
   \Oh\!\left(\frac{\sigma^2R}{\sqrt L\,\eps^{3/2}}\right)\right\},
\end{equation}
\begin{equation}\label{eq:batch2}
   T=N\cdot r=\Oh\!\left(\sqrt{\frac{LR^2}{\eps}}
   +\frac{\sigma^2R^2}{\eps^{2}}\right).
\end{equation}
The condition $r\ge1$ is essential: for $\sigma=0$, a formal calculation
would give $r=0$, even though the method must still perform at least $N$
gradient steps; in~\cite{AvetEtAl2026}, this is specified by the rule
$T:=\max\{T,N\}$. The total number $T$ of stochastic oracle calls is
optimal as the sum of deterministic and variance contributions, while
the number of \emph{sequential} iterations $N$ has decreased from
$\eps^{-2}$ (ordinary stochastic gradient descent) to $\eps^{-1/2}$,
with the entire batch computed in parallel. This explains the practical
efficiency of large batches. For a systematic presentation, see
\cite[\S\,4.5]{AvetEtAl2026}; the same technique extends to
gradient-free (\S\,\ref{sec:zo}), distributed, and tensor
(\S\,\ref{sec:tensor}) optimization.

\subsection{Extending accelerated methods to structured saddle-point problems}
\label{sec:sliding}

The final topic of this section returns us to \S\,\ref{sec:complexity}.
Acceleration also extends to problems in which oracle calls are of
different types and have different costs, with \emph{separate} optimal
bounds for each type.

Consider problem~\eqref{eq:spp} and the condition numbers introduced
in~\eqref{eq:kappas}. Lower bound~\eqref{eq:lbspp} says that one cannot
use fewer than $\sqrt{\kappa_x}$ evaluations of $\grad f$,
$\sqrt{\kappa_y}$ evaluations of $\grad g$, and
$\sqrt{\kappa_{xy}}$ multiplications by $B$ (up to a logarithmic factor).

The optimal method in~\cite{BorGasKov2025} is constructed as follows.
The problem is reduced to a monotone variational inequality with a
finite-sum structure on a closed convex set $C_z$ of the combined variable
$z=(x,y)$ (for problem~\eqref{eq:spp}, this is the whole space
$\R^{d_x}\times\R^{d_y}$):
\[
   \text{find } z^{*}\in C_z:\qquad
   p(z^{*})-p(z)+\dotp{Q(z)}{z^{*}-z}\le0\quad \forall z\in C_z,
\]
where
\[
   p=\sum_{i=1}^{n}p_i,\qquad Q=\sum_{i=1}^{n}Q_i,
\]
with $p_i$ smooth convex functions with constants $L_i$, and $Q_i$
monotone Lipschitz operators with constants $M_i$. Next,
reparameterization~\eqref{eq:ht}---the same one that yielded
Theorem~\ref{th:accsimple}---is applied \emph{recursively} to each term,
resulting in a scheme with nested loops (\emph{sliding}). The resulting
bound~\cite[Corollary~4.6]{BorGasKov2025} is stated in a relative
normalization. Let $z_{\text{in}}$ be the initial point, and let
$\nrm{\cdot}_P$ be the norm induced by a fixed positive definite operator
$P$ specifying the geometry of the method. Then the output point
$z_{\text{out}}$ satisfies
\begin{equation}\label{eq:slidingcrit}
   p(z_{\text{out}})-p(z)+\dotp{Q(z)}{z_{\text{out}}-z}
   \ \le\ \eps\,n\,\nrm{z_{\text{in}}-z}_P^2
   \qquad\text{for all } z\in C_z,
\end{equation}
and this requires no more than
\begin{equation}\label{eq:slidingbound}
   6n\,\max\Bigl\{\sqrt{L_i/\eps},\ M_i/\eps,\ 1\Bigr\}
   \quad\text{evaluations of }\grad p_i\text{ and }Q_i
   \qquad (1\le i\le n).
\end{equation}
It is essential that~\eqref{eq:slidingcrit} uses relative accuracy
(weighted by distance): when $C_z$ is unbounded, taking the supremum of the
left-hand side over all $z$ is not possible---that supremum may be
$+\infty$. The substantive meaning of~\eqref{eq:slidingbound} is that
\emph{each} term is processed with its own complexity rather than the
worst of them. For problem~\eqref{eq:spp}, $n=3$ terms suffice, and
\eqref{eq:slidingbound}, together with restarting, leads to bounds
corresponding to the lower bounds~\eqref{eq:lbspp}; literal agreement
between the upper and lower bounds additionally requires the norms,
radii, factor $n$, and logarithmic accuracy to be matched.

\emph{Complexity separation} is an interesting phenomenon. If a problem
has several types of oracle calls with different costs (an expensive
gradient and a cheap matrix multiplication; local computation and
communication in distributed optimization; gradients with respect to
``large'' and ``small'' blocks of variables), then an optimal method
should access each source \emph{at its own frequency}, skipping the
others. We will encounter the same idea in \S\,\ref{sec:cd} (coordinate
methods) and \S\,\ref{sec:superfast} (one Hessian evaluation per outer
iteration).

Both lines of complexity separation---by oracle type and by blocks of
variables (the min-min problem to be discussed in \S\,\ref{sec:cd})---were
recently unified in~\cite{NguyenEtAl2026}, and at the same time extended
to saddle-point problems with \emph{nonbilinear} coupling,
\begin{equation}\label{eq:sppH}
   \min_{x\in\R^{d_x}}\max_{y\in\R^{d_y}}\ \bigl\{f(x)+H(x,y)-g(y)\bigr\},\qquad
   x=(x_1,\dots,x_I),\quad y=(y_1,\dots,y_J).
\end{equation}
Here $f$ and $g$ are smooth and strongly convex, with block strong
convexity constants $\mu^f_i$, $\mu^g_j$ and smoothness constants $L^f_i$,
$L^g_j$ (separability is not assumed); the convex-concave coupling
function $H$ has block-Lipschitz partial gradients; and each partial
gradient---$\grad_{x_i}f$, $\grad_{y_j}g$, $\grad_{x_i}H$,
$\grad_{y_j}H$---is counted as a separate oracle with its own cost. The
randomized extragradient method in~\cite{NguyenEtAl2026} calls
$\grad_{x_i}f$ only
$\Oh\bigl(\sqrt{L^f_i/\mu^f_i}\,\ln\frac{\kappa_{\max}}{\eps}\bigr)$
times (and analogously for $\grad_{y_j}g$), where $\kappa_{\max}$ is the
largest of the block condition numbers $L^f_i/\mu^f_i$ and
$L^g_j/\mu^g_j$. The numbers of calls to the coupling's partial gradients
are determined not by these condition numbers but by the dimensionless
quantities
\[
   \kappa^x_{H,i}=\Bigl(\sum_{t=1}^{I}\frac{(L^{xx}_{it})^2}{\mu^f_i\mu^f_t}
   +\sum_{s=1}^{J}\frac{(L^{xy}_{is})^2}{\mu^f_i\mu^g_s}\Bigr)^{1/2},
\]
where $L^{xx}_{it}$ and $L^{xy}_{is}$ are the Lipschitz constants of
$\grad_{x_i}H$ with respect to blocks $x_t$ and $y_s$, respectively.
Each constant is divided by the geometric mean of the strong convexity
parameters of the two blocks it couples (block $x_i$ and block $x_t$ or
$y_s$); in this sense, the coupling constants are normalized by the
curvature of $f$ and $g$: the stronger the convexity of the blocks, the
less their coupling constrains the method. The method calls
$\grad_{x_i}H$
$\Oh\bigl((1+\kappa^x_{H,i}\sqrt{\Delta_i})\ln\frac{\kappa_{\max}}{\eps}\bigr)$
times, where $\Delta_i\ge1$ is a correction for block overlap: for each
block on which $\grad_{x_i}H$ depends, count how many partial gradients
of the coupling depend on it, and take the largest such number (with no
overlap, $\Delta_i=1$). The quantities $\kappa^y_{H,j}$ for
$\grad_{y_j}H$ are defined symmetrically. In the single-block bilinear
case $H(x,y)=\dotp{y}{Bx}$, we obtain
$\kappa^x_{H}=\kappa^y_{H}=\nrm{B}/\sqrt{\mu_x\mu_y}$, that is,
$\sqrt{\kappa_{xy}}$ from~\eqref{eq:kappas} with
$\mu_{xy}=\mu_{yx}=0$, recovering the bounds
of~\cite{BorGasKov2025} for strongly convex $f$ and $g$; matching lower
bounds have been proved for bilinear coupling between disjoint pairs of
blocks (in which case all $\Delta_i=1$). Without strong convexity,
regularization gives $\Ot(\eps^{-1/2})$ calls to the gradients of $f$
and $g$, and $\Ot(\eps^{-1})$ calls to the coupling gradients, also
consistent with the lower bounds for the smooth and bilinear parts.

Even the complexity scales themselves can be separated
(\S\,\ref{sec:cutting}). Suppose that in
$\min_{x\in\R^n}\max_{y\in Y}f(x,y)$ the block $y$ is low-dimensional:
$Y\subset\R^m$, $m\ll n$, and $D_y$ is the diameter of $Y$. Suppose also
that $f$ is $\mu$-strongly convex and $L$-smooth in $x$
($\kappa=L/\mu$), but only concave and Lipschitz in $y$, with constant
$M_y$. It is then natural to use a cutting-plane method for $y$ and an
accelerated method for $x$. If each inner problem is solved to high
accuracy (of order $\eps$), such a nested scheme requires
$\Oh\bigl(m\sqrt\kappa\,\ln^2\frac{M_yD_y}{\eps}\bigr)$ evaluations of
$\grad_xf$ (here and below, terms independent of $\eps$ are omitted).
Paper~\cite{FominGasnikov2026} shows that high-accuracy inner solutions
are unnecessary: it suffices to maintain a strongly convex lower model
of one of the sections $f(\cdot,y)$ and ``transfer'' it to new sections,
so that it continuously provides a computable accuracy certificate.
Then $\Oh\bigl(m\ln\frac{M_yD_y}{\eps}\bigr)$ oracle calls in $y$ and
$\Oh\bigl((m+\sqrt{m\kappa})\ln\frac{M_yD_y}{\eps}\bigr)$ in $x$
suffice; for $m=1$, both bounds match the lower bounds for the
corresponding blocks.

%% file: en/r4_ad.tex

\section{Fast automatic differentiation (1984)}
\label{sec:ad}

\subsection{Formulation of the question}

In \S\,\ref{sec:complexity}, complexity was measured by the number of
oracle calls. But the oracle itself must be implemented, and here a
question arises whose answer seemed almost obvious until the 1980s,
but turned out to be unexpected.

Suppose a function $f\colon\R^n\to\R$ is given by a \emph{program}: a
finite sequence of elementary operations. How many operations are needed
to compute the \emph{entire} gradient $\grad f(x)$, a vector of $n$
numbers?

The ``obvious'' answer is $n$ times as many as for evaluating the function
itself. Indeed, the standard approach, finite differences,
\[
   \frac{\partial f}{\partial x^{(i)}}(x)\approx
   \frac{f(x+\theta e_i)-f(x)}{\theta},\qquad i=1,\dots,n,
\]
requires $n+1$ function evaluations and introduces an unavoidable error:
small $\theta$ increases rounding error, while large values increase
approximation error. The correct answer is: \emph{the entire gradient
costs as much as the function itself}, up to an absolute constant, and is
computed \emph{exactly}. This result was established in the work of
K.\,V.~Kim, Yu.\,E.~Nesterov, and
B.\,V.~Cherkassky~\cite{KimNesCher1984} (see
also~\cite{KimNesSkokCher1984} and the exposition
in~\cite{NesRus1989}); in the English-language literature it is known as
the \emph{cheap gradient principle}, the technique as reverse-mode
automatic differentiation, and, in the context of neural networks, as
\emph{backpropagation}.

Related ideas were expressed independently and in other forms by
R.~Wengert (1964)~\cite{Wengert1964},
B.~Speelpenning (1980)~\cite{Speelpenning1980}, and
W.~Baur and V.~Strassen (1983)~\cite{BaurStrassen1983}; the relationship
between these results is discussed in Remark~\ref{rem:bs}. In the Russian
literature, this subject, under the name \emph{fast automatic
differentiation}, was systematically developed by the school of
Academician Yu.\,G.~Evtushenko~\cite{Evtushenko1982,Evtushenko2013},
with applications integral to the work from the outset: fast
differentiation formulas were applied to optimal control, inverse
problems, and sensitivity analysis of complex engineering systems.

\subsection{The computational model}

A precise statement requires fixing the computational model: without
one, the assertion that ``the constant does not depend on the program's
structure'' is meaningless.

Fix a finite set $\Phi$ of \emph{elementary operations}: the arithmetic
operations $+,-,\times,/$ and a finite list of functions
($\exp$, $\ln$, $\sin$, $\sqrt{\;\cdot\;}$, etc.). Assume that all
operations in $\Phi$ have bounded arity (number of arguments) and satisfy
the following property: \emph{all partial derivatives of an operation
$\phi\in\Phi$ can be computed in $\Oh(1)$ elementary operations from the
already computed values of its arguments and result}. For example, for
$\phi(u,v)=uv$, these are $\partial\phi/\partial u=v$ and
$\partial\phi/\partial v=u$, numbers already known; for
$\phi(u)=\exp u$, the derivative equals the value $\phi(u)$ itself; and
for $\phi(u)=\sqrt u$, it equals $1/(2\phi(u))$.

A \emph{program for evaluating} $f$ is a collection of
\emph{intermediate variables}
\[
   \underbrace{v_1,\dots,v_n}_{\text{inputs: } v_i=x^{(i)}},\quad
   v_{n+1},\dots,v_{n+m},
\]
in which each variable with index $i>n$ is computed as
\begin{equation}\label{eq:graph}
   v_i=\phi_i\bigl(v_j:\ j\prec i\bigr),\qquad \phi_i\in\Phi,
\end{equation}
and the result is $v_{n+m}=f(x)$. The notation $j\prec i$ means:
\emph{the variable $v_j$ is one of the arguments of the operation that
computes $v_i$}. The relation $\prec$ is acyclic ($j\prec i$ implies
$j<i$), so~\eqref{eq:graph} defines a directed acyclic graph on the index
set $\{1,\dots,n+m\}$, the \emph{computational graph}: vertices are
intermediate variables, and edges are ``argument--result'' relations.
The number of elementary operations required to compute the $m$ additional
variables is the program's computational cost; denote it by $T$.

\begin{example}\label{ex:speel}
Speelpenning's example: $f(x)=x^{(1)}x^{(2)}\cdots x^{(n)}$.
The program is $v_{n+1}=v_1v_2$, $v_{n+2}=v_{n+1}v_3$, \dots,
$v_{2n-1}=v_{2n-2}v_n$, requiring $T=n-1$ multiplications; here
$1\prec n+1$, $2\prec n+1$, then $n+1\prec n+2$, $3\prec n+2$, and so
on. The partial derivatives are
$\partial f/\partial x^{(i)}=\prod_{j\ne i}x^{(j)}$; evaluating each
independently takes $\Oh(n^2)$ operations. Yet the entire gradient can be
obtained in $\Oh(n)$ operations by computing prefix and suffix products.
Reverse mode finds such a procedure \emph{automatically}, knowing
nothing about prefixes and suffixes.
\end{example}

\subsection{The reverse sweep}

Introduce the \emph{adjoint variables}
\[
   \bar v_i=\frac{\partial v_{n+m}}{\partial v_i},
   \qquad i=1,\dots,n+m,
\]
where the derivative is taken with $v_i$ regarded as an independent
argument of the part of the graph that depends on it. The chain rule
gives the recurrence
\begin{equation}\label{eq:adjoint}
   \bar v_{n+m}=1,\qquad
   \bar v_j=\sum_{i:\ j\prec i}\ \bar v_i\cdot
   \frac{\partial\phi_i}{\partial v_j}\bigl(v_k:\ k\prec i\bigr),
\end{equation}
which must be applied in decreasing index order
$j=n+m-1,\dots,1$, that is, traversing the computational graph
\emph{in reverse}. At the end, we obtain
$\partial f/\partial x^{(i)}=\bar v_i$, $i=1,\dots,n$.

Let us estimate the cost. The sum in~\eqref{eq:adjoint} is over all edges
leaving vertex~$j$; thus every edge is processed exactly once during the
entire reverse sweep, with $\Oh(1)$ operations per edge (one
multiplication, one addition, and evaluation of the local derivative
$\partial\phi_i/\partial v_j$, costing $\Oh(1)$ by the assumption on
$\Phi$). The number of edges is at most $\mathrm{const}\cdot m$, since
the arity of the operations is bounded. The result is:

\begin{theorem}[\cite{KimNesCher1984}]\label{th:ad}
In the computational model described above, there exists an absolute
constant $C$, depending only on the set of elementary operations $\Phi$
and not on the function, the dimension $n$, or the structure of the
computational graph, such that if $f$ is evaluated by a program in $T$
elementary operations, then the entire vector $\grad f(x)$, together
with the value $f(x)$, can be computed in at most $C\cdot T$ elementary
operations.
\end{theorem}

The qualification about a fixed model is essential. If $\Phi$ is allowed
to include operations whose derivatives are themselves expensive,
Theorem~\ref{th:ad} no longer holds; this is exactly what happens when an
``elementary operation'' is a call to an external package
(see \S\,\ref{sec:adnograd} below).

The price of Theorem~\ref{th:ad} is \emph{memory}: the reverse sweep needs
the values $v_i$ obtained in the forward sweep, requiring $\Oh(T)$ memory
cells. This limitation is the bottleneck in training large neural
networks, and addressing it is precisely the purpose of checkpointing,
which trades memory for recomputation.

\begin{remark}[on the Baur--Strassen constant]\label{rem:bs}
The paper of W.~Baur and V.~Strassen~\cite{BaurStrassen1983} obtains a
result with the explicit constant $C=3$, but under a \emph{different}
complexity measure: for rational functions evaluated by a straight-line
program, only \emph{nonscalar} operations are counted (multiplications
and divisions in which neither argument is a constant). The universal
constant for the total number of arithmetic operations is different.
Thus ``$C=3$'' cannot be compared directly with the constant in
Theorem~\ref{th:ad}: these are statements about different complexity
measures.
\end{remark}

\subsection{What changed}

The tempting conclusion (``gradients are cheap, so one should use
gradients'') explains little by itself. There are three substantive
consequences, and they are subtler.

\smallskip
\textbf{(i) The oracle model of \S\,\ref{sec:complexity} becomes an
adequate measure of arithmetic work.} In a typical iteration of a
first-order method, the cost splits into gradient evaluation and
``everything else'': a few operations on vectors of length $n$ and,
possibly, a prox-mapping. In real-world problems, the first term generally
dominates: $T\gg n$. Thus the number of gradient-oracle calls is
\emph{proportional} to the total arithmetic work, and bounds such as
$N(\eps)=\Oh(\sqrt{LR^2/\eps})$ can be read as computation-time bounds.
This is why the oracle complexity of first-order methods is customarily
measured in gradient evaluations. Without Theorem~\ref{th:ad}, this
convention would be unjustified: if the gradient cost $\Oh(nT)$, the
separation ``number of iterations $\times$ iteration cost'' would have
to be handled differently, as it must be in
\S\S\,\ref{sec:cd},~\ref{sec:superfast}.

\smallskip
\textbf{(ii) An explanation emerges for why first order is special.}
Compare the information ``density'' of different oracles. A function
value is one number for $T$ operations. A gradient is $n$ numbers for
$\Oh(T)$ operations, or $n$ times as much information at the same cost.
The full Hessian is $n^2$ numbers, and reverse mode no longer rescues us:
standard bounds give $\Oh(nT)$. The full $p$th derivative is $n^p$
numbers for $\Oh(n^{p-1}T)$ operations. In other words, it is precisely
first order that is ``free'' in the sense of Theorem~\ref{th:ad}, which
is why first-order methods have become the principal tool for
high-dimensional problems.

An important qualification is needed here, to which we return in
\S\,\ref{sec:tensor}. One must distinguish the \emph{full tensor}
$D^pf(x)$, with $n^p$ components, from a \emph{fixed number of its
contractions} with specified directions,
\[
   \underbrace{D^pf(x)[h_1,\dots,h_p]\in\R}_{\text{all arguments
   filled}},
   \qquad
   \underbrace{D^pf(x)[h_1,\dots,h_{p-1},\,\boldsymbol\cdot\,]
   \in(\R^n)^{*}\simeq\R^n}_{\text{one argument free}} .
\]
For a scalar function $f$, a fully evaluated $p$-linear form is a number;
a vector results when one argument is left free. For fixed $p$, each such
contraction can be computed in $\Oh(T)$ operations by the same reverse
mode (applied to an auxiliary function of one variable). The tensor
methods of \S\,\ref{sec:tensor} require precisely contractions of the
form $D^pf(x)[h,\dots,h,\boldsymbol\cdot\,]$, rather than the full
tensor, which is one reason they are implementable. There is no
contradiction with the bound $\Oh(n^{p-1}T)$: by enumerating $n^{p-1}$
vector contractions along basis directions, we recover the full tensor.

\smallskip
\textbf{(iii) Practice.} Reverse mode underlies the training of all
modern neural networks: relation~\eqref{eq:adjoint} is applied to the
network's computational graph, and all the technical work of automatic
differentiation libraries consists in constructing this graph and
organizing the reverse sweep under memory constraints. The monograph of
A.~Griewank and A.~Walther~\cite{GriewankWalther2008} is the standard
reference on algorithmic differentiation.

\subsection{The nonsmooth case}

The reverse sweep~\eqref{eq:adjoint} is a formal application of the chain
rule and ``knows'' nothing about differentiability. What happens if the
program contains $\max$, $\min$, $|\cdot|$, or $\mathrm{ReLU}$? Libraries
return something in this case, and this ``something'' works in
practice---but what exactly is it?

The answer ``an element of the Clarke subdifferential'' is incorrect:
for the Clarke subdifferential, the chain rule holds only as an inclusion,
and formal ``automatic subdifferentiation'' can return a vector that is
\emph{not} a generalized gradient of the composition. Three related but
nonidentical objects must be distinguished here:

\begin{itemize}[leftmargin=1.6em,itemsep=3pt]
\item the \emph{Clarke subdifferential} $\partial^{\mathrm C}f(x)$:
      the convex hull of limits of gradients at nearby differentiability
      points; its chain rule holds only as an inclusion;
\item the \emph{lexicographic derivative} in the sense of
      Yu.\,E.~Nesterov~\cite{Nes2005lex}: for the class of
      \emph{lexicographically smooth} functions (which includes all
      compositions of smooth functions with $\max$, $\min$, and
      $|\cdot|$), a single-valued object is defined for which the chain
      rule holds \emph{with equality} and which is correctly computed
      by forward and reverse sweeps;
\item a \emph{conservative field} in the sense of J.~Bolte and
      E.~Pauwels~\cite{BoltePauwels2021}: a set-valued mapping consistent
      with $f$ along almost all curves. It turns out that automatic
      differentiation of nonsmooth programs generates precisely a
      conservative field, and this suffices to prove almost-sure
      convergence of the stochastic gradient method.
\end{itemize}

Thus the question ``what does backpropagation return on nonsmooth
functions?'' has received two different substantive answers, the first
of which,~\cite{Nes2005lex}, was given fifteen years before the question
became widely relevant.

\subsection{When no gradient is available}
\label{sec:adnograd}

Finally, we note a topic to which we return in \S\,\ref{sec:zo}.
Theorem~\ref{th:ad} assumes that the function is given by a program in a
fixed computational model. If this is not the case---if $f(x)$ is
obtained by physical measurement, running a third-party simulation
package, or solving an auxiliary problem to a prescribed accuracy---then
the reverse sweep is inapplicable in principle: there is no computational
graph to which it could be applied.

Crucially, in such a situation, \emph{the function value itself is
usually inexact}. This changes the formulation qualitatively: the oracle
returning an exact value is replaced by one with controllable accuracy
$\delta$ and cost $c(\delta)$ increasing as $\delta\to0$, and one must
optimize not the number of calls but the total time $N\cdot c(\delta)$.
This three-parameter formulation is the subject of a separate line of
research~\cite{SuvorikovaEtAl2026}; it is discussed at the end of
\S\,\ref{sec:zo}.

%% file: en/r5_ipm_pi2.tex

\section{Lexicographic differentiation of nonsmooth functions
(1986 -- 1987)}
\label{sec:lexdiff}

\subsection{The question: what remains of the chain rule}

In \S\,\ref{sec:ad}, the function was specified by a program, and its entire gradient
was computed by a reverse pass through the graph of that program. The algorithm
was based on the ordinary rule for differentiating a composition.
But actual computational graphs contain the operations
$\max$, $\min$, $|\cdot|$, and piecewise-defined functions. Even
$f(t)=|t|$ has no gradient at zero. Can it be replaced by some
generalized gradient while retaining the main algorithmic property---
an exact chain rule?

The best-known answer in nonsmooth analysis uses the Clarke subdifferential
$\partial^{\mathrm C}f(x)$---the convex hull of limits of
gradients at nearby differentiability points. This is a powerful
theoretical object, but it is inconvenient for automatic differentiation
for two reasons. First, sets must be propagated through the graph
instead of a single vector. Second, the chain rule for the Clarke
subdifferential generally takes the form of an inclusion rather than
an equality. Simply substituting arbitrary elements of local
subdifferentials can therefore give an incorrect answer.

For example, the function
\[
        q(t)=|t|-|t|\equiv0
\]
has $\partial^{\mathrm C}q(0)=\{0\}$. If, however, the two occurrences of the absolute value
are differentiated independently, choosing respectively the numbers
$1$ and $-1$ from $\partial^{\mathrm C}|\,\cdot\,|(0)=[-1,1]$, the formal
rule gives $1-(-1)=2$. The correlation between two identical
branches of the program has been lost. Thus it is not enough to decide \emph{which
set} should be regarded as the derivative: a consistent way of selecting
elements from it throughout the computational graph is needed.

This was precisely the problem addressed by Yu.\,E.~Nesterov's work in the mid-1980s.
A detailed treatment appeared in a 65-page CEMI preprint in 1986; in 1987,
the paper ``The technique of nonsmooth differentiation''~\cite{NesNonsmooth1987} and
a joint short communication with B.\,V.~Cherkassky giving explicit calculus rules
were published~\cite{NesCherLex1987}. The theory was subsequently
presented in the book~\cite{NesRus1989}, and later systematized and
substantially extended in~\cite{Nes2005lex}. Below we first
reproduce the original language of cones and then its modern formulation,
which retains the same idea of resolving ambiguities successively.

\subsection{From derivatives along cones to ordered directions}

In~\cite{NesCherLex1987}, the starting point was a linear
approximation within a cone rather than throughout a neighborhood of a point. Let
$K\subset\R^n$ be a cone with vertex at zero. A function
$f\colon\R^n\to\R$ is called differentiable at $x$ along the cone $K$
if there exists a vector $g=g(x;K)$ such that
\begin{equation}\label{eq:lexcone}
   f(x+h)=f(x)+\dotp{g}{h}+o(\nrm{h}),
   \qquad h\to0,\quad h\in K .
\end{equation}
The vector $g(x;K)$ was called a generalized gradient along the cone. For a smooth
function it coincides with the ordinary gradient; for a convex function every
such vector is an extreme point of the subdifferential. If $K$ has
nonempty interior, then $g(x;K)$ is unique. For a nonconvex function,
under the assumptions of the original paper, it is an element of the
Clarke subdifferential.

The main construction of the 1987 paper specifies a cone using
an \emph{ordered basis}
$S=(s_1,\dots,s_n)$. Consider the curve
\begin{equation}\label{eq:lexcurve}
   z_S(t)=t\,s_1+t^2s_2+\cdots+t^n s_n,\qquad 0<t<t_0,
\end{equation}
and the convex cone
\begin{equation}\label{eq:lexcone-def}
   K(S,t_0)=\operatorname{conv}
   \bigl\{\alpha z_S(t):\alpha>0,\ 0<t<t_0\bigr\}.
\end{equation}
It is called a \emph{lexicographic cone of width $t_0$}.
The cone has nonempty interior; it narrows as the width decreases,
and every open cone contains some
$K(\widetilde S,t_0)$ of positive width. The first direction $s_1$
has the highest priority; if it does not resolve the ambiguity between
active branches, $s_2$ is used, then $s_3$, and so on. The powers of
$t$ in~\eqref{eq:lexcurve} explain the word
``lexicographic'': the sign of the first nonzero quantity determines the outcome
of the comparison, and subsequent quantities are used only when the preceding
ones are equal.

This is particularly transparent in the convex case. Set
\begin{equation}\label{eq:lexfaces}
   G_0=\partial f(x),\qquad
   G_i=\operatorname*{Arg\,max}_{u\in G_{i-1}}
         \dotp{s_i}{u},
   \qquad i=1,\dots,n .
\end{equation}
First, only those subgradients that are maximal in direction
$s_1$ are retained; among these, those maximal in direction $s_2$ are selected,
and so on.

\begin{proposition}[\cite{NesCherLex1987}]
\label{pr:lexselect}
For a finite convex function, the sets in~\eqref{eq:lexfaces}
are nonempty, and $G_n$ consists of a single vector. If the function has
a generalized gradient $g(x;K(S,t_0))$ in the sense of~\eqref{eq:lexcone}, then
$g(x;K(S,t_0))\in G_n$.
\end{proposition}

Thus, whenever the cone derivative exists, a finite
sequence of scalar comparisons identifies it uniquely.
Permuting the basis vectors may change the result:
uniqueness holds once an \emph{ordered} basis is fixed,
not independently of it. The converse of Proposition~\ref{pr:lexselect}
is generally false. For example, for
$f(x)=\nrm{x}_2$ at zero and $S=(e_1,e_2)$, procedure
\eqref{eq:lexfaces} selects $G_2=\{e_1\}$, but $e_1$ does not give a linear
approximation on any $K(S,t_0)$ of positive width. Indeed,
along $h=\alpha(te_1+t^2e_2)$ the relative error is
\[
 \frac{\nrm{h}_2-\dotp{e_1}{h}}{\nrm{h}_2}
 =
 \frac{\sqrt{1+t^2}-1}{\sqrt{1+t^2}}>0
\]
for every fixed $t\in(0,t_0)$ and does not vanish as
$\alpha\downarrow0$.

In the later systematization~\cite{Nes2005lex}, the idea of ordered
directions was substantially extended: the recursive construction
was also defined for functions that may have no linear approximation
along a cone. Let $X$ be open. For a locally Lipschitz function
$f\colon X\subset\R^n\to\R^m$, denote the one-sided directional derivative by
\[
   f'(x;d)=\lim_{\tau\downarrow0}
   \frac{f(x+\tau d)-f(x)}{\tau}.
\]
Fix a matrix of directions
$M=[m_1\ \cdots\ m_k]\in\R^{n\times k}$ and define
a \emph{sequence of homogenizations}
\begin{equation}\label{eq:lexhom}
 \begin{aligned}
   f^{(0)}_{x,M}(d)&=f'(x;d),\\
   f^{(j)}_{x,M}(d)&=
   \bigl[f^{(j-1)}_{x,M}\bigr]'(m_j;d),
   \qquad j=1,\dots,k .
 \end{aligned}
\end{equation}
In other words, the function itself is first differentiated at $x$,
then the resulting positively homogeneous function at $m_1$,
the next one at $m_2$, and so on. Here ``order'' does not mean
the order of an ordinary derivative: all operations in~\eqref{eq:lexhom} remain
first-order derivatives, but they are applied to successively
obtained functions of directions. Formula~\eqref{eq:lexhom} does not merely
abbreviate the notation; it gives a more general object. It should therefore
not be presented as the literal definition from the three-page note of 1987:
there the main language was that of the cones~\eqref{eq:lexcone-def}.

\subsection{Lexicographic smoothness}

\begin{definition}[\cite{Nes2005lex}]
\label{def:lsmooth}
A locally Lipschitz mapping $f\colon X\to\R^m$ is called
\emph{lexicographically smooth} at $x\in X$ if, for every
$k\ge1$ and every collection of directions
$M\in\R^{n\times k}$, all functions in the sequence
\eqref{eq:lexhom} exist and are differentiable in every direction.
\end{definition}

The requirement appears strong, yet the class it defines is quite broad.
It contains all continuously differentiable and all finite convex
functions on open sets, is a linear space, and is
closed under componentwise composition of its elements.
It therefore includes finite compositions of smooth operations with
$\max$, $\min$, absolute values, and norms. It was later shown, in particular,
that all piecewise continuously differentiable functions are
lexicographically smooth~\cite{KhanBarton2015}.

The meaning of the definition is revealed by the following finite-dimensional fact.

\begin{theorem}[linearization;~\cite{Nes2005lex}]
\label{th:lexlinear}
Let $f\colon X\subset\R^n\to\R^m$ be lexicographically smooth at $x$, and
let the columns of the matrix $M=[m_1\ \cdots\ m_k]$ span $\R^n$. Then
the last function $f^{(k)}_{x,M}$ in~\eqref{eq:lexhom} is linear in
its argument $d$.
\end{theorem}

The mechanism of the theorem can be seen without technical details. After
differentiation at $m_1$, the function $f^{(1)}_{x,M}$ becomes linear
on the line $\operatorname{span}\{m_1\}$; the next step extends
the subspace of linearity to
$\operatorname{span}\{m_1,m_2\}$, and so on. Once the directions
span the entire space, no nonlinear part can remain.
Finite dimensionality is essential here: in $\R^n$, no more than $n$
linearly independent directions suffice.

Now let $k=n$ and let $M$ be nonsingular. By
Theorem~\ref{th:lexlinear}, there exists a unique matrix
\begin{equation}\label{eq:lexjac}
   J_Lf(x;M)\in\R^{m\times n},\qquad
   f^{(n)}_{x,M}(d)=J_Lf(x;M)d .
\end{equation}
It is called the \emph{lexicographic derivative} (for
a scalar function, the lexicographic gradient) with respect to
$M$. The words ``with respect to $M$'' are essential: for a fixed
ordered basis, the matrix~\eqref{eq:lexjac} is unique, but
replacing or permuting the directions generally changes it.
The collection of all such matrices
\begin{equation}\label{eq:lexsubdiff}
   \partial_L f(x)=
   \bigl\{J_Lf(x;M):M\in\R^{n\times n},\ \det M\ne0\bigr\}
\end{equation}
is called the lexicographic subdifferential.

If $f$ is smooth, then
$\partial_Lf(x)=\{Jf(x)\}$. For a scalar convex function, every
lexicographic gradient is an ordinary subgradient. For
piecewise smooth mappings, the more precise inclusion
\[
   \partial_Lf(x)\subseteq\partial_Bf(x)
   \subseteq\partial^{\mathrm C}f(x),
\]
holds, where $\partial_B$ is the set of limits of Jacobians at nearby
smoothness points, and $\partial^{\mathrm C}$ is its Clarke convex
hull~\cite{KhanBarton2015}. For arbitrary vector-valued
lexicographically smooth mappings, the relationship with the Clarke generalized
Jacobian is more subtle, and this inclusion should not be asserted
without additional assumptions.

\subsection{The lexicographic directional derivative and the exact
chain rule}

Definition~\eqref{eq:lexjac} is convenient for analysis, but computationally
it is simpler to store the results of successive directional derivatives.
Set
\begin{equation}\label{eq:ldder}
   D_Lf(x;M)=
   \bigl[
   f^{(0)}_{x,M}(m_1)\ \
   f^{(1)}_{x,M}(m_2)\ \cdots\
   f^{(k-1)}_{x,M}(m_k)
   \bigr]\in\R^{m\times k}.
\end{equation}
In the later algorithmic literature, this object is called
the lexicographic directional derivative
(LD-derivative)~\cite{KhanBarton2015}. If $k=n$ and $M$ is nonsingular,
then
\begin{equation}\label{eq:ldfactor}
       D_Lf(x;M)=J_Lf(x;M)M,\qquad
       J_Lf(x;M)=D_Lf(x;M)M^{-1}.
\end{equation}
Thus, the sequence of directional derivatives
indeed recovers the linear operator.

\begin{theorem}[exact chain rule;
{\cite{Nes2005lex,KhanBarton2015}}]
\label{th:lexchain}
Let $h\colon X\subset\R^n\to\R^p$ be lexicographically smooth at $x$, and
let $g\colon Y\subset\R^p\to\R^q$ be lexicographically smooth at $h(x)$.
Then $g\circ h$ is lexicographically smooth at $x$, and for every collection
of directions $M\in\R^{n\times k}$,
\begin{equation}\label{eq:lexchain}
   D_L(g\circ h)(x;M)=
   D_Lg\bigl(h(x);D_Lh(x;M)\bigr).
\end{equation}
\end{theorem}

Unlike Clarke's rule, formula~\eqref{eq:lexchain} is an
\emph{equality}. The directions selected by the inner mapping $h$
are passed to the outer function $g$ together and in the same order; therefore,
the correlation between branches is not lost. The proof proceeds
by induction on the number of columns of $M$: at each step, the ordinary
rule for the one-sided directional derivative is applied, after which
one passes to the next homogenization.

The lexicographic derivative also preserves another fundamental
property of the ordinary derivative.

\begin{proposition}[Newton--Leibniz formula;
{\cite{Nes2005lex}}]\label{pr:lexnl}
Let $f\colon X\to\R^m$ be lexicographically smooth in a neighborhood of the segment
$[x,y]\subset X$, $x\ne y$, and let $h=y-x$. Extend $h$ to a basis
$M=[h\ m_2\ \cdots\ m_n]$. Then, componentwise,
\begin{equation}\label{eq:lexnl}
   f(y)-f(x)=
   \int_0^1 J_Lf(x+t h;M)\,h\,dt .
\end{equation}
\end{proposition}

\begin{proof}
The function $\phi(t)=f(x+th)$ is Lipschitz and hence absolutely continuous.
From~\eqref{eq:ldder}--\eqref{eq:ldfactor}, since the first column of
$M$ is $h$, we have
$J_Lf(z;M)h=f'(z;h)$. Consequently, almost everywhere,
$\phi'(t)=J_Lf(x+th;M)h$, and~\eqref{eq:lexnl} follows by integration
with respect to $t$.
\end{proof}

In other words, the selected lexicographic Jacobians are not
arbitrary local slopes: with a consistent choice of
directions, they recover the finite increment of the function.

\subsection{Rules on the computational graph}

In the 1987 note, the rules were formulated directly in
the language of cones~\cite{NesCherLex1987}. Let $u$ and $v$ have
generalized gradients $g_u$ and $g_v$ along a lexicographic cone
constructed from the same ordered basis
$S=(s_1,\dots,s_n)$. Each ``function--gradient'' pair is assigned
a row
\begin{equation}\label{eq:lexrow}
  Y_S(u,g_u)=
  \bigl(u(x),\dotp{s_1}{g_u},\dots,\dotp{s_n}{g_u}\bigr).
\end{equation}
Such rows are compared lexicographically: first the function values,
then the responses of the gradients to $s_1$, then to $s_2$, and so on.

For a smooth outer operation
$w(x)=\varphi(u(x),v(x))$, the ordinary formula is retained:
\begin{equation}\label{eq:lexsmoothop}
   g_w=\varphi_1(u(x),v(x))\,g_u+
       \varphi_2(u(x),v(x))\,g_v .
\end{equation}
In particular, $g_{u+v}=g_u+g_v$ and
$g_{uv}=v(x)g_u+u(x)g_v$. For nonsmooth operations, rather than transforming
sets, one consistently selected row is used:
\begin{equation}\label{eq:lexmax}
 g_{\max\{u,v\}}=
 \begin{cases}
   g_u,&Y_S(u,g_u)\succ_{\rm lex}Y_S(v,g_v),\\
   g_v,&Y_S(u,g_u)\prec_{\rm lex}Y_S(v,g_v),\\
   g_u,&Y_S(u,g_u)=Y_S(v,g_v).
 \end{cases}
\end{equation}
For a minimum, the lexicographically smaller row is selected. Finally,
\begin{equation}\label{eq:lexabs}
 g_{|u|}=
 \begin{cases}
    g_u,&Y_S(u,g_u)\succ_{\rm lex}0,\\
   -g_u,&Y_S(u,g_u)\prec_{\rm lex}0,\\
    0,&Y_S(u,g_u)=0 .
 \end{cases}
\end{equation}
When the two rows in~\eqref{eq:lexmax} coincide completely, we have
$g_u=g_v$, because $S$ is a basis, and the choice is again unique.

The substantive part of the result is not merely the plausibility
of the formulas. Nesterov and Cherkassky proved that after each such operation
there exists a positive cone width on which
the computed vector actually gives the approximation
\eqref{eq:lexcone}. In general, the width must be reduced;
the ordered basis remains unchanged. Thus a finite
superposition of smooth operations, absolute values, maxima, and minima
can be differentiated without propagating multidimensional sets through the graph.

In modern notation, the row~\eqref{eq:lexrow} is replaced by
\[
   \widehat D_Lu(x;M)=\bigl(u(x),D_Lu(x;M)\bigr),
\]
and formulas~\eqref{eq:lexsmoothop}--\eqref{eq:lexabs} are special
cases of the exact chain rule~\eqref{eq:lexchain}.

Now suppose that the program contains $T$ elementary scalar operations.
In the modern vector forward mode, each intermediate variable
carries its value and $k$ directional coefficients; the smooth rule
\eqref{eq:lexsmoothop} and lexicographic comparisons
\eqref{eq:lexmax}--\eqref{eq:lexabs} require $\Oh(k)$ arithmetic
operations per node. Thus the vector forward variant computes
$D_Lf(x;M)$ in $\Oh(kT)$ operations, without enumerating active regions or
constructing multidimensional subdifferentials. For $k=n$,~\eqref{eq:ldfactor}
yields the full selected Jacobian; for $M=I$, no additional
linear system needs to be solved.

\begin{remark}[on the cost of nonsmooth differentiation]
\label{rem:lexcost}
The bound $\Oh(kT)$ should not be confused with the ``cheap gradient''
of Theorem~\ref{th:ad}. There, a single reverse pass gives the entire smooth gradient
in $\Oh(T)$ operations. Here, the forward algorithm simultaneously propagates
$k$ levels of lexicographic information. For piecewise smooth
programs, once the branches have been chosen consistently, more economical
combinations of forward and reverse modes are possible, but this is an additional
algorithmic result, not an automatic consequence of
Definition~\ref{def:lsmooth}. The modern vector forward mode and its
implementation are developed in detail in~\cite{KhanBarton2015}. In the
1987 note, the authors asserted only that the complexity of the new rules was comparable
to that of the formulas of classical differential calculus, without
assigning them a specific asymptotic bound.
\end{remark}

\subsection{An example}

\begin{example}\label{ex:lex}
Consider the piecewise linear function
\[
    F(x_1,x_2)=\max\{x_1,x_2\}+|x_1-x_2|
\]
at $x=0$. In the two smooth regions,
\[
  \grad F(x)=
  \begin{cases}
     (2,-1),&x_1>x_2,\\
     (-1,2),&x_2>x_1.
  \end{cases}
\]
Choose $m_1=(1,1)^\top$ as the first direction. It lies exactly along
the boundary of the two regions: both branches of the maximum are equal, and the argument
of the absolute value is zero. Thus the first level does not select a side, and
\[
    F'(0;m_1)=1.
\]

Let the second direction be $m_2^{+}=(1,-1)^\top$, and let
$M_{+}=[m_1\ m_2^{+}]$. For small $\tau>0$,
\[
  \max\{1+\tau,1-\tau\}+|(1+\tau)-(1-\tau)|
  =1+3\tau,
\]
whence
\[
     D_LF(0;M_{+})=[\,1\ \ 3\,].
\]
By~\eqref{eq:ldfactor}, the lexicographic gradient $g_{+}$ is determined
by $g_{+}^{\top}M_{+}=[\,1\ \ 3\,]$, that is,
\[
     g_{+}=(2,-1)^\top .
\]
If instead we take $m_2^{-}=(-1,1)^\top$ and
$M_{-}=[m_1\ m_2^{-}]$, the second level selects the other region, and
\[
     g_{-}:=J_LF(0;M_{-})^\top=(-1,2)^\top .
\]

Thus the first direction left a tie, and the second broke
it---literally in lexicographic order. At the same time,
\[
 \partial_BF(0)=\{(2,-1),(-1,2)\},\qquad
 \partial^{\mathrm C}F(0)=
 \operatorname{conv}\{(2,-1),(-1,2)\}.
\]
Lexicographic differentiation returns one of the endpoints rather than the whole
segment, and the choice is fully reproducible once
$M$ is fixed.
\end{example}

This last qualification also matters for optimality conditions.
The convex hull sometimes conceals a descent direction. For example,
$r(t)=-|t|$ has
$0\in\partial^{\mathrm C}r(0)=[-1,1]$, although zero is a maximum and
$r'(0;d)=-|d|<0$ for every $d\ne0$. The lexicographic derivatives
themselves also form the set $\{-1,1\}$, whose convex hull
contains zero. The stronger information lies in the
\emph{pair} ``ordered basis--selected derivative'':
\[
       J_Lr(0;M)m_1=r'(0;m_1)<0 .
\]
It is the preservation of this dependence, rather than a simple replacement of one
subdifferential by another, that strengthens the screening of candidates for a minimum
in~\cite{Nes2005lex}.

\subsection{What changed and the current state}

\smallskip
\textbf{(i) Exact calculus in place of a calculus of inclusions.}
The lexicographic derivative showed that, for a broad class of
nonsmooth functions, the central features of classical
analysis can be retained: a linear Jacobian after finitely many steps, an exact chain
rule~\eqref{eq:lexchain}, and the
Newton--Leibniz formula~\eqref{eq:lexnl}. The price is additional
structure---an ordered collection of directions; without it,
the lexicographic gradient is not an intrinsically unique object.

\smallskip
\textbf{(ii) Another way to look inside the black box.}
An ordinary nonsmooth-function oracle returns a value and possibly a
subgradient, but does not indicate how choices should be coordinated under composition.
If the function is specified by a computational graph, however, rules
\eqref{eq:lexsmoothop}--\eqref{eq:lexabs} use its internal
structure and propagate the same priority system through all
nodes. As in \S\,\ref{sec:ad}, the advantage comes not from a new
lower bound in the old information model, but from a richer
model of the function.

\smallskip
\textbf{(iii) Algorithmic developments.}
K.\,A.~Khan and P.\,I.~Barton showed that, for piecewise smooth programs,
lexicographic derivatives are elements of the
$B$-subdifferential, and constructed a practically implementable vector
forward mode of automatic differentiation
\cite{KhanBarton2015}. The same construction was extended to the sensitivity
analysis of parametric nonlinear programs when the active set of
constraints changes~\cite{StechlinskiKhanBarton2018}.

\smallskip
\textbf{(iv) Limits of applicability.}
Lexicographic smoothness covers smooth, convex, and many
piecewise smooth compositions, but not all locally Lipschitz functions.
Moreover, an arbitrary result returned by a modern
library at a $\max$ or $\mathrm{ReLU}$ node does not automatically become
a lexicographic derivative: branch-selection rules
must be consistent throughout the graph (we note that Yu.\,E.~Nesterov pointed out in personal conversations that the justification of backpropagation for neural networks with the ReLU function had, on the whole, not been worked out very carefully mathematically in the AI literature). A more general answer to the question
``what nonsmooth backpropagation computes'' is provided by
conservative fields~\cite{BoltePauwels2021}. These two theories do not
compete: lexicographic differentiation gives a sharper
object and an exact calculus on a special structured class, while
conservative fields provide a broader model sufficient for analyzing
algorithmic trajectories.

Historically, this line of work has remained somewhat in the shadow of accelerated methods and
interior-point methods. Yet its formulation of the question is just as
characteristic of Yu.\,E.~Nesterov's work: instead of passing an algorithm
a set-valued answer from a nonsmooth oracle, it proposes opening up the structure
of the function and making the choice of derivative part of the computation.

\section{Self-concordant barriers and interior-point methods
(1988 -- 1994)}
\label{sec:ipm}

\subsection{Background: why the logarithm?}

By the mid-1980s, a paradoxical situation had developed in the theory of linear
programming (LP). Dantzig's simplex method, known since
1947, worked impeccably in practice, but in
1972 V.~Klee and G.~Minty~\cite{KleeMinty1972} constructed examples on which it
enumerates exponentially many vertices. Conversely, the ellipsoid method,
applied to LP by L.\,G.~Khachiyan~\cite{Khachiyan1979} in 1979, gave
the \emph{first} proof of polynomial-time solvability of linear
programming---a theoretical breakthrough of the first order---but
proved impractical: it does not exploit the structure of the problem and, on
real examples, is orders of magnitude slower than the simplex method.

In 1984, N.~Karmarkar~\cite{Karmarkar1984} proposed a projective
method for LP that was both polynomial and competitive
in practice. It soon became clear that his method was related to
the long-known idea of \emph{barriers}: instead of a constrained problem
$x\in Q$, solve a sequence of smooth unconstrained problems
\[
   \min_{x\in\operatorname{int}Q}\ \bigl\{\, t\,\dotp{c}{x} + F(x) \,\bigr\},\qquad t\to+\infty,
\]
where $F$ is a function tending to $+\infty$ at the boundary of $Q$. For
$Q=\{x:\ Ax\le b\}$, the natural choice is the logarithmic barrier
$F(x)=-\sum_j\ln(b_j-\dotp{a_j}{x})$, proposed as early as R.~Frisch (1955)
and studied by A.~Fiacco and G.~McCormick (1968); J.~Renegar
(1988)~\cite{Renegar1988} showed that following the ``central path''
of this function yields a bound of $\Oh(\sqrt{m}\,\ln(1/\eps))$ Newton
iterations for LP.

A fundamental question nevertheless remained. In
general, Newton's method does \emph{not} have global guarantees: its convergence
rate depends on Lipschitz constants, which for barrier functions
grow without bound near the boundary. Why, then, does everything work for the
logarithmic barrier? Is there something special about the logarithm---or is there
an entire class of ``good'' barriers, and if so, how can it be described?

The answer given by Yu.\,E.~Nesterov and A.\,S.~Nemirovski in a series of papers
in the late 1980s and in the monograph~\cite{NesNem1994} is
as follows. The ``goodness'' of a barrier is not a property of the logarithm, but
of the \emph{concordance of the third derivative with the second}. It is
this property that makes Newton's method an affine-invariant tool with
a global complexity bound, and it is this property that makes it possible to construct
polynomial-time methods not only for LP, but for any convex problem
admitting a computable barrier with a small parameter.

\subsection{Definitions}

Throughout this section, $D\subset\R^n$ is an \emph{open proper
convex domain}: an open nonempty convex set, $D\ne\R^n$, containing
no lines; put $Q=\overline D$. The requirement that it contain no
lines is essential not for the definitions (positive definiteness
of the Hessian is stipulated separately), but for the \emph{existence} of a barrier
of the type considered and, in particular, for the universal-barrier
formula~\eqref{eq:univbar}: if $Q$ contains a line, the polar set
lies in a proper subspace and has zero volume.

The presence of lines in $D$ is not a severe practical restriction. In that case, the barrier is defined on the quotient of $D$ by the linear subspace generated by the lines contained in $D$. Thus the resulting barrier function is constant along these lines.

\begin{definition}[self-concordance]\label{def:sc}
A three-times continuously differentiable convex function
$F\colon D\to\R$ with $\grad^2F(x)\succ0$ for every $x\in D$ is called
\emph{self-concordant} if
\begin{equation}\label{eq:sc}
   \bigl|D^3F(x)[h,h,h]\bigr|\ \le\ 2\,\bigl(D^2F(x)[h,h]\bigr)^{3/2}
   \qquad\forall x\in D,\ \forall h\in\R^n .
\end{equation}
If, in addition, $F(x_k)\to+\infty$ for every sequence
$x_k\in D$, $x_k\to x\in\partial D$, then $F$ is called
a \emph{self-concordant barrier} for $Q$.
\end{definition}

Both conventions occur in the literature: the
barrier property may be included in the definition of self-concordance or
stated separately. We adopt the latter convention because different
classes of objects are compared below.

Membership in $C^3$ is assumed for convenience. Methods using self-concordant barriers do not require a third-order oracle. Their theory also works for a limit of a sequence of functions self-concordant in the sense of Definition~\ref{def:sc}, even though such a limit need not be three-times continuously differentiable.

The constant~$2$ in~\eqref{eq:sc} is a normalization: the replacement $F\mapsto
\alpha F$ for $\alpha>1$ only improves the inequality, so any
function satisfying $|D^3F[h,h,h]|\le M(D^2F[h,h])^{3/2}$ can be brought
to the form~\eqref{eq:sc}. The essential point is different: inequality~\eqref{eq:sc}
is \emph{affine-invariant}---both sides transform in the same way under
a change of variables $x=By$. Methods built on it will therefore
also be affine-invariant: they do not depend on the coordinate system.
This is a fundamental difference from gradient methods, where everything is determined
by a Euclidean norm chosen in advance.

Associate with the second derivative of the barrier a local Euclidean structure
(which varies with $x$!):
\[
   \nrm{h}_x=\dotp{\grad^2F(x)h}{h}^{1/2},\qquad
   \nrm{g}_{x}^{*}=\dotp{[\grad^2F(x)]^{-1}g}{g}^{1/2}.
\]
An important geometric consequence of self-concordance of a \emph{barrier} is that
the \emph{Dikin ellipsoid}
\[
   W(x)=\{y:\ \nrm{y-x}_x<1\}
\]
lies entirely in $D$ (and its closure is contained in $Q$).
The barrier qualification is indispensable here: inequality
\eqref{eq:sc} alone is not enough. Indeed, for
$D=(-\frac12,\frac12)$ and $F(x)=\frac{x^2}{2}$, we have $F''\equiv1>0$,
$F'''\equiv0$, so~\eqref{eq:sc} holds, but the Dikin ellipsoid at
zero is the interval $(-1,1)\not\subset D$; the reason is that the domain
has been artificially truncated and $F$ does not tend to $+\infty$ at its boundary.
For a self-concordant \emph{barrier}, the inclusion is valid, and its meaning is this:
the barrier Hessian itself ``senses'' the shape of the set---at a point close to
the boundary, the Dikin ellipsoid flattens in the corresponding direction
(Fig.~\ref{fig:cpath}).

\begin{definition}[$\nu$-self-concordant barrier]\label{def:scb}
A self-concordant barrier $F$ for $Q$ is called
\emph{$\nu$-self-concordant} if, in addition,
\begin{equation}\label{eq:nu}
   \dotp{\grad F(x)}{h}^2\ \le\ \nu\,\dotp{\grad^2F(x)h}{h}
   \qquad\forall x\in D,\ \forall h\in\R^n,
\end{equation}
that is, $\nrm{\grad F(x)}^{*}_x\le\sqrt{\nu}$. The number $\nu\ge1$
is called the barrier parameter.
\end{definition}

The meaning of~\eqref{eq:nu} is that the barrier gradient, measured in the local
metric induced by its own Hessian, is bounded. It is the parameter $\nu$
(rather than the dimension, the number of constraints, or Lipschitz constants) that enters
the final complexity bound. The problem dimension does, however, indirectly affect the parameter: a barrier with a small parameter may not exist.

\begin{example}\label{ex:barriers}
Direct verification gives the catalogue in Table~\ref{tab:barr}. For $n=1$,
$D=(0,+\infty)$, $F(x)=-\ln x$, we have $F''=x^{-2}$, $F'''=-2x^{-3}$,
so $|F'''|=2(F'')^{3/2}$---\emph{equality} is attained in~\eqref{eq:sc};
further, $(F')^2=x^{-2}=F''$, that is, $\nu=1$. Both bounds
are additive under direct sums, giving the rows for the orthant and
polyhedron in the table.
\end{example}

\begin{table}[htbp]
\centering
\caption{Standard self-concordant barriers. Throughout, $D$ is the interior
of the indicated set.}
\label{tab:barr}
\footnotesize
\renewcommand{\arraystretch}{1.2}
\begin{tabular}{@{}lll@{}}
\toprule
Set $Q$ & Barrier $F$ & $\nu$\\
\midrule
$\R^n_+$ & $-\sum_{i=1}^n\ln x_i$ & $n$\\
$\{x:\ \dotp{a_j}{x}\le b_j,\ j=1,\dots,m\}$
   & $-\sum_{j=1}^m\ln\bigl(b_j-\dotp{a_j}{x}\bigr)$ & $m$\\
Lorentz cone $\{(x,\tau):\ \nrm{x}_2\le\tau\}$
   & $-\ln\bigl(\tau^2-\nrm{x}_2^2\bigr)$ & $2$\\
Cone $\Ss^n_+$ of positive semidefinite matrices
   & $-\ln\det X$ & $n$\\
Exponential cone $\{y\,e^{x/y}\le z,\ y>0\}$
   & $-\ln\bigl(y\ln(z/y)-x\bigr)-\ln y-\ln z$ & $3$\\
\emph{Any} proper convex domain in $\R^n$
   & Universal barrier \eqref{eq:univbar} & $n$\\
\bottomrule
\end{tabular}
\end{table}

The last row of the table is a separate result of
Nesterov and Nemirovski, showing that the theory is nonvacuous for every
proper convex domain. More precisely, the classical Nesterov--Nemirovski result states that a barrier with parameter $\nu = \Oh(n)$ always exists. The \emph{universal barrier}
is the function
\begin{equation}\label{eq:univbar}
   \varphi(x)\ =\ \ln\vol\bigl(P(x)\bigr),\qquad
   P(x)=\{s\in\R^n:\ \dotp{s}{y-x}\le1\ \ \forall y\in Q\},
\end{equation}
where $P(x)$ is the polar set of $Q-x$. Yu.\,E.~Nesterov and
A.\,S.~Nemirovski proved that $\varphi$ (up to multiplication by
an absolute constant needed to obtain the normalization in
Definition~\ref{def:sc}) is a $\nu$-self-concordant barrier with
$\nu=\Oh(n)$. Multiplying a barrier by a constant $C>1$
also scales its parameter, so a formula with an unspecified factor
is incompatible with an exact value of $\nu$; the exact statements below
use normalization~\eqref{eq:univbar}.
Determining the exact constant took more than twenty years. First, S.~Bubeck
and R.~Eldan~\cite{BubeckEldan2019} constructed the \emph{entropic
barrier}---the logarithm of the Laplace transform of the indicator of the set---and
proved that it is $\nu$-self-concordant with $\nu=(1+o(1))n$;
thus an asymptotically optimal constant was achieved by an explicitly
specified function that also has a meaningful
probabilistic interpretation (the barrier Hessian is the covariance matrix
of the corresponding exponential family). Then Y.\,T.~Lee and
M.-C.~Yue showed that the universal barrier~\eqref{eq:univbar} itself,
with precisely this normalization and no additional factor, is
$n$-self-concordant~\cite{LeeYue2021}, so the exact upper
bound on the parameter is $n$.

\begin{remark}[existence versus computability]\label{rem:univcomp}
The statement ``the barrier parameter can always be made of the order of the dimension''
concerns the \emph{existence} of a barrier and says nothing about its
\emph{computability}. Computing the volume of the polar set
in~\eqref{eq:univbar}, let alone its gradient and Hessian, is in general
no easier than the original problem. Thus the universal-barrier theorem
shows that the theory is nonvacuous, but does not automatically give a polynomial-time
algorithm in a particular input model. The practical value
of the theory lies in specific barriers for which $F$, $\grad F$, and $\grad^2F$ can be computed
explicitly and cheaply. Such efficiently computable barriers include those in the first five rows
of Table~\ref{tab:barr}. In recent years, the list of sets for which an efficiently computable barrier is known has grown further~\cite{RoyXiao22}.
\end{remark}

The requirement that the domain be proper cannot be omitted either. If $Q$ contains
a line, then the polar set $P(x)$ lies in a proper subspace
and has zero volume, so formula~\eqref{eq:univbar} does not define
a finite barrier. The standard remedy is to pass to the quotient space by
the maximal linear subspace contained in $Q$.

\subsection{The central path and the path-following method}

\begin{figure}[h]
\centering
\includegraphics[width=0.70\textwidth]{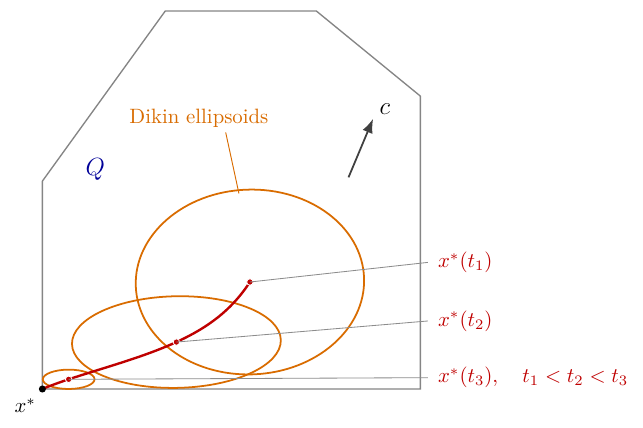}
\caption{The central path for
$\min\{\dotp{c}{x}:x\in Q\}$ on a hexagon $Q$ with the
logarithmic barrier
$F(x)=-\sum_j\ln\bigl(b_j-\dotp{a_j}{x}\bigr)$.
The arrow $c$ indicates the direction of increase of the linear objective;
the minimum $x^{*}$ is attained at the opposite boundary vertex.
The red curve consists of the points
$x^{*}(t)=\argmin_{x\in\operatorname{int}Q}
\{t\dotp{c}{x}+F(x)\}$ and, as $t$ increases, approaches $x^{*}$.
Orange indicates the closures of the unit Dikin ellipsoids
$\{y:\nrm{y-x^{*}(t_i)}_{x^{*}(t_i)}\le1\}$ at three points on the path,
where $t_1<t_2<t_3$. Their open parts lie in
$\operatorname{int}Q$ and shrink in ordinary coordinates as
the boundary is approached.}
\label{fig:cpath}
\end{figure}

Consider a conic (or simply convex) programming problem in
the form
\begin{equation}\label{eq:lincvx}
   \min_{x\in Q}\ \dotp{c}{x},
\end{equation}
to which every convex optimization problem can be reduced (by moving
the objective into the constraints through its epigraph). Set
\[
   f_t(x)=t\,\dotp{c}{x}+F(x),\qquad
   x^{*}(t)=\argmin_{x\in D} f_t(x),\qquad t>0 .
\]
Here $F$ is a self-concordant barrier on $D = Q^{\circ}$. The curve $\{x^{*}(t)\}_{t>0}$ is called the \emph{central path}; as
$t\to\infty$, it tends to a solution of~\eqref{eq:lincvx}, if one exists
(Fig.~\ref{fig:cpath}). As $t\to0^{+}$, it emanates from
the \emph{analytic center} $\argmin_{x\in D}F(x)$---but only provided
that this minimum exists. For unbounded domains, and in particular for
cones, it does not exist: the function $F(x)=-\ln x$ on $(0,+\infty)$
has no minimum. In general, small values of $t$ should be understood not as
the beginning of the path, but as a chosen initial central neighborhood
(see Remark~\ref{rem:phase1}).

Accuracy along the path has an elementary bound.

\begin{lemma}\label{lem:cpath}
Let $F$ be a $\nu$-self-concordant barrier for $Q$. Then, for every
$t>0$,
\begin{equation}\label{eq:gap}
  \dotp{c}{x^{*}(t)}-\min_{x\in Q}\dotp{c}{x}\ \le\ \frac{\nu}{t}.
\end{equation}
\end{lemma}

\begin{proof}
First we show that $\dotp{\grad F(x)}{y-x}<\nu$ for all
$x,y\in D$. Set $\psi(\tau)=\dotp{\grad F(x+\tau(y-x))}{y-x}$,
$\tau\in[0,1]$. Then $\psi'(\tau)=\nrm{y-x}^2_{x+\tau(y-x)}\ge
\psi(\tau)^2/\nu$ by~\eqref{eq:nu}. If $\psi(0)\le0$, there is nothing
to prove. If $\psi(0)>0$, then $\psi$ is increasing and positive, and integrating
$\nu\psi'/\psi^2\ge1$ gives
$\nu/\psi(0)-\nu/\psi(1)\ge1$, whence $\nu/\psi(0)>1$, that is,
$\psi(0)<\nu$.

Now the point $x=x^{*}(t)$ satisfies $tc+\grad F(x)=0$, so for
every $y\in D$,
\[
   \dotp{c}{x-y}=-\tfrac1t\dotp{\grad F(x)}{x-y}
   =\tfrac1t\dotp{\grad F(x)}{y-x}\le\frac{\nu}{t}.
\]
Passing to the limit in $y$ (the interior is dense in $Q$)
yields~\eqref{eq:gap}.
\end{proof}

Thus, to attain accuracy $\eps$, it suffices to follow the path to $t\approx\nu/\eps$.
The question remains: how many Newton steps are needed to move
along the path? Proximity to the path is measured by the \emph{Newton decrement}
\[
   \lambda(t,x)=\nrm{\grad f_t(x)}^{*}_x
   =\bigl\|[\grad^2F(x)]^{-1/2}\bigl(tc+\grad F(x)\bigr)\bigr\|_2 .
\]

\begin{theorem}[path-following method;~\cite{NesNem1994}]
\label{th:pf}
\emph{Class:} problem~\eqref{eq:lincvx} with a $\nu$-self-concordant barrier
$F$ for $Q$; the minimum $\min_{x\in Q}\dotp{c}{x}$ is assumed to be attained.
\emph{Oracle:} values of $\grad F$ and $\grad^2F$. \emph{Initial
condition:} a pair $(t_0,x_0)$, $t_0>0$, $x_0\in D$, is given with
$\lambda(t_0,x_0)\le\beta$. Then there exist absolute constants
$\beta>0$ and $\gamma>0$ such that the process
\[
   t_{k+1}=t_k\Bigl(1+\frac{\gamma}{\sqrt{\nu}}\Bigr),\qquad
   x_{k+1}=x_k-\bigl[\grad^2F(x_k)\bigr]^{-1}\grad f_{t_{k+1}}(x_k)
\]
preserves $\lambda(t_k,x_k)\le\beta$ for every $k$, and after
\[
   N\ =\ \Oh\!\left(\sqrt{\nu}\,
   \ln\frac{\nu}{t_0\,\eps}\right)
\]
iterations (when $t_0\eps<\nu$), a point $x_N\in D$ with
$\dotp{c}{x_N}-\min_{x\in Q}\dotp{c}{x}\le\eps$ is obtained. Each iteration consists of
one evaluation of the barrier gradient and Hessian and the solution of one linear
system.
\end{theorem}

\begin{remark}[on the starting point; phase~I]\label{rem:phase1}
The initial condition $\lambda(t_0,x_0)\le\beta$ does not come for free: finding
such a pair is a separate problem, called phase~I. Standard
approaches include minimizing the barrier itself (if the analytic center
exists), solving an auxiliary problem with an artificial variable,
and homogenization (passing to a self-dual embedding). All of them
fit within the same bound $\Oh(\sqrt\nu\,\ln(\cdot))$, so the final
complexity order is unchanged, but this stage cannot be ignored when
stating the theorem.
\end{remark}

The origin of the factor $\sqrt{\nu}$ is clear from
inequality~\eqref{eq:nu}: increasing $t$ by a factor of
$(1+\gamma/\sqrt{\nu})$ changes the gradient $\grad f_t$ at the current point
by an amount of order
$\gamma\nrm{\grad F(x)}^{*}_x/\sqrt{\nu}\le\gamma$ in the local metric,
so the perturbation remains within the Dikin ellipsoid---the region where
Newton's method converges quadratically and one step suffices to
restore proximity to the path. Thus the step size in the parameter
$t$ is controlled \emph{automatically} by the geometry of the barrier. The total number
of iterations $\sqrt{\nu}\,\ln(\cdot)$ is essentially the ``number of
Dikin ellipsoids that fit along the central path on a
logarithmic scale.''

The geometry of the central path and the change in the local barrier metric
are shown in Fig.~\ref{fig:cpath}.

\begin{example}\label{ex:lp}
For the LP problem $\min\{\dotp{c}{x}:\ Ax\le b\}$ with $A\in\R^{m\times n}$,
we have $\nu=m$, so $\Oh(\sqrt{m}\,\ln(\cdot))$ Newton steps are needed,
each assembling the matrix $A^{\!\top}\!DA$ ($\Oh(mn^2)$ operations) and
computing its Cholesky factorization ($\Oh(n^3)$). For semidefinite
programming $\min\{\dotp{C}{X}:\ \dotp{A_i}{X}=b_i,\ X\succeq0\}$,
$X\in\Ss^n$, the barrier $-\ln\det X$ gives $\nu=n$ and a bound of
$\Oh(\sqrt{n}\,\ln(\cdot))$ iterations, even though the dimension
of the variable is $n(n+1)/2$. It is this ``square-root'' saving that
makes the method practical.
\end{example}

\subsection{What changed}

\textbf{(i) The boundary of computability shifted.} Before 1988, ``efficiently
solvable'' meant ``linear'' in practice. The theory of self-concordant
barriers replaced this criterion by another: a problem is efficiently solvable
if its feasible set has a known barrier with a small parameter
$\nu$ whose value, gradient, and Hessian are computable in
polynomial time. This immediately brought semidefinite programming,
geometric programming, eigenvalue-constrained problems, and second-order
cone constraints into the ``efficiently solvable'' class.

\smallskip
\textbf{(ii) Conic programming emerged as a language.} Nesterov and
Todd's work on \emph{self-scaled} barriers~\cite{NesTodd1997}
identified the class of cones for which primal-dual
symmetric interior-point methods exist. It later became clear
(O.~G\"uler~\cite{HauserGuler02}; the classification uses the theory of Euclidean
Jordan algebras, see the monograph by J.~Faraut and
A.~Kor\'anyi~\cite{FarautKoranyi1994}) that self-scaled cones (that is, cones admitting a self-scaled barrier) are
exactly the \emph{symmetric} (homogeneous self-dual) cones:
the nonnegative orthant, Lorentz cones, cones of Hermitian positive
semidefinite matrices of arbitrary size over $\R$, $\mathbb{C}$, and the quaternions, and of size 3 over the octonions, as well as their direct products.
The octonionic cone is also called the Albert cone; it is an exceptional 27-dimensional cone that does not occur in applications. This explains why
the triad LP\,--\,SOCP\,--\,SDP became fundamental in the practice of convex optimization:
up to standard real representations, these three cones cover all classical families
of symmetric cones.

\smallskip
\textbf{(iii) Semidefinite programming became a working
tool.} The ability to solve medium-sized SDP problems numerically
transformed control theory (the Lyapunov-function method through linear
matrix inequalities), approximation theory in combinatorial optimization
(Goemans--Williamson relaxations~\cite{GoemansWilliamson1995},
Lasserre and Parrilo hierarchies, sums of squares), quantum information,
and statistics (matrix recovery problems). The software packages SeDuMi,
SDPT3, MOSEK, and the ``disciplined'' convex
programming interfaces CVX/CVXPY implement precisely this theory.

\smallskip
\textbf{(iv) Current state.} The bound $\Oh(\sqrt{\nu}\,\ln(\cdot))$
remains a benchmark, but it has been improved through a more refined treatment of
the path: Y.\,T.~Lee and A.~Sidford~\cite{LeeSidford2014} constructed an LP
method requiring $\Ot(\sqrt{n}\,\ln(1/\eps))$ iterations (in terms of the number of
\emph{variables}, rather than constraints), and M.~Cohen, Y.\,T.~Lee, and
Z.~Song~\cite{CohenLeeSong2021} obtained a solution of an LP problem in time
$\Ot(n^{\omega})$, that is, up to logarithmic factors, the time needed for matrix
multiplication. These results in turn led to breakthroughs in network
flow algorithms. In parallel, there is active development of methods for
\emph{nonsymmetric} cones (power cones, exponential cones, relative
entropy cones), where the Nesterov--Todd primal-dual theory
does not apply, but barrier calculus works; the packages
\texttt{alfonso} and \texttt{Hypatia} should be mentioned here. Finally,
barrier techniques have unexpectedly returned to first-order theory:
smoothness of a barrier ``relative to itself'' is an early precursor
of the modern notion of relative smoothness, to be discussed in
\S\,\ref{sec:superfast}.

\section{The \texorpdfstring{$\pi/2$}{pi/2} theorem: accuracy of semidefinite relaxation (1998)}
\label{sec:pi2}

\subsection{Problem statement and background}

Having mastered techniques for solving SDP problems, it is natural to ask: how
good are semidefinite relaxations of \emph{nonconvex} problems? A classical
object is nonconvex quadratic maximization
\begin{equation}\label{eq:bqp}
   f^{*}=\max_{x}\ \bigl\{\dotp{Ax}{x}:\ |x_i|\le1,\ i=1,\dots,n\bigr\},
   \qquad A\in\Ss^n,\quad A\succeq0 .
\end{equation}
Since, for $A\succeq0$, the function $\dotp{Ax}{x}$ is convex, its maximum on
the cube is attained at a vertex, so~\eqref{eq:bqp} is equivalent to
the \emph{Boolean} problem
\[
   f^{*}=\max\bigl\{\dotp{Ax}{x}:\ x\in\{-1,1\}^n\bigr\},
\]
which is NP-hard (it includes, for example, the maximum-cut
problem). In the Boolean problem, the condition $A \succeq 0$ entails no loss of generality, since adding a multiple of the identity matrix to $A$ shifts the objective by an additive constant.

The standard relaxation replaces $x$ by the matrix
$X=xx^{\!\top}$ and drops the condition $\rank X=1$:
\begin{equation}\label{eq:sdprel}
   \mathrm{SDP}=\max\bigl\{\dotp{A}{X}:\ X\succeq0,\ X_{ii}=1,\
   i=1,\dots,n\bigr\}.
\end{equation}
This is a semidefinite programming problem solvable by an
interior-point method. Clearly, $f^{*}\le\mathrm{SDP}$
(the matrix $xx^{\!\top}$ is feasible for~\eqref{eq:sdprel}). The question is:
how large can the gap be?

By the 1990s, it was known that the gap could be controlled for \emph{some} problems:
M.~Goemans and D.~Williamson~\cite{GoemansWilliamson1995} showed that, for
the maximum-cut problem, randomized rounding of an SDP solution gives
a guarantee of $0{,}878\ldots$, with the qualification that this bound concerns the cut weight rather than the objective in~\eqref{eq:sdprel}. Their analysis, however, relied essentially on
nonnegative graph weights and did not extend to an arbitrary matrix
$A$. Yu.\,E.~Nesterov's result~\cite{Nes1998} settled the general
case, giving a \emph{universal and sharp} guarantee for the whole class
of positive semidefinite quadratic forms on the cube.

\begin{theorem}[$\pi/2$ theorem;~\cite{Nes1998}]\label{th:pi2}
Let $A\succeq0$. Then
\begin{equation}\label{eq:pi2}
   \frac{2}{\pi}\,\mathrm{SDP}\ \le\ f^{*}\ \le\ \mathrm{SDP},
\end{equation}
Moreover, if $X^{*}$ is a solution of~\eqref{eq:sdprel} and
$x=\sgn\xi$, where $\xi\sim\mathcal N(0,X^{*})$, then
\begin{equation}\label{eq:roundguar}
   \Exp\bigl[\dotp{Ax}{x}\bigr]\ \ge\ \frac{2}{\pi}\,\mathrm{SDP}.
\end{equation}
\end{theorem}

It is useful to distinguish three quantities obtained from~\eqref{eq:pi2}, since
they are often confused:
\begin{align}
   \text{rounding guarantee ratio:}\qquad
   &\frac{\Exp\dotp{Ax}{x}}{\mathrm{SDP}}\ \ge\ \frac{2}{\pi}
   \approx0{,}6366;\label{eq:pi2a}\\
   \text{integrality gap of the relaxation:}\qquad
   &\frac{\mathrm{SDP}}{f^{*}}\ \le\ \frac{\pi}{2}
   \approx1{,}5708;\label{eq:pi2b}\\
   \text{maximum relative overestimate:}\qquad
   &\frac{\mathrm{SDP}-f^{*}}{f^{*}}\ \le\ \frac{\pi}{2}-1
   \approx0{,}5708.\label{eq:pi2c}
\end{align}
The standard \emph{approximation ratio} is the quantity
$2/\pi$ in~\eqref{eq:pi2a}: it states the fraction of the relaxation value
that rounding is guaranteed to deliver. The quantity $\pi/2-1$
in~\eqref{eq:pi2c} is not ``relative accuracy'' in the usual sense, but
the maximum relative overestimate of the optimum by the relaxation value;
calling it accuracy would promise the reader a quality guarantee of
$57\%$ instead of the correct $63{,}66\%$.

\subsection{Proof}

The proof is short, and worth giving in full: it exemplifies
how convex analysis and a probabilistic argument combine in a single
line.

\begin{proof}
Let $X^{*}$ be a solution of~\eqref{eq:sdprel}. Since $X^{*}\succeq0$ and
$X^{*}_{ii}=1$, the matrix $X^{*}$ is a correlation matrix; take
a Gaussian vector $\xi\sim\mathcal N(0,X^{*})$ and put
$x=\sgn\xi$ (componentwise), so $x\in\{-1,1\}^n$ almost surely.
By Sheppard's formula (Grothendieck's identity), for standard normal
$\xi_i,\xi_j$ with correlation coefficient $\rho=X^{*}_{ij}$,
\[
   \Exp\bigl[\sgn\xi_i\cdot\sgn\xi_j\bigr]
   =\frac{2}{\pi}\arcsin\rho .
\]
Thus, writing $\arcsin[X]$ for the matrix with entries
$\arcsin X_{ij}$,
\begin{equation}\label{eq:round}
   \Exp\bigl[\dotp{Ax}{x}\bigr]
   =\frac{2}{\pi}\,\dotp{A}{\arcsin[X^{*}]}
   =\frac{2}{\pi}\,\dotp{A}{X^{*}}
    +\frac{2}{\pi}\,\dotp{A}{\arcsin[X^{*}]-X^{*}} .
\end{equation}
The key observation is that, for every correlation matrix $X$,
\begin{equation}\label{eq:arcsinpsd}
   \arcsin[X]-X\ \succeq\ 0 .
\end{equation}
Indeed, the Taylor series $\arcsin u=\sum_{k\ge0}c_k u^{2k+1}$ with
$c_k=\dfrac{(2k)!}{4^k(k!)^2(2k+1)}>0$ converges for $|u|\le1$, and
$|X_{ij}|\le\sqrt{X_{ii}X_{jj}}=1$, so
\[
   \arcsin[X]-X=\sum_{k\ge1}c_k\,X^{\circ(2k+1)},
\]
where $X^{\circ m}$ is the $m$th Hadamard (entrywise) power. By the
Schur product theorem, the Hadamard product of positive semidefinite
matrices is positive semidefinite; hence $X^{\circ m}\succeq0$ for every
$m\ge1$, and the entire sum is positive semidefinite. Finally, for
$A\succeq0$ and $B\succeq0$, we have $\dotp{A}{B}=\Tr(AB)\ge0$, so
the second term in~\eqref{eq:round} is nonnegative. This gives
\eqref{eq:roundguar}, and with it
\eqref{eq:pi2}, since $f^{*}\ge\Exp[\dotp{Ax}{x}]$.
\end{proof}

Observe exactly where the condition $A\succeq0$ was used: only at the very
end, to discard the term $\dotp{A}{\arcsin[X^{*}]-X^{*}}$.
This is why the theorem is so robust: no assumptions about the signs
of the entries of $A$, the graph structure, or sparsity are required.

\subsection{From a bound in expectation to an algorithmic guarantee}

Bound~\eqref{eq:roundguar} is a guarantee \emph{in expectation} for a single
random trial. It immediately implies the existence of a suitable
sign vector, but does \emph{not} automatically imply that a
polynomial number of independent trials will attain the exact threshold
$2/\pi$ with high probability. The correct statement is as follows.

\begin{proposition}\label{pr:round}
Set $Z=\dotp{Ax}{x}/\mathrm{SDP}$; then $Z\in[0,1]$ and
$\Exp Z\ge a:=2/\pi$. For every $\delta\in(0,a)$,
\begin{equation}\label{eq:markov}
   \Prob\{Z\ge a-\delta\}\ \ge\ \frac{\delta}{1-a+\delta}.
\end{equation}
Consequently, $\Oh\bigl(\frac{1-a+\delta}{\delta}\ln\frac1\gamma\bigr)$
independent trials produce, with probability at least $1-\gamma$,
a vector $x\in\{-1,1\}^n$ with $\dotp{Ax}{x}\ge(2/\pi-\delta)\,\mathrm{SDP}$.
\end{proposition}

\begin{proof}
Let $p=\Prob\{Z\ge a-\delta\}$. Since $Z\le1$,
$\Exp Z\le(a-\delta)(1-p)+p=a-\delta+p(1-a+\delta)$. Together with
$\Exp Z\ge a$, this gives $p\ge\delta/(1-a+\delta)$.
\end{proof}

Thus any factor $2/\pi-\delta$ can be attained in polynomial time;
for the \emph{exact} factor $2/\pi$, one should either
retain the ``in expectation'' formulation or use deterministic
rounding by the method of conditional expectations.

\subsection{Sharpness of the constant and an example}

The constant $2/\pi$ cannot be improved. Two different
statements must be distinguished here.

\smallskip
\emph{Integrality gap of the relaxation.} The worst-case ratio
$\mathrm{SDP}/f^{*}$ over all $A\succeq0$ is exactly the quantity
known as the \emph{positive Grothendieck constant} $K_G^{+}$, and
it equals $\pi/2$. The key bound goes back to
R.~Rietz's work~\cite{Rietz1974} on the proof of Grothendieck's
inequality; the equality $K_G^{+}=\pi/2$ in precisely the
normalization used in~\eqref{eq:bqp}\,--\,\eqref{eq:sdprel},
together with the history of the question, is given in the survey by
Khot and Naor~\cite[\S\,2]{KhotNaor2012}. In other words, no analysis
of \emph{this} relaxation can yield a constant better than $2/\pi$.

\smallskip
\emph{Hardness of approximation.} A separate question is whether
a better factor can be obtained by \emph{any} polynomial-time algorithm.
Assuming the Unique Games Conjecture, known hardness-of-approximation
results show that it cannot; for a survey of these
results with precise statements, see
\cite{KhotNaor2012}.

\smallskip
Finally, if the condition $A\succeq0$ is dropped, the quality of the relaxation
deteriorates sharply: the gap becomes of order
$\Theta(\ln n)$~\cite{NemRoosTerlaky1999}.

It is useful to see that the gap in~\eqref{eq:pi2} already occurs in a very simple
geometric setting.

\begin{example}[vectors on a semicircle]\label{ex:semicircle}
Let $\theta_i=\pi(i-1)/n$, $u_i=(\cos\theta_i,\sin\theta_i)\in\R^2$,
and $A=UU^{\!\top}$, where the rows of $U\in\R^{n\times2}$ are the vectors
$u_i$; then $A\succeq0$, $A_{ij}=\cos(\theta_i-\theta_j)$, and $A_{ii}=1$.
For this matrix,
\[
   \dotp{Ax}{x}=\Bigl\|\sum_{i=1}^n x_iu_i\Bigr\|_2^2 .
\]
\emph{The exact value of $f^{*}$.} Interchanging the maxima,
\[
   \max_{x\in\{-1,1\}^n}\Bigl\|\sum_i x_iu_i\Bigr\|_2^2
   =\left(\max_{\nrm{v}_2=1}\ \max_{x\in\{-1,1\}^n}
   \sum_i x_i\dotp{u_i}{v}\right)^{\!2}
   =\left(\max_{\nrm{v}_2=1}\ \sum_{i=1}^n
   \bigl|\dotp{u_i}{v}\bigr|\right)^{\!2}.
\]
For $v=(\cos\phi,\sin\phi)$, the inner sum is
$\sum_i|\cos(\theta_i-\phi)|$---a Riemann sum for
a $\pi$-periodic Lipschitz function on a uniform grid with spacing
$\pi/n$, and hence it equals
$\frac{n}{\pi}\int_0^{\pi}|\cos\theta|\,d\theta+\Oh(1)=\frac{2n}{\pi}+\Oh(1)$
uniformly in $\phi$. Therefore,
\[
   f^{*}=\left(\frac{2n}{\pi}+\Oh(1)\right)^{\!2}
   =\frac{4n^2}{\pi^2}+\Oh(n).
\]
\emph{A lower bound for $\mathrm{SDP}$.} The matrix $X=A$ is feasible
in~\eqref{eq:sdprel}, so
\[
   \mathrm{SDP}\ \ge\ \dotp{A}{A}=\sum_{i,j}\cos^2(\theta_i-\theta_j)
   =\frac{n^2}{2}+\frac12\Bigl|\sum_{i=1}^n e^{2\mathrm{i}\theta_i}\Bigr|^2
   =\frac{n^2}{2},
\]
since the $e^{2\mathrm{i}\theta_i}$ are all the $n$th roots of unity,
and their sum is zero for $n\ge2$.
\emph{Conclusion:}
\[
   \frac{\mathrm{SDP}}{f^{*}}\ \ge\ \frac{\pi^2}{8}+o(1)
   \approx1{,}2337,
\]
so the relaxation errs by at least $23{,}4\%$, while
Theorem~\ref{th:pi2} guarantees at most $57{,}1\%$. Both quantities
have the same origin: in both the bound and the example, everything reduces to the mean
value of $|\cos|$, that is, to $2/\pi$.
\end{example}

\subsection{What changed}

Before~\cite{Nes1998}, the quality of semidefinite relaxation
was studied separately for each particular problem, and substantial
results already existed, including the Goemans--Williamson
theory. The contribution of~\cite{Nes1998} lies elsewhere: a
\emph{universal} theorem on the quality of relaxation became available for the entire class
of nonconvex quadratic problems with a positive semidefinite matrix, with
a sharp constant independent of everything else. This is one of
the starting points of modern approximation theory: the combination ``convex
relaxation $+$ randomized rounding $+$ proof of sharpness of the
constant'' became a standard methodology, and subsequent
hardness-of-approximation results turned such constants from
``best known'' into ``optimal under known assumptions.''
The paper~\cite{Nes1998} also considers the more general case of maximizing
a quadratic form over an intersection of ellipsoids, with bounds depending on
the number of constraints; this line continues to this day in problems involving
quadratic mappings, upper bounds for control problems, and
the theory of sums of squares.

%% file: en/r6_traffic.tex

\section{The stable dynamics model of traffic flows
(1998 -- 2003)}
\label{sec:traffic}

\subsection{Problem statement: the network, origin--destination pairs, and equilibrium}

An urban transportation network is described by a directed graph
$\Gamma=(V,E)$: vertices are network nodes, and arcs are road sections. In
modern models of large metropolitan areas, $|V|\sim10^4$\,--\,$10^5$, and
$|E|$ is three to four times larger.

A set of \emph{origin--destination pairs}
$W\subseteq V\times V$ is given; for a pair $w=(i,j)$, the quantity $d_w>0$ is
the demand rate (vehicles per hour), $P_w$ is the set of paths from
$i$ to $j$, and $P=\bigcup_wP_w$. The assignment of flows to paths is a vector
$x=(x_p)_{p\in P}$ in the set
\[
   X=\Bigl\{x\ge0:\ \sum_{p\in P_w}x_p=d_w,\ w\in W\Bigr\},
\]
and the arc flow is $f_e(x)=\sum_{p\ni e}x_p$, or, in matrix notation,
$f=\Theta x$, where $\Theta$ is the arc--path incidence matrix.
Let $\chi_p\in\{0,1\}^{E}$ denote the incidence vector of path $p$.

In 1952, J.~Wardrop~\cite{Wardrop1952} formulated the principle that has underpinned transportation modeling ever since: \emph{at equilibrium, for each origin--destination pair, the travel times on all paths actually used are equal and do not exceed the travel time on any unused path}. This is exactly the Nash condition in a game with a continuum of infinitesimal players, each neglecting the effect of their choice on the total
network load.

\subsection{Background: the Beckmann potential and its weaknesses}

The remarkable observation of M.~Beckmann, C.~McGuire, and
C.~Winsten~\cite{BMW1956} is as follows. If the travel time on arc $e$ is an increasing function $\tau_e(f_e)$ of its own flow alone, then the game is a \emph{potential game}: there exists a function
\begin{equation}\label{eq:beckmann}
   \Psi(f)=\sum_{e\in E}\int_0^{f_e}\tau_e(z)\,dz,
\end{equation}
for which $\partial\Psi/\partial x_p$ equals the travel time on path~$p$, and Wardrop equilibrium coincides with a solution of the convex optimization problem
\begin{equation}\label{eq:beckprob}
   \min\bigl\{\Psi(f):\ f=\Theta x,\ x\in X\bigr\}.
\end{equation}
The function $\Psi$ has no direct physical interpretation in its own right, but its partial derivatives with respect to path variables give precisely the route travel times; this is what turned equilibrium computation into an optimization problem amenable to numerical methods. The empirical formula of the US Bureau
of Public Roads (BPR) is used almost universally for $\tau_e$:
\begin{equation}\label{eq:bpr}
   \tau_e(f_e)=\bt_e\Bigl(1+\gamma\bigl(f_e/\bff_e\bigr)^{1/\mu}\Bigr),
   \qquad \mu=0{,}25,\quad \gamma\in[0{,}15;\,2],
\end{equation}
where $\bt_e$ is the free-flow travel time on the arc and $\bff_e$ is its
capacity.

By the late 1990s, the weaknesses of this construction were becoming increasingly clear~\cite{NesdePalma2003,GasnDornNesShp2014}.

\begin{itemize}[leftmargin=1.6em,itemsep=3pt]
\item \emph{Unobservability of parameters.} The values of $\gamma$ and $\mu$
      in~\eqref{eq:bpr} are fitted rather than measured; equilibrium is sensitive to their choice, particularly when flows are close to capacity.
      Worse still, excessive freedom in calibration creates a risk of overfitting: the model will learn to describe historical data, but its predictive value becomes questionable. The parameters
      in~\eqref{eq:bpr} are empirical and calibrated; the power-law form itself is not derived from microscopic queue
      dynamics, but fitted to data.
\item \emph{Derivatives of unknown functions are required.} Optimal tolls in the Beckmann model have the form
      $\bar\tau_e(f_e)=\tau_e(f_e)+f_e\tau'_e(f_e)$: computing them requires knowledge of the \emph{derivative} of $\tau_e(\cdot)$. Derivative estimation is even more susceptible to the problems listed above than estimation of the function itself.
\item \emph{Mismatch with the mechanism of queue formation.}
      Model~\eqref{eq:bpr} permits
      $f_e>\bff_e$: more vehicles ``pass through'' an arc than it can accommodate, merely traveling more slowly. In reality, exceeding capacity leads not to slower travel on the arc, but to a
      \emph{queue in front of} it. Queues are absent as an object in the model, although the main time losses in overloaded networks are concentrated precisely at bottlenecks and are associated with congestion.
\end{itemize}

\subsection{The Nesterov--de Palma model}

In the work of Yu.\,E.~Nesterov and A.~de Palma~\cite{NesdePalma2003} (the first
version was a 1998 preprint), a radical revision of the
arc description was proposed. Each arc $e$ is assigned \emph{only two}
parameters: $\bt_e$, the minimum arc traversal time, and $\bff_e$, the maximum capacity. Both can be measured directly: the former from the length and speed limit, the latter from the number of lanes and intersection characteristics. The model is thus specified by the graph $\Gamma(V,E,\bt,\bff)$.

The underlying description is as follows. As long as the arc flow is below capacity, all vehicles travel at maximum speed, and the arc traversal time equals $\bt_e$. Once the flow reaches capacity, a queue starts to form before the arc; the arc traversal time $t_e$ is the sum of $\bt_e$ and the queueing time and may be arbitrarily large, while the flow on the arc itself no longer increases. Formally, a state $(f,t)$ is called \emph{stable} (an equilibrium) if it satisfies three conditions:

\begin{equation}\label{eq:sdpost}
   \text{(P1)}\ \ t_e\ge\bt_e,\  f_e\le\bff_e; \ \ \ \
   \text{(P2)}\ \  (t_e-\bt_e)(\bff_e-f_e)=0;\ \ \ \
   \text{(P3)}\ \  f\in\partial\!\!\sum_{w\in W} d_w T_w(t),
\end{equation}
where
\begin{equation}\label{eq:Tw}
   T_w(t)=\min_{p\in P_w}\ \dotp{t}{\chi_p}
\end{equation}
is the shortest-path length for origin--destination pair $w$ with arc lengths $t$.

Each condition has a clear interpretation. (P1) gives the physical constraints: flow cannot exceed capacity, and travel time cannot be less than the time required at maximum speed.
(P2) is the complementary slackness condition: congestion can occur only where flow has reached capacity. (P3) is Wardrop's principle: all demand is assigned to shortest paths (with respect to the \emph{total} times $t$).

\begin{remark}[on the concavity of $T_w$]\label{rem:concave}
The function $T_w(t)$ is the minimum of finitely many linear functions;
it is therefore \emph{concave} and generally nonsmooth. If
$p^{*}$ is a shortest path for the given $t$, then, for every $s$,
\[
   T_w(s)\ \le\ T_w(t)+\dotp{\chi_{p^{*}}}{s-t},
\]
so $\chi_{p^{*}}$ is a \emph{super}gradient, and $\partial$ in~\eqref{eq:sdpost} denotes the superdifferential. It is efficiently computed by Dijkstra's algorithm (the structure of the superdifferential follows from the Dubovitskii--Milyutin theorem; when there are several shortest paths, it is the convex hull of their incidence vectors). An equivalent convex problem is obtained by minimizing the negative of the function, in which case the subgradient $-\chi_{p^{*}}$ appears.
\end{remark}

The main result of~\cite{NesdePalma2003} is that the three
conditions~\eqref{eq:sdpost} defining equilibrium hold at a solution of a single convex optimization problem.

\begin{theorem}[\cite{NesdePalma2003}]\label{th:sd}
Suppose that demand is compatible with the capacities, that is, the set
$\{f=\Theta x:\ x\in X\}\cap\{f\le\bff\}$ is nonempty. Then a pair $(f^{*},t^{*})$
satisfies~\eqref{eq:sdpost} if and only if $t^{*}$ solves the \emph{concave} maximization problem
\begin{equation}\label{eq:sddual}
   \max_{t\ge\bt}\ \Bigl\{\ \sum_{w\in W}d_w\,T_w(t)
   -\dotp{\bff}{t-\bt}\ \Bigr\},
\end{equation}
and $f^{*}=\bff-s^{*}$, where $s^{*}\ge0$ is the vector of Lagrange multipliers for the constraints $t\ge\bt$ in~\eqref{eq:sddual}. Problem~\eqref{eq:sddual} is dual to the linear programming problem
\begin{equation}\label{eq:sdprimal}
  \min_{f}\ \bigl\{\dotp{\bt}{f}:\ f=\Theta x,\ x\in X,\
   f\le\bff\bigr\}
\end{equation}
and has a finite solution precisely under the stated compatibility
assumption.
\end{theorem}

\begin{proof}[Proof sketch]
Condition (P3) is the first-order optimality condition for
$\max_{t}\{\sum_wd_wT_w(t)-\dotp{f}{t-\bt}\}$. Condition (P2) is complementary slackness for
the constraint $t\ge\bt$ with Lagrange multipliers $s=\bff-f\ge0$. There is no need to impose $f\le\bff$ separately: Lagrange multipliers for inequality constraints are automatically nonnegative. This gives~\eqref{eq:sddual}. It is an LP problem written compactly:
\[
   \max_{t\ge\bt,\ \{T_w\}}\Bigl\{\sum_w d_wT_w-\dotp{\bff}{t-\bt}:\
   T_w\le\dotp{t}{\chi_p},\ p\in P_w,\ w\in W\Bigr\},
\]
and its dual is~\eqref{eq:sdprimal}. Crucially,~\eqref{eq:sddual} can be solved by a \emph{primal-dual} method: the variable $t$ itself gives the travel times, but flows $f$ are available only through the dual multipliers.
\end{proof}

Thus equilibrium computation reduces to an LP problem. The objective
in~\eqref{eq:sdprimal} has an interpretation: $\dotp{\bt}{f}$ is the total time drivers would spend \emph{in motion} (excluding congestion). The model's authors~\cite{NesdePalma2003}
interpret the objective in~\eqref{eq:sddual} as maximizing the number of
freely moving vehicles, that is, vehicles not standing in
queues. Most importantly, the multipliers $\eta_e:=t_e-\bt_e\ge0$ have
a direct physical meaning as the \emph{time spent waiting in congestion on arc~$e$}: the model yields not only flows, but also delays.

\begin{remark}[on nonuniqueness]\label{rem:nonuniq}
In the Beckmann model with strictly increasing $\tau_e$, equilibrium arc flows are unique. In the stable dynamics model, this is generally not so: uniqueness is typical (in general position), but in degenerate cases the solution of the LP problem~\eqref{eq:sdprimal} and/or its dual may be nonunique. Ways of selecting the equilibrium that ``materializes'' are already discussed in~\cite{NesdePalma2003}; a modern approach is to add a strictly convex (for example, entropic) term to the linear objective, simultaneously restoring uniqueness and allowing the bounded rationality of
drivers to be described (stochastic equilibrium).
\end{remark}

\subsection{The stable dynamics model as a limit of the Beckmann model}

The passage from~\eqref{eq:beckmann} to~\eqref{eq:sdprimal} can be made rigorous, and this is useful for two reasons: first, it shows that the new model is not ``different,'' but a limiting case of its predecessor; second, a regularized (stochastic) variant can be extracted from the proof.
We give the derivation following~\cite{GasnDornNesShp2014}.

For $\tau_e$, take a family of \emph{interior barriers} for the domain $f_e<\bff_e$:
\begin{equation}\label{eq:barlat}
   \tau^{\mu}_e(f_e)=\bt_e\Bigl(1-\mu\ln\bigl(1-f_e/\bff_e\bigr)\Bigr),
   \qquad \mu>0 .
\end{equation}
These are increasing, convex, smooth functions tending to $+\infty$ as $f_e\to\bff_e$. Write $\sigma_e(f_e)=\int_0^{f_e}\tau^\mu_e(z)\,dz$ and compute the conjugate function. The condition
$\sigma'_e(f_e)=t_e$ gives
$f_e=\bff_e\bigl(1-e^{-(t_e-\bt_e)/(\bt_e\mu)}\bigr)$, and direct
calculation yields
\begin{equation}\label{eq:sigmastar}
   \sigma^{*}_e(t_e)=\bff_e\,(t_e-\bt_e)
   -\mu\,\bff_e\bt_e\Bigl(1-e^{-(t_e-\bt_e)/(\bt_e\mu)}\Bigr),
   \qquad t_e>\bt_e .
\end{equation}
Passing to the dual problem in~\eqref{eq:beckprob} (the interchange
of $\min$ and $\max$ is justified by the minimax theorem in Sion's form), we obtain
\begin{align}
   \min_{f=\Theta x,\,x\in X}\sum_e\sigma_e(f_e)
   &=\max_{t}\Bigl\{\min_{f=\Theta x,\,x\in X}\dotp{t}{f}
     -\sum_e\sigma^{*}_e(t_e)\Bigr\}\notag\\
   &=\max_{t}\Bigl\{\sum_{w}d_wT_w(t)-\dotp{\bff}{t-\bt}
     +\mu\sum_e\bff_e\bt_e
     \Bigl(1-e^{-\frac{t_e-\bt_e}{\bt_e\mu}}\Bigr)\Bigr\}\notag
\\
   &\xrightarrow[\ \mu\to0+\ ]{}\
     \max_{t\ge\bt}\Bigl\{\sum_{w}d_wT_w(t)-\dotp{\bff}{t-\bt}\Bigr\},
   \notag
\end{align}
using
$\min\limits_{f=\Theta x,\,x\in X}\dotp{t}{f}=\sum\limits_wd_wT_w(t)$: for
fixed arc lengths $t$, the cheapest choice is to send all demand along
shortest paths. We have obtained exactly~\eqref{eq:sddual}.

\begin{figure}[htbp]
\centering
\includegraphics[width=0.7\textwidth]{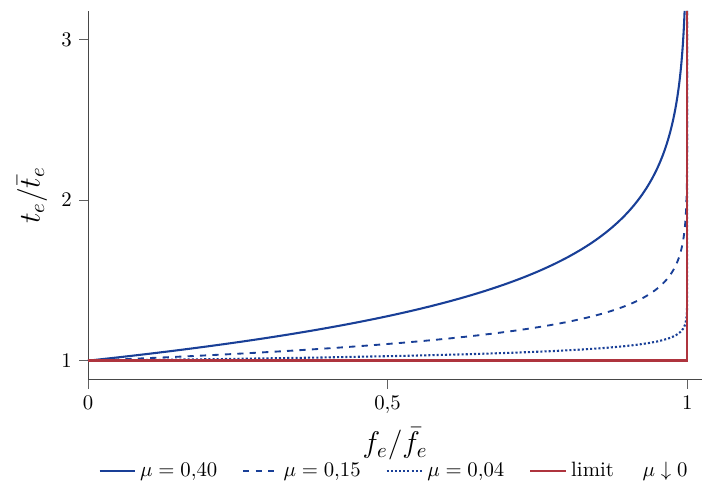}
\caption{The arc delay function and its limiting graph. The blue curves
are given by
$t_e/\bt_e=1-\mu\ln(1-f_e/\bff_e)$ for $0\le f_e<\bff_e$.
As $\mu\downarrow0$, they converge graphically to the stable dynamics model shown
in red: $t_e=\bt_e$ for
$0\le f_e\le\bff_e$, while at $f_e=\bff_e$ any
$t_e\ge\bt_e$ is allowed. The vertical part corresponds to the waiting time
$\eta_e=t_e-\bt_e$ on a saturated arc; values $f_e>\bff_e$
are infeasible.}
\label{fig:latency}
\end{figure}

The particular choice of barrier~\eqref{eq:barlat} is inessential. If
$\mu>0$ is treated as a parameter in the BPR
function~\eqref{eq:bpr}, an analogous calculation gives
\[
   \sigma^{*}_e(t_e)=\bff_e
   \Bigl(\frac{t_e-\bt_e}{\bt_e\gamma}\Bigr)^{\!\mu}
   \frac{t_e-\bt_e}{1+\mu}
   \xrightarrow[\mu\downarrow0]{}
   \bff_e\,(t_e-\bt_e).
\]
Thus, as $\mu\downarrow0$, that is, as the exponent
$1/\mu\to\infty$, the same limiting term in the dual
problem emerges. Figure~\ref{fig:latency} illustrates the corresponding transition for
the barrier family~\eqref{eq:barlat}.

\subsection{Braess's paradox}

Braess's paradox (1968)~\cite{Braess1968} is one of the best-known phenomena in transportation theory: adding a new road to a network can make road users worse off because Wardrop equilibrium is not a social optimum. We give an example from~\cite{NesdePalma2003} itself (following its presentation in~\cite[Section~6]{GasnDornNesShp2014}); it shows that the paradox takes a particularly clear form in the stable dynamics model.

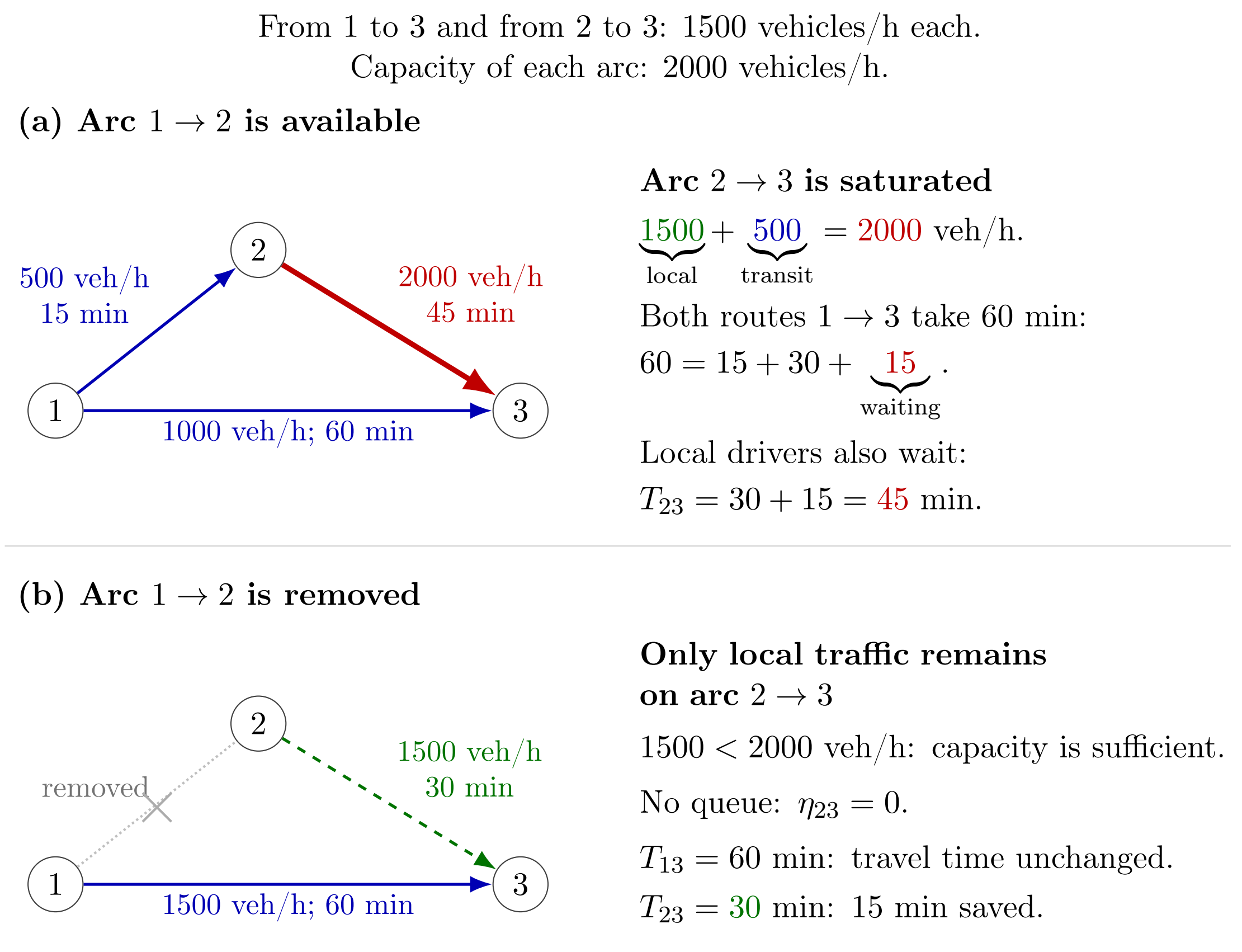
\begin{figure}[h]
\centering
{\small From $1$ to $3$ and from $2$ to $3$: $1500$ vehicles/h each.\par
Capacity of each arc: $2000$ vehicles/h.\par}
\medskip
\begin{tikzpicture}[>=Latex,
   vx/.style={circle,draw=black!75,fill=white,minimum size=6.5mm,
              inner sep=1pt,font=\small},
   transit/.style={->,line width=1pt,blue!70!black},
   localflow/.style={->,line width=1pt,green!45!black,dashed},
   congested/.style={->,line width=1.8pt,red!75!black},
   roadlabel/.style={font=\footnotesize,align=center,inner sep=2pt},
   explanation/.style={anchor=north west,font=\small,align=left,
                       text width=7.0cm,inner sep=0pt}]
  \path[use as bounding box] (0,-6.85) rectangle (14.55,3.95);
  \node[anchor=west,font=\small\bfseries] at (0,3.75)
    {(a) Arc $1\to2$ is available};
  \node[vx] (a1) at (0.6,0.35) {$1$};
  \node[vx] (a2) at (3.0,2.25) {$2$};
  \node[vx] (a3) at (6.1,0.35) {$3$};
  \draw[transit] (a1)--(a2)
    node[roadlabel,midway,above left,blue!70!black] {$500$ veh/h\\$15$ min};
  \draw[transit] (a1)--(a3)
    node[roadlabel,midway,below,blue!70!black] {$1000$ veh/h; $60$ min};
  \draw[congested] (a2)--(a3)
    node[roadlabel,midway,above right,red!75!black] {$2000$ veh/h\\$45$ min};
  \node[explanation] at (7.5,3.22) {
    \textbf{Arc $2\to3$ is saturated}\\[3pt]
    $\underbrace{\textcolor{green!45!black}{1500}}_{\text{local}}
     +\underbrace{\textcolor{blue!70!black}{500}}_{\text{transit}}
      =\textcolor{red!75!black}{2000}$ veh/h.\\[5pt]
    Both routes $1\to3$ take $60$ min:\\[2pt]
    $60=15+30+\underbrace{\textcolor{red!75!black}{15}}_{\text{waiting}}$.\\[5pt]
    Local drivers also wait:\\[2pt]
    $T_{23}=30+15=\textcolor{red!75!black}{45}$ min.
  };
  \draw[black!15,line width=0.5pt] (0,-1.25)--(14.5,-1.25);
  \begin{scope}[yshift=-5.6cm]
    \node[anchor=west,font=\small\bfseries] at (0,3.75)
      {(b) Arc $1\to2$ is removed};
    \node[vx] (b1) at (0.6,0.35) {$1$};
    \node[vx] (b2) at (3.0,2.25) {$2$};
    \node[vx] (b3) at (6.1,0.35) {$3$};
    \draw[gray!50,densely dotted,line width=0.8pt] (b1)--(b2)
      node[roadlabel,midway,above left,text=black!55] {removed};
    \draw[gray!65,line width=0.8pt] (1.63,1.09)--(1.97,1.43);
    \draw[gray!65,line width=0.8pt] (1.63,1.43)--(1.97,1.09);
    \draw[transit] (b1)--(b3)
      node[roadlabel,midway,below,blue!70!black] {$1500$ veh/h; $60$ min};
    \draw[localflow] (b2)--(b3)
      node[roadlabel,midway,above right,green!45!black] {$1500$ veh/h\\$30$ min};
    \node[explanation] at (7.5,3.22) {
      \textbf{Only local traffic remains\\on arc $2\to3$}\\[4pt]
      $1500<2000$ veh/h: capacity is sufficient.\\[5pt]
      No queue: $\eta_{23}=0$.\\[5pt]
      $T_{13}=60$ min: travel time unchanged.\\[3pt]
      $T_{23}=\textcolor{green!45!black}{30}$ min: $15$ min saved.
    };
  \end{scope}
\end{tikzpicture}
\caption{The mechanism of Braess's paradox in the stable dynamics model
(Example~\ref{ex:braess}, network from~\cite{NesdePalma2003}).
Each origin generates $1500$ vehicles/h. The blue flow of pair $(1,3)$
splits in the original network: $1000$ vehicles/h travel directly, and $500$ vehicles/h
join the $1500$ local drivers on arc $(2,3)$.
Equality of the two route travel times gives $15$ minutes of waiting
on the saturated arc for both groups.
After arc $(1,2)$ is removed, transit traffic takes the direct route;
pair $(2,3)$ saves $15$ minutes, while the travel time for pair $(1,3)$ is unchanged.
The labels on the arcs give flow and total travel time, including waiting.}
\label{fig:braess}
\end{figure}

\begin{example}[\cite{NesdePalma2003}]\label{ex:braess}
Consider the network in Fig.~\ref{fig:braess}: three vertices and three arcs with
parameters
\[
   \bt_{13}=60\ \text{min},\quad
   \bt_{12}=15\ \text{min},\quad
   \bt_{23}=30\ \text{min},
\]
all capacities being $2000$~vehicles/hour. The origin--destination
pairs are $(1,3)$ and $(2,3)$, with demand $d_{13}=d_{23}=1500$~vehicles/hour. For pair $(2,3)$ there is only one path; for pair $(1,3)$ there are two: the direct path $1\!\to\!3$ (free-flow time $60$~min) and the indirect path $1\!\to\!2\!\to\!3$ (free-flow time $15+30=45$~min).

\emph{Primal problem.} Let $y$ be the flow for pair $(1,3)$ passing through vertex~$2$. Then $f_{12}=y$, $f_{13}=1500-y$, $f_{23}=1500+y$, and~\eqref{eq:sdprimal} becomes
\[
   \min_{0\le y\le1500}\ \bigl\{45y+60(1500-y)+30\cdot1500\bigr\}
   \quad\text{subject to}\quad 1500+y\le2000 .
\]
The objective is $135\,000-15y$, so $y$ must be maximized; the constraint on arc $(2,3)$ gives $y\le500$. Thus,
\[
   y^{*}=500:\qquad f_{12}=500,\quad f_{13}=1000,\quad
   f_{23}=2000=\bff_{23}\ \text{(saturated)} .
\]

\emph{Dual problem.} Arcs $(1,2)$ and $(1,3)$ are not saturated, so
$\eta_{12}=\eta_{13}=0$. Both paths for pair $(1,3)$ are used, so their
travel times are equal:
\[
   \underbrace{60}_{1\to3}
   =\underbrace{15+30+\eta_{23}}_{1\to2\to3}
   \quad\Longrightarrow\quad
   \eta^{*}_{23}=15\ \text{min}.
\]
Consequently, the equilibrium travel times are
\[
   T_{13}=60\ \text{min}=1\ \text{hour},\qquad
   T_{23}=30+15=45\ \text{min}.
\]

\emph{The paradox.} Remove arc $(1,2)$. All demand for pair $(1,3)$ then
travels directly: $f_{13}=1500<2000$ and $f_{23}=1500<2000$, so no
arc is saturated, $\eta^{*}=0$, and
\[
   T_{13}=60\ \text{min}\ \text{(unchanged)},\qquad
   T_{23}=30\ \text{min}\ \text{(previously 45)} .
\]
In other words, \emph{prohibiting travel on arc $(1,2)$ improves the situation
of one group of road users and harms no one}: drivers
traveling from~2 to~3 save 15~minutes, while nothing changes for those traveling from~1 to~3.
The presence of arc $(1,2)$ causes
``transit'' drivers of pair $(1,3)$ to saturate arc $(2,3)$ and create
a queue in which the ``local'' drivers of pair $(2,3)$ also wait.

\emph{Comparison with the Beckmann model.} The
classical Braess paradox~\cite{Braess1968} itself is formulated precisely for Wardrop equilibrium with separable increasing delay functions, that is, within the standard Beckmann model family, and removing an arc there can also improve the situation of all users. It would be incorrect to claim that ``this does not happen in the Beckmann model.''
\end{example}

\subsection{Tolls}

The dual variables $\eta^{*}=t^{*}-\bt$ also have another interpretation: optimal tolls. The corresponding statement from~\cite{GasnDornNesShp2014} is as follows.

\begin{proposition}\label{pr:tolls}
Let $\eta^{*}=t^{*}-\bt\ge0$ be the equilibrium delays, and set the toll on arc $e$, expressed in time units, equal to $\eta^{*}_e$. Then the equilibrium flow assignment $f^{*}$ is unchanged.
\end{proposition}

\begin{proof}
By the method of Lagrange multipliers, problem~\eqref{eq:sdprimal} has the same solution $f^{*}$ as the problem without the constraint $f\le\bff$:
\[
   \min_{f=\Theta x,\,x\in X}\ \sum_e\bigl(\bt_ef_e
   +\eta^{*}_e(f_e-\bff_e)\bigr)
   \quad\Longleftrightarrow\quad
   \min_{f=\Theta x,\,x\in X}\ \dotp{\bt+\mu^{*}}{f},
\]
since the term $-\dotp{\mu^{*}}{\bff}$ is constant. But
$\bt+\mu^{*}=t^{*}$, so the same vector $f^{*}$ also minimizes $\min\{\dotp{t^{*}}{f}:\ f=\Theta x,\ x\in X\}$ and, since $f^{*}\le\bff$, also $\min\{\dotp{t^{*}}{f}:\ f=\Theta x,\ x\in X,\ f\le\bff\}$. Moreover, the same holds for the entire ``segment'' of weights $\tilde t=\alpha\bt+(1-\alpha)t^{*}$, $\alpha\in[0,1]$.
\end{proof}

\practice The substantive conclusion is this: without tolls, the case
$\alpha=0$ is realized, so users' actual costs equal $t^{*}$ and
include time spent in queues. If tolls
$\eta^{*}_e$ are imposed, the route assignment remains unchanged, while the actual
\emph{time} costs become $\bt$. The procedure can
be implemented adaptively: gradually increase the toll on arcs where
queues are observed until the queues disappear. We emphasize, however, that
this is an economic interpretation requiring additional
modeling assumptions: about converting monetary costs into time and about
the possibility of maintaining flow exactly at capacity
without accumulating a queue.

\subsection{A multistage model}

The stable dynamics model covers only one component of transportation modeling: assigning given origin--destination demands to paths. In practice, \emph{multistage} models are used to determine simultaneously (1)~the origin--destination matrix $\{d_{ij}\}$,
(2)~the split of demand between modes of travel (private and public transport), and (3)~the assignment of flows to paths. Remarkably, all three stages can be combined into \emph{one} convex optimization problem~\cite{GasnDornNesShp2014}.

The first stage is described by an \emph{entropy model}: the most likely origin--destination matrix solves the
entropy-linear programming problem
\begin{equation}\label{eq:entropy}
   \max\Bigl\{-\sum_{i,j}d_{ij}\ln\frac{d_{ij}}{e}
   -\beta\sum_{i,j}d_{ij}T_{ij}:\
   d\ge0,\ \sum_jd_{ij}=l_i,\ \sum_id_{ij}=w_j\Bigr\},
\end{equation}
where $l_i$ and $w_j$ are prescribed departure and arrival rates for districts, $T_{ij}$ is travel time, and $\beta>0$ is a parameter (the reciprocal of a characteristic travel time). The justification of~\eqref{eq:entropy} uses Markovian ``exchange'' dynamics and the ergodic theorem: the stationary measure concentrates exponentially near the solution of~\eqref{eq:entropy}, and the functional in~\eqref{eq:entropy} is both a Lyapunov function for the dynamics and a Sanov functional. A substantial difficulty is that $T_{ij}$ themselves depend on $d$; separation of time scales (years for origin--destination demand formation, weeks for route equilibrium) allows, by Tikhonov's theorem, the dependence $T_{ij}(d)$ to be substituted into~\eqref{eq:entropy}.

The key device of~\cite{GasnDornNesShp2014} is to embed problem~\eqref{eq:sddual} itself in~\eqref{eq:entropy}, rather than substituting an explicit
dependence $T_{ij}(d)$, which cannot be written down. As a result (after eliminating the inner maximum over $d$ and introducing dual variables $\lambda^L,\lambda^W$ for the balance constraints), one obtains a single nonsmooth convex optimization problem of the form
\begin{multline}\label{eq:threestage}
   \min_{\substack{t^{\text{pr}}\ge\bt^{\text{pr}},\
   t^{\text{pu}}\ge\bt^{\text{pu}}\\ \lambda^L,\lambda^W}}
   \Bigl\{\
   \ln\!\Bigl(\sum_{i,j}
   e^{-\beta\min\{T^{\text{pr}}_{ij}(t^{\text{pr}}),
   T^{\text{pu}}_{ij}(t^{\text{pu}})\}-\lambda^L_i-\lambda^W_j}\Bigr)
   +\dotp{\lambda^L}{l}+\dotp{\lambda^W}{w}\\
   {}+\beta\dotp{\bff^{\text{pr}}}{t^{\text{pr}}-\bt^{\text{pr}}}
   +\beta\dotp{\bff^{\text{pu}}}{t^{\text{pu}}-\bt^{\text{pu}}}
   \ \Bigr\},
\end{multline}
where the indices ``pr'' and ``pu'' refer to private and public transport, and the $\min$ operator in the exponent represents mode choice. The origin--destination matrix and flows are recovered from the solution: $d_{ij}$ is proportional to the corresponding exponential, and $f^{\text{pr}}=\bff^{\text{pr}}-s^{\text{pr}}$, where $s^{\text{pr}}$ are the dual multipliers of the constraints $t^{\text{pr}}\ge\bt^{\text{pr}}$.

Problem~\eqref{eq:threestage} is nonsmooth and convex, with a bounded Lipschitz constant but a poorly controlled initial residual, and, most importantly, requires simultaneous recovery of \emph{primal and dual} solutions. This makes it suitable for the methods of \S\,\ref{sec:pd} (primal-dual schemes) and \S\,\ref{sec:univ} (the universal method, which automatically adapts to the actual smoothness of the piecewise linear function $T_w$). Problem~\eqref{eq:threestage} also admits a stochastic variant in which $\min$ is replaced by a ``soft minimum'' (log-sum-exp) with parameter $T>0$: this is exactly the smoothing of \S\,\ref{sec:smooth}, and it makes the problem smooth.

\subsection{What changed and the current state}

\textbf{(i) Observability and computability.} The model retains exactly two
parameters per arc, both measurable. The nonlinear
problem~\eqref{eq:beckprob} is replaced by an LP problem, and its dual
form~\eqref{eq:sddual} is a nonsmooth \emph{concave}
maximization problem whose supergradient is computed by Dijkstra's algorithm. This
makes computations on networks of realistic size possible.

\smallskip
\textbf{(ii) Stability---and exactly what it means.} The model's
name is related to stability, but this must be formulated locally.
The optimal value of~\eqref{eq:sdprimal}, as a function of the data
$(\bt,\bff,d)$, is piecewise linear and locally Lipschitz on regions where
feasibility and boundedness of dual solutions are preserved; at
the feasibility boundary, a small change in demand can make the problem
infeasible, in which case equilibrium simply does not exist
(Theorem~\ref{th:sd}). The comparison with the Beckmann model is as follows: there, for
$\beta=4$, a small change in demand near capacity leads to
a sharp change in predicted travel times \emph{within} the region
where equilibrium exists; here, all the sharpness is shifted to the feasibility
boundary itself, which is easy to monitor in applications.

\smallskip
\textbf{(iii) Delays as an output of the model.} The model predicts
\emph{where} delays occur and \emph{how long the waiting time is}, and
simultaneously gives optimal tolls, without computing
derivatives of unknown functions as required in the Beckmann model.
We emphasize that the dual variable $\eta_e=t_e-\bt_e$ has
the dimension of time and specifies precisely a delay; recovering
the spatial queue length or the number of queued vehicles requires
an additional dynamic service model.

\smallskip
\textbf{(iv) Current state.} The stable dynamics model has become
an important tool in transportation modeling, although it has not displaced
the classical static assignment models with smooth delay functions
that still dominate commercial packages. Systematic
accounts of the theory and computational methods are given in the
monographs~\cite{GasnTransp2013,GasnTranspMono2020}; the three-stage
model and its stochastic variant are presented in~\cite{GasnDornNesShp2014};
primal-dual methods for the dual problem and their numerical
comparison are given in~\cite{GasDvurDornMax, GasnGasnNes2018,Baimurzina2019,
KubentayevaGasnikov2021}. A separate, actively developing
direction is that of \emph{combined} models in which capacity
constraints from~\eqref{eq:sdprimal} are combined with
entropy or Beckmann terms and a nested choice
structure~\cite{Kubentayeva2024}. Finally, the phenomenon in
Example~\ref{ex:braess} has itself become an object of study~\cite{DornShit}:
dual variables provide a criterion for \emph{detecting}
Braess's paradox in large networks without enumerating subnetworks.

%% file: en/r7_smooth_pd.tex

\section{Smoothing technique (2003 -- 2005)}
\label{sec:smooth}

Throughout, by \emph{smoothing technique} we mean the specific construction of
Yu.\,E.~Nesterov that uses the algebraic model~\eqref{eq:algmodel};
smoothing in a broader sense (for example, by Gaussian convolution)
will be discussed separately in \S\,\ref{sec:zo}.

\subsection{Background: the impasse in nonsmooth optimization}

By the late 1970s, the theory of nonsmooth convex optimization appeared
complete. For the class of convex functions with Lipschitz constant $M$ on
a set of diameter $R$, N.\,Z.~Shor's subgradient method achieves accuracy
$\eps$ in $\Oh(M^2R^2/\eps^2)$ iterations, and, as shown by
A.\,S.~Nemirovski and D.\,B.~Yudin, \emph{this bound cannot be improved} in
the oracle model (Table~\ref{tab:map}). Despite the enormous gap compared with the smooth case ($\eps^{-2}$ versus $\eps^{-1/2}$), there seemed to be nothing left to improve.

The observation made in~\cite{Nes2005smooth} (a CORE preprint of
2003) is as follows. The lower bound concerns a \emph{black
box}. But nonsmooth functions encountered in practice almost never
come as black boxes: they are usually specified by an \emph{algebraic
model} in the sense of \S\,\ref{sec:models}, most often in the form
\begin{equation}\label{eq:maxstruct}
   f(x)=\max_{u\in U}\ \bigl\{\dotp{Ax}{u}-\hat\phi(u)\bigr\},
   \qquad x\in Q,
\end{equation}
where $U$ is a compact convex set, $\hat\phi$ is a convex function, and $A$ is
a linear operator. The maximum of the coordinates, the norm
$\nrm{\cdot}_1$, total variation, the objective of a matrix game,
and the function $T_w(t)$ of \S\,\ref{sec:traffic} all have such representations. If the model is known, one can look inside the black box.

\subsection{Smoothing technique}

Fix the geometry on both sets---this is essential, and we
state the requirements in full.

\begin{itemize}[leftmargin=1.6em,itemsep=2pt]
\item On $U$, choose a \emph{prox-function} $d_U$, strongly convex with
      parameter $\sigma_U>0$ with respect to the norm $\nrm{\cdot}_U$, with
      $\min_{u\in U}d_U(u)=0$; set $D_U=\max_{u\in U}d_U(u)$.
\item On $Q$, choose a prox-function $d_Q$, strongly convex with parameter
      $\sigma_Q>0$ with respect to the norm $\nrm{\cdot}$, with
      $\min_{x\in Q}d_Q(x)=0$; set $R^2=\max_{x\in Q}d_Q(x)$. For
      definiteness, throughout the following we assume $\sigma_Q=1$ (this can always
      be achieved by the replacement $d_Q\mapsto d_Q/\sigma_Q$; when $\sigma_Q\ne1$,
      $R^2$ in the bounds~\eqref{eq:smoothrate} must be replaced by
      $R^2/\sigma_Q$).
\end{itemize}
Normalization is necessary because $R^2$ alone is not invariant:
the replacement $d_Q\mapsto cd_Q$ changes $R^2$, but should not change the problem's complexity.

Define the smoothed function
\begin{equation}\label{eq:fmu}
   f_\mu(x)=\max_{u\in U}\ \bigl\{\dotp{Ax}{u}-\hat\phi(u)
   -\mu\,d_U(u)\bigr\},\qquad \mu>0 .
\end{equation}

\begin{theorem}[\cite{Nes2005smooth}]\label{th:smooth}
Let a function $f$ of the form~\eqref{eq:maxstruct} be given on a compact convex set $Q$, with known $A$, $U$, $\hat\phi$, and $\sigma_Q=1$. Suppose an \emph{oracle} is available for computing
$\grad f_\mu$, that is, solving the auxiliary problem~\eqref{eq:fmu}. Then $f_\mu$ is convex and continuously differentiable,
\begin{equation}\label{eq:sandwich}
   f_\mu(x)\ \le\ f(x)\ \le\ f_\mu(x)+\mu D_U
   \qquad\forall x\in Q,
\end{equation}
its gradient $\grad f_\mu$ is Lipschitz continuous with constant
\begin{equation}\label{eq:Lmu}
   L_\mu=\frac{\nrm{A}^2}{\mu\,\sigma_U},\qquad
   \nrm{A}=\max\bigl\{\dotp{Ax}{u}:\ \nrm{x}\le1,\ \nrm{u}_U\le1\bigr\},
\end{equation}
and applying the \emph{proximal} version of the fast gradient
method to $f_\mu$ (that is, scheme~\eqref{eq:fgm} with the gradient step replaced by
the composite mapping~\eqref{eq:composite} with $\Psi=\delta_Q$), with
\[
   \mu=\frac{2\nrm{A}R}{N}\,\frac{1}{\sqrt{\sigma_U D_U}}
\]
gives, after $N$ iterations, a point $x_N\in Q$ satisfying
\begin{equation}\label{eq:smoothrate}
   f(x_N)-\min_{x\in Q}f(x)\ \le\
   \frac{4\nrm{A}R}{N}\sqrt{\frac{D_U}{\sigma_U}},
   \qquad\text{that is,}\qquad
   N(\eps)=\Oh\!\left(\frac{\nrm{A}R}{\eps}
   \sqrt{\frac{D_U}{\sigma_U}}\right).
\end{equation}
\end{theorem}

The logic behind~\eqref{eq:smoothrate} is that the total error is the sum of the smoothing error $\mu D_U$ and the method's error $4L_\mu R^2/N^2=4\nrm{A}^2R^2/(\mu\sigma_U N^2)$; minimizing their sum over $\mu$ gives $\mu^{*}=2\nrm{A}R/(N\sqrt{\sigma_U D_U})$ and the value $4\nrm{A}R\sqrt{D_U/\sigma_U}/N$. Crucially, \emph{the smoothing parameter is chosen in advance based on the planned number of iterations}: the larger the budget, the smaller $\mu$ and the more accurate the approximation. Replacing the gradient step by a proximal step cannot be omitted: when $Q\ne\R^n$, iterations~\eqref{eq:fgm} may leave $Q$.

The result is \emph{$\Oh(\eps^{-1})$ instead of $\Oh(\eps^{-2})$}---a quadratic
speedup for a very broad class of nonsmooth problems. The
Nemirovski--Yudin lower bound is not violated: the method is not a
black-box method. This is perhaps the purest example of the program of
\S\,\ref{sec:lookinside}: the bound is not refuted but bypassed through
knowledge of the model.

\begin{example}[smoothing the maximum: log-sum-exp]\label{ex:lse}
Let $f(x)=\max_{1\le i\le m}\dotp{a_i}{x}$, so in
\eqref{eq:maxstruct}, $U=\Delta_m=\{u\in\R^m_+:\sum_iu_i=1\}$ is
the standard simplex,
$\hat\phi\equiv0$, and the rows of $A$ are the vectors $a_i$. Take the entropy
prox-function $d_U(u)=\ln m+\sum_{i=1}^m u_i\ln u_i$; it is nonnegative
on $\Delta_m$, vanishes at the center of the simplex, and is strongly convex with
$\sigma_U=1$ with respect to $\nrm{\cdot}_1$ (Pinsker's inequality), with
$D_U=\ln m$. Explicit evaluation of the maximum in~\eqref{eq:fmu} gives
\[
   f_\mu(x)=\mu\ln\Bigl(\frac1m\sum_{i=1}^m
   e^{\dotp{a_i}{x}/\mu}\Bigr),
\]
and~\eqref{eq:sandwich} becomes the familiar two-sided inequality
$f(x)-\mu\ln m\le f_\mu(x)\le f(x)$. The Lipschitz constant~\eqref{eq:Lmu}
equals $\max_i\nrm{a_i}_2^2/\mu$. Thus log-sum-exp
is not a heuristic ``soft'' replacement for the maximum, but an exact
instance of the general construction, with controlled error and an explicitly
optimized parameter. This is precisely the construction that turns $\min$
in~\eqref{eq:threestage} into a smooth function.
\end{example}

\begin{example}[matrix games and optimal transport]\label{ex:game}
For $\min_{x\in\Delta_n}\max_{u\in\Delta_m}\dotp{Ax}{u}$, the smoothing
technique with entropy on both simplices gives a method requiring
$\Oh\bigl(\max_{i,j}|A_{ij}|\sqrt{\ln n\,\ln m}\,/\eps\bigr)$
iterations; before~\cite{Nes2005smooth}, methods with complexity
$\Oh(\eps^{-2})$ were standard. Each iteration requires one multiplication by $A$
\emph{and} one by $A^{\!\top}$ (the gradient of the smoothed function with respect to $x$
requires computing $Ax$, while recovery of the inner maximum requires
multiplication by $A^{\!\top}$), that is, $\Oh(\mathrm{nnz}(A))$ operations, where
$\mathrm{nnz}(\cdot)$ is the number of nonzero entries of a matrix.

The same construction underlies modern optimal transport
algorithms: entropic regularization of the
Monge--Kantorovich problem~\cite{Cuturi2013} is precisely
smoothing~\eqref{eq:fmu}. Applying accelerated methods to the smoothed
problem gives theoretical complexity bounds better than those for the
Sinkhorn algorithm~\cite{DvurGasnKrosh2018}. We emphasize that this comparison is
\emph{theoretical}: the practical balance depends on the required
accuracy, data sparsity, numerical stability, and the cost of
matrix--vector multiplication, and neither approach unconditionally
dominates the other.
\end{example}

\subsection{The excessive gap condition}
\label{sec:excessgap}

The scheme of Theorem~\ref{th:smooth} has a drawback: the smoothing parameter
$\mu$ must be fixed in advance, knowing the budget $N$. The remedy found
in~\cite{Nes2005excess} is to use a symmetric scheme that
alternates steps in the primal and dual spaces.

Alongside~\eqref{eq:maxstruct}, consider the dual problem
$\max_{u\in U}\phi(u)$, where
$\phi(u)=-\hat\phi(u)+\min_{x\in Q}\dotp{Ax}{u}$, and its smoothed
version $\phi_{\mu_1}$, defined analogously to~\eqref{eq:fmu} using
a prox-function on $Q$. Weak duality gives
$\phi(u)\le f(x)$ for every $x\in Q$, $u\in U$; but for the \emph{smoothed}
functions, $f_{\mu_2}(x)\le f(x)$ and $\phi(u)\le\phi_{\mu_1}(u)$,
so the inequality can be reversed. The \emph{excessive gap}
condition is
\begin{equation}\label{eq:excessgap}
   f_{\mu_2}(\bar x)\ \le\ \phi_{\mu_1}(\bar u)
\end{equation}
for some $\bar x\in Q$, $\bar u\in U$. Equation~\eqref{eq:excessgap}
immediately gives a quality bound for the primal-dual pair:
\begin{equation}\label{eq:excessgapbound}
   0\ \le\ \max\bigl\{f(\bar x)-f^{*},\,f^{*}-\phi(\bar u)\bigr\}
   \ \le\ f(\bar x)-\phi(\bar u)\ \le\ \mu_1D_Q+\mu_2D_U ,
\end{equation}
where $D_Q=\max_{x\in Q}d_Q(x)$ and $D_U=\max_{u\in U}d_U(u)$.

Remarkably, condition~\eqref{eq:excessgap} can be \emph{maintained}
while decreasing the smoothing parameters. Let $T_{\mu_2}$ denote the primal
gradient mapping~\eqref{eq:gradmap} for $f_{\mu_2}$, and let
$x_{\mu_1}(u)$ and $u_{\mu_2}(x)$ be the solutions of the corresponding inner
problems.

\begin{theorem}[\cite{Nes2005excess}]\label{th:excess}
Suppose $\bar x\in Q$ and $\bar u\in U$ satisfy~\eqref{eq:excessgap}
for some $\mu_1,\mu_2>0$. Fix $\tau\in(0,1)$, put
$\mu_1^{+}=(1-\tau)\mu_1$, and set
\begin{equation}\label{eq:excessstep}
   \hat x=(1-\tau)\bar x+\tau x_{\mu_1}(\bar u),\qquad
   \bar u_{+}=(1-\tau)\bar u+\tau u_{\mu_2}(\hat x),\qquad
   \bar x_{+}=T_{\mu_2}(\hat x).
\end{equation}
Then the pair $(\bar x_{+},\bar u_{+})$ satisfies~\eqref{eq:excessgap} with
parameters $\mu_1^{+}$ and $\mu_2$ for every sufficiently small $\tau$,
namely when
$\dfrac{\tau^2}{1-\tau}\le\dfrac{\mu_1}{L(f_{\mu_2})}$.
\end{theorem}

Applying rule~\eqref{eq:excessstep} alternately in the primal and dual
spaces yields a method with complexity $\Oh(\eps^{-1})$ that
\emph{does not require advance knowledge of $N$} and produces an accuracy \emph{certificate}:
the computable bound~\eqref{eq:excessgapbound}. In the strongly
convex case, the same technique gives $\Oh(\eps^{-1/2})$. A separate
line of work extends smoothing to symmetric functions of the eigenvalues
of symmetric matrices, yielding methods with complexity
$\Oh(\eps^{-1})$ for minimizing the maximum eigenvalue and
the spectral radius of a family depending linearly on
parameters~\cite{Nes2007sdpsmooth}.

\subsection{Current state}

The smoothing technique has become a standard part of the optimizer's toolkit. An alternative and also a complement is A.\,S.~Nemirovski's Mirror Prox method~\cite{Nemirovski2004} for monotone variational inequalities, which gives the same bound $\Oh(\eps^{-1})$ for convex-concave saddle-point problems; this was followed by the Chambolle--Pock primal-dual algorithm~\cite{ChambollePock2011}, which became a working tool in image processing. Optimality of the bound itself has been established: lower bounds for the class of bilinear saddle-point problems show that $\Oh(\nrm{A}R^2/\eps)$ cannot be improved~\cite{OuyangXu2021}. Thus smoothing completes the picture: exploiting problem structure makes it possible to attain the dependence $\eps^{-1}$.

\section{Primal-dual methods (2005 -- 2011)}
\label{sec:pd}

\subsection{Why duality is needed as a computational tool}

In \S\,\ref{sec:traffic}, we encountered a situation typical of
applications but atypical of textbooks. Problem~\eqref{eq:sddual}
is formulated in dual variables $t$ (travel times), while
the quantities actually needed are the \emph{primal} variables $f$ (flows), available only
as Lagrange multipliers. In~\eqref{eq:threestage}, both are needed:
the origin--destination matrix, arc flows, \emph{and} district
potentials. The same is true in entropy-linear programming,
optimal transport, and pricing problems.

A second, equally important reason is the \emph{accuracy certificate}. A bound
of the form $f(x_N)-f^{*}\le\eps$ is useless in practice if $f^{*}$
is unknown: we do not know when to stop. If, however, a method
produces a \emph{pair} $(\hat x_N,\hat\lambda_N)$ of approximate primal and dual
solutions, the duality gap can be \emph{computed} and
bounds both errors from above. Methods producing such a pair together with
a gap bound are called \emph{primal-dual}; the a posteriori
gap bound itself is called a \emph{primal-dual
certificate}. It is important to distinguish the two: a certificate is a form of bound, not
an algorithmic property, and can be constructed for methods from quite different
families (Remark~\ref{rem:pdeverywhere}).

A systematic theory of such methods was developed
in~\cite{Nes2009pd} (a CORE preprint of 2005) and, independently and in a more
general setting, by A.\,S.~Nemirovski, S.~Onn, and
U.~Rothblum~\cite{NemOnnRothblum2010}. Together, these two papers provide
the most complete modern account of primal-dual numerical
methods for problems with convex structure.

\subsection{The mechanism: averaging the inner solutions}
\label{sec:pdmech}

We examine the mechanism on the model problem
\begin{equation}\label{eq:pdprob}
   f^{*}=\min_{x\in Q}\ \{f(x):\ Ax=b\},
\end{equation}
where $f$ is convex and $Q$ is a compact convex set. The objective of the dual problem,
\begin{equation}\label{eq:dualfun}
   \phi(\lambda)=\min_{x\in Q}
   \bigl\{f(x)+\dotp{\lambda}{Ax-b}\bigr\}
   =f\bigl(x(\lambda)\bigr)+\dotp{\lambda}{Ax(\lambda)-b}
\end{equation}
is concave, and its supergradient is
$\grad\phi(\lambda)=Ax(\lambda)-b$, where $x(\lambda)$ solves
the inner problem. Crucially, \emph{computing one
supergradient of the dual function automatically supplies a point
$x(\lambda)$ in the primal space} (when the functional dependence $x(\lambda)$ is known). It turns out that these points simply need to be averaged, using the same weights as the method.

\begin{lemma}[pointwise primal-dual inequality]\label{lem:pd}
Let $\lambda_0,\dots,\lambda_N$ be arbitrary points,
$\alpha_0,\dots,\alpha_N>0$, $S_N=\sum_{k=0}^N\alpha_k$, and
\[
   \hat x_N=\frac{1}{S_N}\sum_{k=0}^{N}\alpha_k\,x(\lambda_k)\ \in Q .
\]
Then, for every $\lambda$,
\begin{equation}\label{eq:pdkey}
   f(\hat x_N)+\dotp{\lambda}{A\hat x_N-b}
   \ \le\ \Phi_N(\lambda)
   :=\frac{1}{S_N}\sum_{k=0}^{N}\alpha_k
   \Bigl[\phi(\lambda_k)
   +\dotp{\grad\phi(\lambda_k)}{\lambda-\lambda_k}\Bigr].
\end{equation}
\end{lemma}

\begin{proof}
Substituting $\phi(\lambda_k)$ and
$\grad\phi(\lambda_k)=Ax(\lambda_k)-b$ from~\eqref{eq:dualfun}, we obtain
\[
   \Phi_N(\lambda)=\frac{1}{S_N}\sum_k\alpha_k
   \Bigl[f\bigl(x(\lambda_k)\bigr)
   +\dotp{\lambda}{Ax(\lambda_k)-b}\Bigr] .
\]
The second term is linear in $x(\lambda_k)$ and gives exactly
$\dotp{\lambda}{A\hat x_N-b}$; by convexity of $f$, the first is no less than
$f(\hat x_N)$.
\end{proof}

\begin{remark}[why a bounded set $\Lambda$ is needed]
\label{rem:boundedLambda}
Inequality~\eqref{eq:pdkey} is \emph{pointwise}, and one cannot maximize it over
all of $\R^m$: if $A\hat x_N\ne b$, then
$\sup_{\lambda\in\R^m}\{f(\hat x_N)+\dotp{\lambda}{A\hat x_N-b}\}=+\infty$.
Thus a \emph{bounded} set of dual variables must be fixed
in advance:
\[
   \Lambda=\bigl\{\lambda:\ \nrm{\lambda}_{*}\le R_\lambda\bigr\}
\]
and the method must control the \emph{maximum}, rather than the minimum, of
$\Phi_N$ over $\Lambda$. This is precisely how the standard accuracy
certificate is constructed~\cite{NemOnnRothblum2010}.
\end{remark}

\begin{proposition}[accuracy certificate]\label{pr:certificate}
Suppose that strong duality holds for~\eqref{eq:pdprob} and
there exists an optimal multiplier $\lambda^{*}$ with
$\nrm{\lambda^{*}}_{*}\le R_\lambda/2$; set
$\Lambda=\{\nrm{\lambda}_{*}\le R_\lambda\}$. Suppose that the method
applied to maximizing $\phi$ over $\Lambda$ guarantees
\begin{equation}\label{eq:methodguar}
   \max_{\lambda\in\Lambda}\Phi_N(\lambda)
   \ -\ \frac1{S_N}\sum_{k=0}^{N}\alpha_k\phi(\lambda_k)
   \ \le\ \eps_N .
\end{equation}
Then
\begin{equation}\label{eq:certificate}
   \bigl[f(\hat x_N)-f^{*}\bigr]
   +R_\lambda\nrm{A\hat x_N-b}\ \le\ \eps_N,
\end{equation}
and, in particular, $f(\hat x_N)-f^{*}\le\eps_N$ and
$\nrm{A\hat x_N-b}\le2\eps_N/R_\lambda$. The left-hand side of~\eqref{eq:methodguar}
is computed during the run and therefore serves as a stopping criterion.
\end{proposition}

\begin{proof}
Maximize~\eqref{eq:pdkey} over $\lambda\in\Lambda$; the left-hand side
becomes $f(\hat x_N)+R_\lambda\nrm{A\hat x_N-b}$ by the definition
of the dual norm. On the right, using~\eqref{eq:methodguar} and
$\phi(\lambda_k)\le\phi(\lambda^{*})=f^{*}$, we obtain
$\max_{\lambda\in\Lambda}\Phi_N(\lambda)\le f^{*}+\eps_N$. This is~\eqref{eq:certificate}.
Next, strong duality implies the exact-penalty inequality
$f(x)+\nrm{\lambda^{*}}_{*}\nrm{Ax-b}\ge f^{*}$ for every $x\in Q$;
subtracting it from~\eqref{eq:certificate}, we find
$(R_\lambda-\nrm{\lambda^{*}}_{*})\nrm{A\hat x_N-b}\le\eps_N$, whence
$\nrm{A\hat x_N-b}\le2\eps_N/R_\lambda$; then~\eqref{eq:certificate}
also yields $f(\hat x_N)-f^{*}\le\eps_N$.
\end{proof}

We emphasize that no additional work beyond solving
the dual problem is required: the points $x(\lambda_k)$ are already computed
at each iteration when evaluating the supergradient.

\subsection{Dual averaging}

The specific scheme in~\cite{Nes2009pd} is \emph{dual
averaging}. Let $Q$ be equipped with a prox-function $d$, and let
$\pi_\beta(s)=\argmin_{x\in Q}\{-\dotp{s}{x-x_0}+\beta d(x)\}$. The scheme is
\begin{equation}\label{eq:da}
   s_0=0;\qquad
   g_k\in\partial f(x_k),\quad s_{k+1}=s_k+\alpha_kg_k,\quad
   x_{k+1}=\pi_{\beta_{k+1}}(-s_{k+1}),
\end{equation}
where $g_k$ is a subgradient returned by the oracle, $\alpha_k>0$ are weights, and
$\beta_{k+1}\ge\beta_k$ are parameters (the notation is consistent
with~\eqref{eq:pdkey}: the same weights are used for averaging). The distinction from the classical
subgradient method is fundamental: subgradients are accumulated
\emph{in the dual space} (as a sum of linear functionals), and
a single prox-mapping takes us back to the primal space. The paper
shows that \emph{two} control sequences are needed---one
for movement in the primal space and the other in the dual space---and that
this is precisely what eliminates superfluous logarithmic factors in
complexity bounds.

The quality of the sequence $\{x_i\}$ is described by the \emph{gap function}
\begin{equation}\label{eq:gapfun}
   \delta_k(D)=\max_x\Bigl\{\sum_{i=0}^k\alpha_i
   \dotp{g_i}{x_i-x}:\ x\in Q,\ d(x)\le D\Bigr\},
   \qquad D\ge0,
\end{equation}
which is \emph{computable} because $Q$ is simple. If a solution
$x^{*}$ exists, then for $D\ge d(x^{*})$, the quantity $\delta_k(D)$
is nonnegative regardless of which sequences are used
to define it, and, normalized by $S_k=\sum_{i\le k}\alpha_i$, it
bounds the residual of the averaged point from above~\cite{Nes2009pd}. This yields
reliable stopping criteria. Scheme~\eqref{eq:da} allows
unbounded $Q$ and readily extends to composite and stochastic
settings, minimax problems, problems with functional
constraints, and variational inequalities.

\practice It was precisely~\eqref{eq:da} that directly preceded
regularized dual averaging in machine
learning~\cite{Xiao2010}. Follow-the-Regularized-Leader constructions in
online learning and adaptive methods such as AdaGrad~\cite{DuchiHazanSinger2011}
are \emph{related} regularized dual schemes.

\subsection{The barrier subgradient method}
\label{sec:bsm}

A separate and, in our view, underappreciated construction is
the \emph{barrier subgradient method}~\cite{Nes2011barrier}, in which
the gap function~\eqref{eq:gapfun} is smoothed using not an arbitrary
prox-function, but a \emph{self-concordant barrier} from
\S\,\ref{sec:ipm}. Let a closed convex set $Q$ contain no lines and
be equipped with a $\nu$-self-concordant barrier $F$, let $\hat P$ be another
convex set, and suppose their intersection $P=\hat P\cap Q$ is bounded;
write $P_0=\hat P\cap\operatorname{int}Q$, and let
$x_0=\argmin_{x\in P_0}F(x)$ be the \emph{constrained analytic center}.
Instead of~\eqref{eq:gapfun}, use a smooth approximation of the support
function of $P$:
\begin{equation}\label{eq:bsmsupport}
   U_\beta(s)=\max_{x\in\hat P}
   \bigl\{\dotp{s}{x-x_0}-\beta\bigl[F(x)-F(x_0)\bigr]\bigr\},
\end{equation}
where $\beta>0$ is the smoothing parameter.
The resulting method solves $\max\{f(x):x\in P\}$ for concave $f$ and,
particularly interestingly, allows positive concave
functions to be maximized \emph{to relative accuracy}, with complexity independent
of the properties of the particular problem, in particular the size of
the domain~\cite{Nes2011barrier}. Applications include computing fractional covers,
maximum multicommodity flows, minimax problems with
nonnegative coefficients, semidefinite relaxation of Boolean
quadratic problems (\S\,\ref{sec:pi2}), and ``online'' optimization
as an alternative to a stochastic formulation.

Thus the barrier subgradient method links two parts of the theory: the geometric model of \S\,\ref{sec:ipm} and the algebraic model of
\S\,\ref{sec:smooth}. We also mention a separate direction here: \emph{optimization to relative accuracy} for homogeneous problems, where complexity bounds are almost independent of the input data~\cite{Nes2009relative,Nes2008rounding}. There an additional improvement is obtained by first computing a suitable Euclidean norm through ellipsoidal approximation of the feasible set.

\subsection{Dual extrapolation for variational inequalities}

The most general class of problems with convex structure is that of monotone
variational inequalities: find $x^{*}\in Q$ such that
$\dotp{g(x)}{x-x^{*}}\ge0$ for every $x\in Q$, where the operator $g$
is monotone. Convex minimization, finding
saddle points of convex-concave functions, and equilibrium computation all fit
this form.
For these problems,~\cite{Nes2007dualextrap} proposed a \emph{dual
extrapolation} step: the linear function representing the problem model
is updated according to a rule dual to the extragradient step.
The scheme gives a bound of $\Oh(\eps^{-1})$ for a Lipschitz monotone operator and
$\Oh(\eps^{-2})$ for an operator with bounded variation, hence in
particular for nonsmooth problems; the stopping criterion is again
computable.

\subsection{Applications}

\begin{example}[entropy-linear programming]\label{ex:elp}
Problem~\eqref{eq:entropy} has the form
\[
   \min_{d\in\Delta}\Bigl\{\sum_{i,j}d_{ij}\ln d_{ij}
   +\beta\dotp{T}{d}\Bigr\}\quad\text{subject to}\quad
   \sum_jd_{ij}=l_i,\ \ \sum_id_{ij}=w_j .
\]
The primal problem has dimension $n^2$ (for a large metropolitan area,
$10^4$\,--\,$10^6$), while the dual has dimension $2n$. Thus one should solve
the dual, but the answer is needed in primal variables. The dual function
is evaluated explicitly (the inner minimum for an entropy objective
has a closed form and gives a Gibbs distribution), and
primal-dual schemes produce both the district potentials
$\lambda^L,\lambda^W$ and the origin--destination matrix itself. A systematic
account is given in~\cite{GasnEtAlELP2016}, which also discusses the relation to
the classical balancing method (Sinkhorn's method), which in these terms
is a special case.
\end{example}

\begin{example}[traffic equilibria]\label{ex:pdtraffic}
The dual problem~\eqref{eq:sddual} and its multistage
variant~\eqref{eq:threestage} provide a benchmark example where
a primal-dual approach is necessary: flows $f$ are available only as
Lagrange multipliers. Corresponding schemes for mixed flow-assignment
models are constructed in~\cite{GasnGasnNes2018} (which also provides
a fairly general and complete account of primal-dual techniques and their
application to transportation problems), and for the nonsmooth dual problem with
a shortest-path oracle in~\cite{Baimurzina2019,KubentayevaGasnikov2021}.
An accessible account of primal-dual techniques oriented toward these
applications is given in \S\,4 of the textbook~\cite{GasnBook}.
\end{example}

\subsection{What changed}

Before~\cite{Nes2005excess,Nes2009pd,NemOnnRothblum2010},
duality in numerical optimization was used mainly as
a way to \emph{reformulate} a problem: to pass to a problem of lower
dimension or simpler structure. Subsequently, it became
a tool of the method itself: dual variables are constructed
\emph{together} with primal ones, the gap between them serves as a computable stopping
criterion, and averaging the inner solutions is a free by-product
of supergradient evaluation. As noted in
Remark~\ref{rem:pdeverywhere}, viewed appropriately,
methods from quite different families also turn out to be primal-dual,
for example, cutting-plane methods. For applications in which
answers in both primal and dual variables are meaningful (traffic
flows and delays, origin--destination demands and district potentials, transportation plans
and prices), this changed the very formulation of the computational problem.

%% file: en/r8_cubic_gn.tex
\section{Cubic regularization of Newton's method (2006)}
\label{sec:cubic}

\subsection{Background: Newton's method without global guarantees}

Newton's method has been known for more than three hundred years and has the remarkable
property of local quadratic convergence. Its global behavior,
however, is poor: it can diverge, and for nonconvex functions it can
converge to a saddle point or a maximum, since the Newton
direction $-[\grad^2f(x)]^{-1}\grad f(x)$ need not be
a descent direction at all if the Hessian is not
positive definite.
By 2005, the standard remedies---backtracking, trust
regions, and spectral shifts---had a rich history of practical use, but lacked
complexity bounds of the form ``how many iterations are needed for accuracy $\eps$?''

The paper by Yu.\,E.~Nesterov and B.\,T.~Polyak~\cite{NesPolyak2006}---
a joint work of student and teacher---filled this gap. Let
$f$ be twice continuously differentiable with a Lipschitz-continuous Hessian:
\begin{equation}\label{eq:hesslip}
   \nrm{\grad^2f(x)-\grad^2f(y)}\ \le\ L_2\nrm{x-y}
   \qquad\forall x,y\in\R^n,
\end{equation}
whence the standard Taylor remainder bound gives
\begin{equation}\label{eq:taylor3}
   \bigl|f(y)-f(x)-\dotp{\grad f(x)}{y-x}
   -\tfrac12\dotp{\grad^2f(x)(y-x)}{y-x}\bigr|
   \ \le\ \frac{L_2}{6}\nrm{y-x}^3 .
\end{equation}
Inequality~\eqref{eq:taylor3} suggests the model
\begin{equation}\label{eq:cubmodel}
   \Omega_{M_2}(x;h)=f(x)+\dotp{\grad f(x)}{h}
   +\tfrac12\dotp{\grad^2f(x)h}{h}+\frac{M_2}{6}\nrm{h}^3 ,
\end{equation}
which, for $M_2\ge L_2$, is a \emph{global upper model}:
$f(x+h)\le\Omega_{M_2}(x;h)$ for all $h$. The cubic regularization method is
\begin{equation}\label{eq:cubstep}
   h_k\in\Argmin_{h\in\R^n}\Omega_{M_2}(x_k;h),\qquad x_{k+1}=x_k+h_k .
\end{equation}

The idea of adding a cubic term had been proposed earlier
(A.~Griewank, an unpublished technical
report of 1981~\cite{Griewank1981}), but it was
in~\cite{NesPolyak2006} that a complexity theory was developed and, no less
importantly, the \emph{computability} of the step was established.

\begin{remark}[solvability of the auxiliary problem; the hard case]\label{rem:cubsub}
The function~\eqref{eq:cubmodel} is nonconvex if
$\grad^2f(x)\not\succeq0$.
Nevertheless, it is coercive (the cubic term dominates as
$\nrm{h}\to\infty$), so its \emph{global} minimum exists and
has an explicit characterization: a point $h$ is optimal if and only if
\begin{equation}\label{eq:cubopt}
   \Bigl(\grad^2f(x)+\frac{M_2}{2}\nrm{h}I\Bigr)h=-\grad f(x),
   \qquad \grad^2f(x)+\frac{M_2}{2}\nrm{h}I\succeq0 .
\end{equation}
Writing $r=\nrm{h}$ yields a scalar (``secular'') equation
in $r$, solved by a one-dimensional search; each search step requires
one Cholesky factorization or several Lanczos iterations. This is
exactly the same computational work as in the trust-region
method---but now with a complexity bound.

As in trust-region theory, there is a so-called
\emph{hard case}: if the shifted matrix in~\eqref{eq:cubopt}
is singular and the gradient is orthogonal to the eigenspace
corresponding to the smallest eigenvalue, the solution cannot
be expressed through the inverse matrix. The global minimum remains
computable, but a pseudoinverse and an additional component in
the eigenspace of the smallest eigenvalue are required; for details, see the theory of
the cubic and trust-region subproblems in~\cite{CartisGouldToint2011}.
\end{remark}

\begin{figure}[htbp]
\centering
\includegraphics[width=0.88\textwidth]{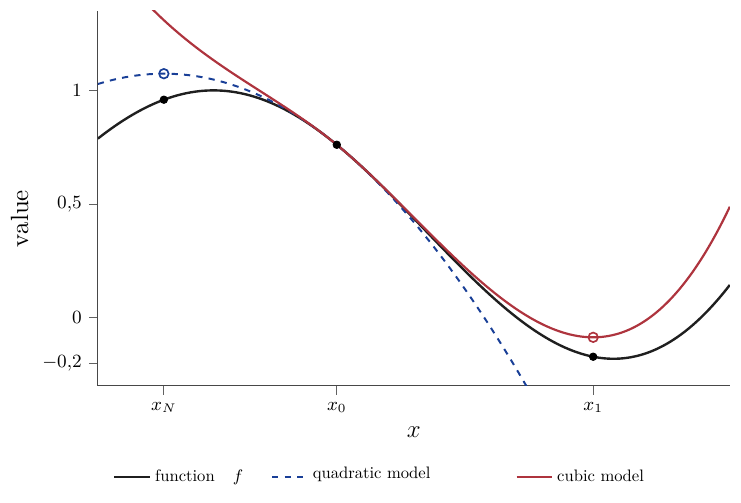}
\caption{Quadratic and cubic models of the function
$f(x)=\cos x+x^2/10$ at $x_0=0{,}8$ with $M_2=L_2=1$. Since
$f''(x_0)=-0{,}49671<0$, the quadratic model $q$ is concave: its stationary
point $x_N=-0{,}32210$ corresponds to the Newton step and gives
$f(x_N)=0{,}95895>f(x_0)=0{,}76071$. The cubic upper model
$\Omega(x)=q(x)+|x-x_0|^3/6$ has its minimum at $x_1=2{,}46351$, with
$f(x_1)=-0{,}17189<\Omega(x_1)=-0{,}08649<f(x_0)$. Open colored
circles mark the extrema of the models, and black dots mark the function
values. The upper bound is global, since
$|f^{(3)}(x)|=|\sin x|\le1$.}
\label{fig:cubic}
\end{figure}

The behavior of the quadratic and cubic models under negative curvature
is shown in Fig.~\ref{fig:cubic}.

\subsection{Convergence rate bounds}

Three basic inequalities do almost all the work.

\begin{lemma}\label{lem:cubkey}
Suppose that~\eqref{eq:hesslip} holds, $M_2 \ge L_2$, and $h$ is
a global minimizer of $\Omega_{M_2}(x;\cdot)$. Then
\begin{align}
   f(x)-f(x+h)&\ \ge\ \frac{M_2}{12}\nrm{h}^3, \label{eq:cubdecr}\\
   \nrm{\grad f(x+h)}&\ \le\ \frac{L_2+M_2}{2}\nrm{h}^2,
   \label{eq:cubgrad}\\
   \lambda_{\min}\bigl(\grad^2f(x+h)\bigr)
   &\ \ge\ -\Bigl(\frac{M_2}{2}+L_2\Bigr)\nrm{h}. \label{eq:cubcurv}
\end{align}
\end{lemma}

\begin{proof}
The first condition in~\eqref{eq:cubopt} implies
$\dotp{\grad f(x)}{h}=-\dotp{\grad^2f(x)h}{h}-\frac{M_2}{2}\nrm{h}^3$,
so
\[
   \Omega_{M_2}(x;h)-f(x)
   =-\tfrac12\dotp{\grad^2f(x)h}{h}-\tfrac{M_2}{3}\nrm{h}^3 .
\]
The second condition in~\eqref{eq:cubopt} gives
$\dotp{\grad^2f(x)h}{h}\ge-\frac{M_2}{2}\nrm{h}^3$, whence
$\Omega_{M_2}(x;h)-f(x)\le\frac{M_2}{4}\nrm{h}^3-\frac{M_2}{3}\nrm{h}^3
=-\frac{M_2}{12}\nrm{h}^3$; together with $f(x+h)\le\Omega_{M_2}(x;h)$ this
gives~\eqref{eq:cubdecr}. Next,
$\nrm{\grad f(x+h)-\grad f(x)-\grad^2f(x)h}\le\frac{L_2}{2}\nrm{h}^2$
by~\eqref{eq:hesslip}, while the first condition in~\eqref{eq:cubopt} yields
$\nrm{\grad f(x)+\grad^2f(x)h}=\frac{M_2}{2}\nrm{h}^2$; the triangle
inequality gives~\eqref{eq:cubgrad}. Finally,
$\grad^2f(x+h)\succeq\grad^2f(x)-L_2\nrm{h}I
\succeq-\frac{M_2}{2}\nrm{h}I-L_2\nrm{h}I$, which
is~\eqref{eq:cubcurv}.
\end{proof}

It is convenient to combine inequalities~\eqref{eq:cubgrad} and~\eqref{eq:cubcurv}.
For $M_2 = L_2$, set
\begin{equation}\label{eq:mudef}
   \mu_{L_2}(x)\ :=\ \max\Bigl\{\sqrt{\nrm{\grad f(x)}/L_2},\;
   -\tfrac{2}{3L_2}\lambda_{\min}\bigl(\grad^2f(x)\bigr)\Bigr\}\ \ge 0 .
\end{equation}
Then~\eqref{eq:cubgrad}--\eqref{eq:cubcurv} mean precisely that
\begin{equation}\label{eq:hmu}
   \nrm{h}\ \ge\ \mu_{L_2}(x+h),
\end{equation}
that is, the step length bounds the ``second-order stationarity residual'' at
the \emph{new} point from above. The quantity $\mu_{L_2}$ has the dimension of length;
$\mu_{L_2}(x)\le\eps$ is equivalent to
$\nrm{\grad f(x)}\le L_2\eps^{2}$ and
$\lambda_{\min}(\grad^2f(x))\ge-\tfrac32 L_2\eps$.

\begin{theorem}[\cite{NesPolyak2006}]\label{th:cubic}
\emph{Class:} $f$ is a twice continuously differentiable function (convexity is not assumed) whose Hessian
satisfies
condition~\eqref{eq:hesslip}, and $f^{\ast}$ is the global minimum
of the function. \emph{Oracle:}
gradient and Hessian. \emph{Method:}~\eqref{eq:cubstep} with $M_2 = L_2$.
\emph{Guarantee:} deterministic.
The method guarantees convergence to an approximately second-order stationary
point in the sense of
the measure~\eqref{eq:mudef}, so that, for every $k\geq 1$,
      \begin{equation}\label{eq:cubnonconv}
      \min\limits_{1\leq i \leq k} \mu_{L_2}(x_i) \leq
      \Bigl(\frac{12\bigl(f(x_0) - f^{\ast}\bigr)}{k L_2}
      \Bigr)^{1/3}.
      \end{equation}
Consequently, taking~\eqref{eq:mudef} into account, we obtain the explicit bound
$$\min\limits_{1\leq i \leq k} \|\nabla f(x_i)\|\ \le \Oh(k^{-2/3})$$
and simultaneously
$\min_{1\le i\le k}\bigl(-\lambda_{\min}(\nabla^2f(x_i))\bigr)
\leq \Oh(k^{-1/3})$.
For comparison, gradient schemes give the rate $\Oh(k^{-1/2})$
and only in the gradient norm.
\end{theorem}

The proof of~\eqref{eq:cubnonconv} follows from the elegant
observation
\begin{align*}
    f(x_0)-f^{\ast} &\geq \sum\limits_{i=0}^{k-1}
    \bigl(f(x_i) - f(x_{i+1})\bigr)\geq
    \sum\limits_{i=0}^{k-1}\tfrac{L_2}{12}\|h_i\|^3
    \\ &\geq \tfrac{L_2}{12}\sum\limits_{i=0}^{k-1}
      \mu_{L_2}(x_{i+1})^3\ \geq\;
      \tfrac{L_2 k}{12}\min\limits_{1\leq i\leq k}\mu_{L_2}(x_i)^3,
\end{align*}
where the second step uses the basic
inequality~\eqref{eq:cubdecr}, and the third uses
inequalities~\eqref{eq:cubgrad} and~\eqref{eq:cubcurv} (in the form~\eqref{eq:hmu}).

For simplicity of presentation we have fixed $M_2 = L_2$, but the original
paper gives results
for general $M_2 \geq L_2$, and this generality will later play an important role in
the development of accelerated
and higher-order methods (\S\,\ref{sec:tensor}).

It should be emphasized that, compared with
$\Oh(\eps^{-3/2})$ for the cubically
regularized Newton method, gradient descent on the same class
gives only $\Oh(\eps^{-2})$, and only for the condition $\nrm{\grad f}\le\eps$, without
distinguishing minima from saddle points.
Cubic regularization is both faster
and \emph{qualitatively better}: it automatically escapes regions of negative
curvature, guaranteeing approximate \emph{second}-order stationarity.
This is the first result of this kind in the literature, and it is what made the question
``how quickly can one escape a saddle point?'' technically meaningful.

\begin{theorem}\label{th:cubic-cvx}
\emph{Class:} $f$ is twice differentiable and convex, its Hessian satisfies
condition~\eqref{eq:hesslip}, and $x^{\ast}$ is a global
minimizer of the function. \emph{Oracle:}
gradient and Hessian. \emph{Method:}~\eqref{eq:cubstep} with $M_2 \geq L_2$.
\emph{Guarantee:} deterministic.
\begin{enumerate}[label=\rm(\alph*),leftmargin=2.2em,itemsep=4pt]
\item \textbf{Nonaccelerated method \cite{NesPolyak2006}.}
    For the basic scheme~\eqref{eq:cubstep}, the standard global bound
      assumes that the initial sublevel set is bounded, with radius
      \[
         D=\sup\bigl\{\dist(x,X^{*}):\ f(x)\le f(x_0)\bigr\}
         < +\infty.
      \]
    Then the method converges at the rate
    \begin{equation*}
         f(x_k)-f^{*}=\Oh\bigl(\tfrac{L_2 D^3}{k^2}\bigr) \qquad k\ge1.
    \end{equation*}
\item \textbf{Accelerated method \cite{Nes2008accel}.}
      The accelerated variant admits a bound in terms of the
      initial distance itself, $R = \|x_0-x^{*}\|$: for $M_2 \geq 2L_2$,
      \begin{equation*}\label{eq:cubaccel}
         f(x_k)-f^{*}\ \le \Oh \left(\tfrac{L_2 R^3}{k^3}\right),
         \qquad k\ge1,
      \end{equation*}
      that is, $N(\eps)=\Oh\bigl((L_2 R^3/\eps)^{1/3}\bigr)$.
\end{enumerate}
\end{theorem}

We complete the picture with a local result: near a nondegenerate
minimum ($\grad^2f(x^{*})\succ0$), method~\eqref{eq:cubstep} retains
the quadratic convergence rate of the ordinary Newton method~\cite{NesPolyak2006,Doikov2022local},
so the global bounds (a)--(b) do not ``spoil'' the classical behavior.

The exponent $1/3$ is \emph{not} optimal
for second-order methods: the lower bound has exponent $2/7$
(Table~\ref{tab:map}), and it is attained by a different, substantially more
complicated construction (\S\,\ref{sec:tensor}).

\subsection{Gradient domination: superlinear convergence}

The advantages of cubic regularization are particularly striking under
the additional assumption of \emph{gradient domination},
introduced by B.\,T.~Polyak back in 1963~\cite{Polyak1963} (now known as
the Polyak--\L{}ojasiewicz condition):
\begin{equation}\label{eq:pl}
   f(x)-f^{*}\ \le\ \frac{1}{2\tau}\nrm{\grad f(x)}^2
   \qquad\forall x,\qquad \tau>0 .
\end{equation}
Condition~\eqref{eq:pl} does not require convexity and holds, for example, for
strongly convex functions (with $\tau=\mu$), as well as for many nonconvex
overdetermined regression problems.

\begin{theorem}\label{th:cubpl}
Suppose that $f$ satisfies~\eqref{eq:hesslip} and~\eqref{eq:pl}, and the method is given
by~\eqref{eq:cubstep} with $M_2 = L_2$. Write
$\delta_k=f(x_k)-f^{*}$. Then
\begin{enumerate}[label=\rm(\roman*),leftmargin=2.2em,itemsep=3pt]
\item \emph{(global bound)}
      $\displaystyle \delta_N\le\Bigl(\frac{\delta_0}{c\,N}\Bigr)^{4/3}$,
      where $c=\dfrac{L_2}{12}\Bigl(\dfrac{2\tau}{L_2^2}\Bigr)^{3/4}$;
\item \emph{(superlinear convergence of order $4/3$)}
      $\delta_{k+1}\le C\,\delta_k^{4/3}$, where
      $C=\dfrac{12^{4/3}L_2^{2/3}}{2\tau}$. In particular, once
      $\delta_k\le\tfrac12 C^{-3}$, attaining
      accuracy $\eps$ requires only
      $\Oh\bigl(\ln\ln(C^{-3}/\eps)\bigr)$ further iterations.
\end{enumerate}
\end{theorem}

\begin{proof}
From~\eqref{eq:pl} at $x_{k+1}$ and~\eqref{eq:cubgrad} with $M_2 = L_2$,
\begin{equation}\label{eq:plchain}
   \delta_{k+1}\ \le\ \frac{1}{2\tau}\nrm{\grad f(x_{k+1})}^2
   \ \le\ \frac{(L_2+M_2)^2}{8\tau}\nrm{h_k}^4
   \ =\ \frac{L_2^2}{2\tau}\nrm{h_k}^4 .
\end{equation}
\emph{(ii)} From~\eqref{eq:cubdecr} we have
$\nrm{h_k}^3\le\frac{12}{L_2}(\delta_k-\delta_{k+1})\le
\frac{12}{L_2}\delta_k$, so
$\nrm{h_k}^4\le(12\delta_k/L_2)^{4/3}$, and substitution
into~\eqref{eq:plchain} gives
$\delta_{k+1}\le\frac{L_2^2}{2\tau}\bigl(\frac{12}{L_2}\bigr)^{4/3}
\delta_k^{4/3}=C\delta_k^{4/3}$. Setting $u_k=C^3\delta_k$, we obtain
$u_{k+1}\le u_k^{4/3}$, that is,
$u_{k_0+N}\le u_{k_0}^{(4/3)^N}$: once $u_{k_0}\le\tfrac12$
(which occurs after a finite, explicitly bounded number of iterations by part~(i)),
attaining $u_{k_0+N}\le C^3\eps$ requires $\Oh(\ln\ln(\cdot))$ iterations.

\emph{(i)} From~\eqref{eq:plchain},
$\nrm{h_k}\ge\bigl(2\tau\delta_{k+1}/L_2^2\bigr)^{1/4}$, and together
with~\eqref{eq:cubdecr},
\[
   \delta_k-\delta_{k+1}\ \ge\ \frac{L_2}{12}
   \Bigl(\frac{2\tau\,\delta_{k+1}}{L_2^{2}}\Bigr)^{3/4}
   \ =\ c\,\delta_{k+1}^{3/4} .
\]
Summing over $k=0,\dots,N-1$ and using the monotonicity of $\delta_k$,
we obtain $\delta_0\ge\delta_0-\delta_N\ge c\,N\delta_N^{3/4}$, which
is equivalent to~(i).
\end{proof}

Thus, under gradient domination, the cubic Newton method has
superlinear convergence \emph{of order $4/3$}---less than quadratic,
but substantially greater than linear---\emph{without} assuming
convexity and \emph{without} localization: bound~(i) guarantees entry into
the region of superlinear convergence after a finite, explicitly bounded number of
iterations.

\hypoth Condition~\eqref{eq:pl} is often mentioned in connection with training
overdetermined (``overparameterized'') neural networks. Care is
needed here: the corresponding results are generally local and
depend on the architecture, initialization, and sample; see the review
in~\cite{KarimiNutiniSchmidt2016} and the discussion
in~\cite[\S\,4.4]{AvetEtAl2026}. Thus, the assertion that ``the
Polyak--\L{}ojasiewicz condition holds in neural network training''
should be regarded as a working hypothesis, not a theorem.

\subsection{Current state of the art}

The paper~\cite{NesPolyak2006} has generated an extensive literature.
C.~Cartis, N.~Gould, and Ph.~Toint~\cite{CartisGouldToint2011} developed the adaptive
cubic method ARC, in which the parameter $M_2$ is chosen during the computation and
inexact solution of the auxiliary problem is allowed.
For the nonconvex case, lower bounds of $\Om(\eps^{-3/2})$ have been established for second-order
methods~\cite{CarmonEtAl2020}, demonstrating the optimality of
Theorem~\ref{th:cubic}; we stress that these bounds refer to
a specific second-order model and a specific stationarity
criterion. Variants with
H\"older-continuous
Hessians~\cite{GrapigliaNes2017}, stochastic and subsampled
versions~\cite{Agafonov2023advancing}, schemes with ``lazy'' Hessians (in which the Hessian is not recomputed at
every iteration)~\cite{Doikov2023Lazy}, and quasi-Newton methods with a proven
superlinear rate have been developed. Here we should mention the series of papers by
A.~Rodomanov and Yu.\,E.~Nesterov~\cite{RodomanovNes2021}, which provided the first
explicit nonasymptotic superlinear convergence bounds for
Broyden-type and DFP methods.

\section{Modified Gauss--Newton scheme (2007)}
\label{sec:gn}

\subsection{Background}

Solving a system of nonlinear equations $F(x)=0$, where
$F\colon\R^n\to\R^m$, is one of the oldest problems in numerical mathematics.
The classical Gauss--Newton method linearizes $F$ and solves the
least-squares problem
\[
   h_k=\argmin_{h\in\R^n}\ \tfrac12\nrm{F(x_k)+F'(x_k)h}_2^2 ,\qquad
   x_{k+1}=x_k+h_k .
\]
It is known to converge quadratically in a neighborhood of a zero-residual solution
under regularity of the Jacobian, and to be capable of diverging outside such a
neighborhood. All the practical remedies---backtracking,
trust regions, and Levenberg--Marquardt regularization---
are reasonable, but did not provide global complexity bounds \emph{in terms of the required
accuracy}.

In~\cite{Nes2007gn}, Yu.\,E.~Nesterov proposed a modification that
at first seems cosmetic but changes everything: minimize not
the squared norm, but \emph{the norm itself},
\begin{equation}\label{eq:gnobj}
   \varphi(x)=\nrm{F(x)},
\end{equation}
and use the model
\begin{equation}\label{eq:gnmodel}
   \Phi_M(x;h)=\nrm{F(x)+F'(x)h}+\frac{M}{2}\nrm{h}^2 .
\end{equation}

\subsection{Why it works}

The function $h\mapsto\nrm{F(x)+F'(x)h}$ is the norm of an affine mapping,
and is therefore \emph{convex}. Hence $\Phi_M(x;\cdot)$ is strongly
convex with parameter $M$, and the auxiliary problem
$\min_{h\in\R^n}\Phi_M(x;h)$ is a well-posed convex problem (for
the Euclidean norm, a second-order cone programming
problem) that can be solved to any prescribed accuracy. This is the first
advantage of the scheme. The second is that the model is a \emph{global
upper model}.

\begin{proposition}\label{pr:gn}
Suppose that $\nrm{F'(x)-F'(y)}\le L\nrm{x-y}$ for all $x,y$. Then, for all
$x,h$,
\begin{equation}\label{eq:gnmaj}
   \varphi(x+h)\ \le\ \Phi_L(x;h),
\end{equation}
and the iterations $x_{k+1}=x_k+h_k$,
$h_k=\argmin_{h\in\R^n}\Phi_M(x_k;h)$ with $M\ge L$ satisfy
\begin{equation}\label{eq:gndecr}
   \varphi(x_k)-\varphi(x_{k+1})\ \ge\ \frac{M}{2}\nrm{h_k}^2 .
\end{equation}
Consequently,
$\min_{0\le k<N}\nrm{h_k}\le\sqrt{2\bigl(\varphi(x_0)
-\varphi^{*}\bigr)/(MN)}$; that is, obtaining a point with
$\nrm{h_k}\le\eps$ requires at most $\Oh(\eps^{-2})$ iterations.
\end{proposition}

\begin{proof}
The Lipschitz continuity of $F'$ implies
$\nrm{F(x+h)-F(x)-F'(x)h}\le\frac{L}{2}\nrm{h}^2$; the triangle
inequality for the norm gives~\eqref{eq:gnmaj}. Next,
$\Phi_M(x_k;\cdot)$ is strongly convex with parameter $M$, so
$\Phi_M(x_k;0)-\Phi_M(x_k;h_k)\ge\frac{M}{2}\nrm{h_k}^2$; but
$\Phi_M(x_k;0)=\varphi(x_k)$ and
$\varphi(x_{k+1})\le\Phi_M(x_k;h_k)$ by~\eqref{eq:gnmaj} and $M\ge L$.
Summing~\eqref{eq:gndecr} over $k=0,\dots,N-1$ completes
the proof.
\end{proof}

\begin{remark}[the meaning of $M\nrm{h_k}$]\label{rem:gnstat}
The quantity $M\nrm{h_k}$ plays the role of the gradient norm: it is the norm
of the \emph{prox-linear mapping}, and its vanishing is equivalent to
\emph{composite stationarity},
\begin{equation}\label{eq:compstat}
   0\in F'(x)^{*}\,\partial\nrm{F(x)} ,
\end{equation}
where $\partial\nrm{\cdot}$ is the subdifferential of the norm, and
$F'(x)^{*}$ is the adjoint operator. It is~\eqref{eq:compstat}, rather than
$\grad\varphi(x)=0$, that is the correct stationarity condition for
the nonsmooth function~\eqref{eq:gnobj}; the corresponding chain rule (with
equality) for a composition of a convex function and a smooth mapping is
a standard fact of convex analysis.
\end{remark}

Thus, a \emph{global} complexity bound of
$\Oh(\eps^{-2})$ has been obtained for an essentially nonconvex and nonsmooth problem. Near
a zero-residual solution, the quadratic convergence of the
classical scheme is retained---but, we emphasize, this requires not only a zero
residual but also regularity of the Jacobian at the solution (for example, full
column rank of $F'(x^{*})$) and local stability: a zero
residual alone is insufficient.

Two further remarks on the precision of the statements are in order. First, the auxiliary
problem is solved not ``exactly'' but to an arbitrary prescribed accuracy;
since the discussion concerns worst-case bounds, the allowable
inexactness of the inner solution must be specified explicitly (typically,
a relative accuracy of order $\nrm{h_k}^2$, which does not change the order of
the bound). Second, replacing the Euclidean norm in~\eqref{eq:gnobj} by
the $\ell_1$-norm (useful for robustness to outliers) preserves
the convexity of the outer function and hence the entire scheme, \emph{but} changes
the smoothness constants and dual norms: the correct wording is not ``without
changes'', but ``within the same prox-linear scheme with the appropriate
geometry''.

\subsection{What changed}

The paper~\cite{Nes2007gn} became the starting point for an entire line of research:
\emph{prox-linear} methods for composite problems of the form
$\min_{x\in\R^n}\ \psi\bigl(c(x)\bigr)$, where $\psi$ is convex (possibly nonsmooth)
and $c$ is a smooth but nonlinear mapping. The general scheme
\begin{equation}\label{eq:proxlin}
   x_{k+1}=\argmin_{y\in\R^n}\Bigl\{\psi\bigl(c(x_k)+c'(x_k)(y-x_k)\bigr)
   +\frac{M}{2}\nrm{y-x_k}^2\Bigr\}
\end{equation}
and its analysis develop precisely the construction~\eqref{eq:gnmodel}. Here
we should mention work by A.~Lewis and S.~Wright~\cite{LewisWright2016},
D.~Drusvyatskiy and C.~Paquette~\cite{DrusPaquette2019},
J.~Duchi and F.~Ruan~\cite{DuchiRuan2019} (robust phase
retrieval), and D.~Davis and D.~Drusvyatskiy~\cite{DavisDrus2019} (weakly convex
stochastic optimization). The range of applications has proved broad:
robust phase retrieval, matrix recovery from incomplete
measurements, and nonlinear least squares robust to outliers.

More recently, N.~Yudin and A.~Gasnikov~\cite{YudinGasnikov2026}
established uniform lower complexity bounds for methods using the oracle
$(F,F')$, with query and output points restricted to $x_0$ plus the span
of previously observed Jacobian rows (residual-only queries within this
span are uncounted). In sufficiently large dimension and the stated
parameter ranges, these bounds match the stationarity guarantees of
Nesterov's norm-model method and the method of three squares, as well as
their global-phase guarantees under uniform nondegeneracy.

%% file: en/r9_cd_zo.tex

\section{Coordinate descent methods and huge-scale problems
(2010--2014)}
\label{sec:cd}

\subsection{Background}

The idea of moving along coordinate directions in turn is as old as
numerical mathematics itself: the Gauss--Seidel method for linear systems,
Southwell's relaxation method, and coordinate descent in nonlinear
programming. Practitioners knew that the method often works well. Theorists had convergence results for particular problem classes, but there was no systematic theory of \emph{global convergence rate bounds}. Moreover, in 1973 M.~Powell~\cite{Powell1973} constructed a smooth (nonconvex) function of three variables for which cyclic coordinate descent converges to a nonoptimal point.

A major breakthrough in the analysis of coordinate descent methods came with Yu.\,E.~Nesterov's paper~\cite{Nes2012cd} (a 2010 CORE preprint). The first important step proposed in the paper was \emph{randomization}: the coordinate is chosen at random, eliminating the cycling that defeated the method in Powell's example and making analysis possible. The second was \emph{the right problem formulation}: one should compare not iteration counts but total computational work, since a coordinate iteration can be $n$ times (and for sparse problems much more) cheaper than a full gradient iteration. The contribution of~\cite{Nes2012cd} is properly stated as follows: it provided the first systematic theory of randomized coordinate descent, including global bounds, importance sampling, acceleration, and accounting for the full cost in huge-scale problems.

Acceleration deserves special attention. To our knowledge, the possibility of accelerating coordinate descent methods had not been considered before~\cite{Nes2012cd}. The paper presents an acceleration scheme attaining the convergence rate $\Oh \bigl(\frac{1}{k^2}\bigr)$ typical of accelerated methods.

\subsection{The method and its bound}

Let $f$ be convex and differentiable on $\R^n$, and let its gradient be coordinatewise Lipschitz continuous: for all $x\in\R^n$, $h\in\R$, $i=1,\dots,n$,
\begin{equation}\label{eq:cdlip}
   \bigl|\grad_if(x+he_i)-\grad_if(x)\bigr|\ \le\ L_i\,|h| ,
\end{equation}
where $\grad_if=\partial f/\partial x^{(i)}$ and $e_i$ is a basis vector.
The coordinate descent method with random coordinate selection is
\begin{equation}\label{eq:cd}
   i_k\sim p,\qquad
   x_{k+1}=x_k-\frac{1}{L_{i_k}}\grad_{i_k}f(x_k)\,e_{i_k} .
\end{equation}

By~\eqref{eq:cdlip}, for the selected coordinate $i$, step~\eqref{eq:cd} gives $f(x_{k+1})\le f(x_k)-\frac{1}{2L_i}\bigl(\grad_if(x_k)\bigr)^2$; averaging over a uniform $i$ and taking $L_i\le L$, we obtain
\begin{equation}\label{eq:cdonestep}
   \Exp_{i}\bigl[f(x_{k+1})\bigr]\ \le\
   f(x_k)-\frac{1}{2nL}\nrm{\grad f(x_k)}_2^2 .
\end{equation}
This is exactly the inequality yielded by gradient descent for an $L$-smooth function, but with constant $nL$ instead of $L$: \emph{coordinate descent is gradient descent with a Lipschitz constant increased by a factor of $n$}.

Passing to a bound on the residual requires care in normalizing the initial error. The standard derivation uses $f(x)-f^{*}\le\dotp{\grad f(x)}{x-x^{*}}\le\nrm{\grad f(x)}\,\dist(x,X^{*})$; that is, the distance to the solution set is controlled not only at the initial point but at every point of the corresponding sublevel set. Thus the bound involves
\begin{equation}\label{eq:Rlev}
   R_{\mathrm{lev}}=\sup\bigl\{\dist(x,X^{*}):\ f(x)\le f(x_0)\bigr\},
\end{equation}
and this quantity is assumed finite. Method~\eqref{eq:cd} is almost surely monotone in function value, so this normalization is legitimate; the equality $R_{\mathrm{lev}}=\nrm{x_0-x^{*}}$ does not hold in general. (The original source~\cite{Nes2012cd} uses specially defined radii in weighted norms depending on $\{L_i\}$, but we give a simplified version.)

\begin{theorem}[\cite{Nes2012cd}]\label{th:cd}
\emph{Class:} $f$ is convex,~\eqref{eq:cdlip} holds, the set
$\{x\in\R^n:\ f(x)\le f(x_0)\}$ is bounded, and $R_{\mathrm{lev}}<\infty$.
\emph{Oracle:} one partial derivative $\grad_if(x)$ per iteration.
\emph{Guarantee:} in expectation. Then, for method~\eqref{eq:cd} with uniform
coordinate selection and $L_i\le L$,
\begin{equation}\label{eq:cdrate}
   \Exp\bigl[f(x_k)\bigr]-f^{*}\ =\
   \Oh\!\left(\frac{n\,L\,R_{\mathrm{lev}}^2}{k}\right),
\end{equation}
and if, in addition, $f$ is strongly convex with parameter $\mu>0$, then
$\Exp[f(x_k)]-f^{*}\le\bigl(1-\frac{\mu}{nL}\bigr)^{k}(f(x_0)-f^{*})$.
If coordinates are selected with probabilities $p_i=L_i/S_1$,
$S_1=\sum_jL_j$, the quantity $nL$ in the bounds is replaced by $S_1$.
\end{theorem}

\subsection{When coordinate descent really gains an advantage}

Formally, we have lost a factor of $n$ in the iteration count. Any gain can come only from the cost per iteration, and here a qualification often omitted is necessary: \emph{a single partial derivative is not always cheaper than the full gradient}. If the function is specified by an arbitrary program, then, by Theorem~\ref{th:ad}, the full gradient costs $\Oh(T)$---the same as one function value, and hence the same as one partial derivative. There is no gain.

A gain arises when there is \emph{structure} that allows one to
maintain a residual or another local aggregate. Let us examine the canonical example (cf.~\cite[pp.~204--205]{GasnBook}).

\begin{example}[a sparse quadratic problem]\label{ex:cdsparse}
Let
\[
   f(x)=\tfrac12\dotp{Ax}{x}-\dotp{b}{x},\qquad A\succeq0,
\]
and let each row and each column of $A$ have at most $s$ nonzero entries, so that $\mathrm{nnz}(A)\le ns$. Then:
\begin{itemize}[leftmargin=1.6em,itemsep=2pt]
\item the full gradient $\grad f(x)=Ax-b$ costs $\Oh(ns)$ operations;
\item one partial derivative $\grad_if(x)=\dotp{A_i}{x}-b_i$ costs
      $\Oh(s)$ operations, where $A_i$ is the $i$th row;
\item step~\eqref{eq:cd} changes exactly one coordinate of $x$, so if
      the vector $r=Ax-b$ is maintained, updating it after the step costs
      $\Oh(s)$ operations (only $s$ components change).
\end{itemize}
Thus a coordinate descent iteration is $n$ times cheaper than a gradient iteration, while the number of iterations is $n$ times larger. The methods seem \emph{comparable}. But the coordinatewise constants differ. For a quadratic function, $L_i=A_{ii}$ and $S_1=\Tr A$, whereas for the full gradient $L=\lambda_{\max}(A)$. We always have $\Tr A\le n\lambda_{\max}(A)$, and the inequality can be very loose.

Take an extreme case: $A=\diag(n,1,1,\dots,1)$ plus a sparse
correction, so that $\lambda_{\max}(A)\approx n$ and
$\Tr A\approx2n$. Then
\[
\begin{aligned}
   \underbrace{\Oh\!\Bigl(\frac{\lambda_{\max}R^2}{\eps}\Bigr)
   \cdot\Oh(ns)}_{\text{gradient descent}}
   &=\Oh\!\Bigl(\frac{n^2sR^2}{\eps}\Bigr),\\[6pt]
   \underbrace{\Oh\!\Bigl(\frac{\Tr A\cdot R^2}{\eps}\Bigr)
   \cdot\Oh(s)}_{\text{coordinate descent with importance sampling}}
   &=\Oh\!\Bigl(\frac{nsR^2}{\eps}\Bigr),
\end{aligned}
\]
so coordinate descent gains a factor of $n$. The reason is simple: gradient descent must take a step determined by the worst direction, whereas coordinate descent uses its own step $1/L_i$ along each coordinate.
\end{example}

There are also three considerations not directly reducible to the convergence rate:
\begin{itemize}[leftmargin=1.6em,itemsep=2pt]
\item for large $n$, storing several dense full-dimensional
      vectors can become the dominant memory constraint;
\item coordinate steps naturally admit parallel and asynchronous execution;
\item there is an \emph{accelerated} coordinate descent method in which $S_1=\sum_iL_i$ is replaced by $\bigl(\sum_i\sqrt{L_i}\bigr)^2$, and the square root of the condition number appears~\cite{NesStich2017}.
\end{itemize}

\subsection{Separating complexities across variable blocks}

The paper by Nesterov and Stich~\cite{NesStich2017} deserves special mention because it studies accelerated coordinate descent specifically for \emph{structured} problems and reveals the following effect: for a problem partitioned into variable blocks with different constants, an optimal method should access each block \emph{at its own frequency}---precisely the idea of complexity separation that we already encountered in \S\,\ref{sec:fgm} for saddle-point problems.

The clearest modern realization of this idea is the
\emph{min-min} problem. Consider
\[
   \min_{x\in\R^{d_x}}\min_{y\in\R^{d_y}}\ f(x,y),
\]
where $f$ is smooth, $\mu_x$-strongly convex in $x$, and $\mu_y$-strongly
convex in $y$; let $\kappa_x$ and $\kappa_y$ be the corresponding block condition numbers. The classical accelerated Nesterov method requires $\Oh\bigl(\sqrt{\max\{\kappa_x,\kappa_y\}}\,\ln(1/\eps)\bigr)$ evaluations of \emph{each} of the partial gradients $\grad_xf$ and $\grad_yf$. Building on~\cite{NesStich2017}, the paper~\cite{KovGasMal2025}
constructed a method requiring only
\[
\begin{gathered}
   \Oh\bigl(\sqrt{\kappa_x}\,\ln(1/\eps)\bigr)
   \ \text{evaluations of }\grad_xf
   \quad\text{and}\\
   \Oh\bigl(\sqrt{\kappa_y}\,\ln(1/\eps)\bigr)
   \ \text{evaluations of }\grad_yf,
\end{gathered}
\]
and these bounds are optimal. In applications, one often has $\kappa_x\gg\kappa_y$, while evaluating $\grad_yf$ is much cheaper than evaluating $\grad_xf$; the gain is then substantial.

The accelerated coordinate scheme~\cite{NesStich2017}, selecting a block with
probability proportional to $\sqrt{L_i}$, also extends to
a version of generalized smoothness related to
$(L_0,L_1)$-smoothness (\S\,\ref{sec:l0l1}), leading to an unexpected application.
In~\cite{PodlipnovaEtAl2026}, the smoothness constant of a block $B$ may grow
with the objective residual:
$\grad^2_{BB}f(x)\preceq\bigl(a_B+b_B[f(x)-f^{*}]\bigr)I$. If
blocks are selected with probabilities proportional to $\sqrt{a_B}+\sqrt{b_B\Delta}$, where
$\Delta$ is the current estimate of the residual, then an accelerated block method
finds an $\eps$-solution with probability at least $1-\beta$ in
\[
   \Ot\Bigl(\beta^{-1/2}R\Bigl(\frac{\sum_B\sqrt{a_B}}{\sqrt\eps}
   +\sum_B\sqrt{b_B}\Bigr)\Bigr)
\]
iterations, where $R$ is the distance from the initial point to a solution. The first
term is exactly the accelerated coordinate descent bound;
the second is the price of generalized smoothness, independent of $\eps$ (cf.\ the additive
structure of the bounds in \S\,\ref{sec:l0l1}). The application is to dual problems of
entropy-regularized linear programming, including
optimal transport (Example~\ref{ex:game}). Instead of the smooth log-sum-exp function,
which ``couples'' all the variables (Example~\ref{ex:lse}),
one can use a sum of exponentials: it splits into sparse terms
suited to block steps, but is smooth only in the generalized sense---hence
the question in the paper's title: ``Smooth or separable?''.

\subsection{Huge-scale problems}

The logical continuation is the paper~\cite{Nes2014huge}, which also introduced
the term \emph{huge-scale}: problems in which even a single pass over the data
is too expensive, so the cost per iteration must depend not on the size of
the problem but on \emph{local} sparsity.

For example, let $f(x)=\max_{1\le j\le m}\dotp{a_j}{x}$ with sparse
$a_j$. A subgradient step changes only a few coordinates of $x$; hence
the inner products $\dotp{a_j}{x}$ can be updated
incrementally, while the maximum over $j$ can be maintained in a tree (heap). If
each column and each row has at most $q$ nonzero entries, then
the cost of one iteration is $\Oh(q^2\log m)$ instead of
$\Oh(\mathrm{nnz}(A))$. Yu.\,E.~Nesterov showed that this technique---
``sparse updating'' with a tree data structure---applies to
an entire family of subgradient and coordinate descent schemes, and performed computations
for problems with $n$ of order $10^9$. This restructures not the method but its
\emph{implementation}: a combinatorial data structure becomes part of
the complexity analysis on an equal footing with the iteration bound.

Two applications deserve special mention.

\smallskip
\emph{Optimal structural design (truss topology design).}
The problem of optimally distributing material in a truss under
given loads is a classical example of a huge-scale problem: the number of
potential bars grows quadratically with the number of nodes, the stiffness
matrix is sparse, and the objective function has a maximum form. It is precisely
this problem that is used in~\cite{Nes2014huge} to demonstrate the effectiveness of
the proposed technique.

\smallskip
\emph{Sparse linear systems and PageRank.} Computing
the PageRank vector reduces to solving a sparse linear system (or a
least-squares problem) with a matrix of size
$10^7$\,--\,$10^9$. The paper by A.\,V.~Anikin and
coauthors~\cite{AnikinEtAl2022} systematically compares
methods for this class of problems---gradient, coordinate descent,
huge-scale, and specialized schemes such as the power
method---and shows which relationships among sparsity,
required accuracy, and graph structure favor each approach.

\subsection{What changed and the current state of the art}

Following~\cite{Nes2012cd}, coordinate descent methods evolved from heuristics into
a theoretically justified and practically dominant tool for
large-scale problems. P.~Richt\'arik and
M.~Tak\'a\v{c}~\cite{RichtarikTakac2014} extended the analysis to composite and
parallel formulations; S.~Shalev-Shwartz and
T.~Zhang~\cite{SSZ2013} developed stochastic dual
coordinate ascent (SDCA) for machine learning problems; accelerated
variants were developed by Lee and Sidford, Allen-Zhu and coauthors, and Nesterov
and Stich~\cite{NesStich2017}. The practical impact is even more evident: the
\texttt{glmnet} package for $\ell_1$-regularized problems and the
\texttt{LIBLINEAR} package for support vector machines implement
coordinate descent and have become standard in applied statistics. The general
idea---``make a step $n$ times cheaper, pay with $n$ times as
many steps, and gain in the constants''---also underlies all of
modern stochastic optimization, where a term of a sum is chosen at random
instead of a coordinate.

\section{Gradient-free methods (2011--2017)}
\label{sec:zo}

\subsection{Background}

Problems in which only function \emph{values} are available are
a classical topic. The direct search methods of Hooke--Jeeves (1961) and
Nelder--Mead (1965), Matyas's random search, the work of
L.\,A.~Rastrigin on random search, the
Kiefer--Wolfowitz procedure (1952)~\cite{KieferWolfowitz1952} with
finite-difference gradient approximation, and Yu.\,M.~Ermoliev's stochastic
approximation had all existed for decades. But there were
no complexity bounds displaying the dependence on dimension;
yet this dependence is the essence of the matter: intuitively,
without a gradient one must ``pay in dimension'', but how much exactly?

The answer was given in the paper by Yu.\,E.~Nesterov and
V.\,G.~Spokoiny~\cite{NesSpok2017} (a 2011 CORE preprint).

\subsection{Gaussian smoothing and a two-point estimator}

The basic construction is smoothing by Gaussian convolution:
\begin{equation}\label{eq:gsmooth}
   f_\mu(x)=\Exp_{u\sim\mathcal N(0,I_n)}\bigl[f(x+\mu u)\bigr],
   \qquad\mu>0 .
\end{equation}
We immediately note an oracle requirement without which~\eqref{eq:gsmooth}
is not well-defined: a Gaussian vector is unbounded, so $f$ is assumed
to be defined (and the oracle available) on all of $\R^n$ or in a sufficiently
large neighborhood of the working region. For a problem on a set $Q$,
one needs either an extension of the function beyond $Q$, or smoothing
with a distribution of bounded support (for example, uniform on a ball),
or a special scheme with projection.

The properties of~\eqref{eq:gsmooth} follow by direct calculation: $f_\mu$
is convex if $f$ is convex; if $f$ has an $L$-Lipschitz gradient, then
\begin{equation}\label{eq:gsm1}
   f(x)\ \le\ f_\mu(x)\ \le\ f(x)+\frac{\mu^2}{2}Ln ,
\end{equation}
since $f_\mu(x)-f(x)=\Exp[f(x+\mu u)-f(x)]\le
\Exp\bigl[\mu\dotp{\grad f(x)}{u}+\tfrac{L\mu^2}{2}\nrm{u}^2\bigr]
=\tfrac{L\mu^2 n}{2}$; if $f$ is only $M$-Lipschitz, then
$|f_\mu(x)-f(x)|\le\mu M\Exp\nrm{u}\le\mu M\sqrt n$. The key
identity is
\begin{equation}\label{eq:twopoint}
   \grad f_\mu(x)=\Exp_u\Bigl[\frac{f(x+\mu u)-f(x)}{\mu}\,u\Bigr].
\end{equation}
In other words, the random variable
$g_\mu(x,u)=\frac{f(x+\mu u)-f(x)}{\mu}u$, whose evaluation requires
\emph{two} calls to the function-value oracle, is an \emph{unbiased} estimator of
the gradient of the smoothed function. This is the central shift in
perspective: we do not approximate $\grad f$ by finite differences (which would require
$n+1$ calls and introduce an unavoidable bias), but evaluate, \emph{exactly in
expectation}, the gradient of a different function $f_\mu$ close to $f$, and
apply a stochastic gradient method to $f_\mu$. The error from
this substitution is controlled by~\eqref{eq:gsm1}.

The price is variance. The estimator~\eqref{eq:twopoint} satisfies
an inequality of the form
\begin{equation}\label{eq:zovar}
   \Exp_u\nrm{g_\mu(x,u)}_2^2\ \le\ \Oh(n)\,\nrm{\grad f(x)}_2^2
   +\Oh\bigl(\mu^2L^2n^3\bigr),
\end{equation}
that is, the second moment of the estimator is $n$ times the squared norm of the true
gradient. It is precisely the factor $n$ in~\eqref{eq:zovar} that carries over into
the final complexity bounds.

\begin{theorem}[\cite{NesSpok2017}]\label{th:zo}
\emph{Class:} $f$ is convex on $\R^n$ with an $L$-Lipschitz gradient,
and $\nrm{x_0-x^{*}}\le R$. \emph{Oracle:} a two-point oracle of
(exact) function values. \emph{Guarantee:} in expectation. Then a method using
the estimator~\eqref{eq:twopoint} with suitable $\mu$ and step size attains
$\Exp[f(x_k)]-f^{*}\le\eps$ in
\[
   \Oh\!\left(\frac{n\,LR^2}{\eps}\right)
   \quad\text{calls to the function-value oracle,}
\]
while its accelerated variant has the rate
$\Exp[f(x_k)]-f^{*}=\Oh\bigl(n^2LR^2/k^2\bigr)$, that is,
\[
   \Oh\!\left(n\sqrt{\frac{LR^2}{\eps}}\right)
   \quad\text{calls.}
\]
For convex $M$-Lipschitz (nonsmooth) $f$, the corresponding bound is
$\Oh\bigl(nM^2R^2/\eps^2\bigr)$ if the oracle returns a directional
derivative, and $\Oh\bigl(n^2M^2R^2/\eps^2\bigr)$ for differences of function values.
\end{theorem}

It is useful to write out the dependence $\Oh(n^2/k^2)$ explicitly: otherwise,
a complexity table is easily misread as $\Oh(n/k^2)$, which is incorrect.

Comparison with Theorem~\ref{th:fgm} answers the original question:
\emph{in the two-point model with exact values, the absence of the gradient
costs a factor of $n$ in the number of oracle calls}. How unavoidable this price
is constitutes a subtler question, and here one must carefully distinguish
between models.

The best-known lower bound~\cite{DuchiEtAl2015} was proved in a
\emph{stochastic} setting: the algorithm receives pairs of observations
$\bigl(F(\theta^t,X^t),F(\tau^t,X^t)\bigr)$, where $X^t$ is an independent
sample from an unknown distribution, and the minimax is taken over pairs
``loss function $+$ distribution''. In this paired-observation model,
the linear dependence on dimension is minimax optimal. The model of
Theorem~\ref{th:zo} is different: values of a single fixed function
are returned \emph{exactly}, and the bound from~\cite{DuchiEtAl2015} cannot
be transferred to it directly.

For an exact (noiseless) function-value oracle, the question long remained
open; the recent preprint~\cite{ZhangEtAl2026} obtained nearly
optimal lower bounds for nonsmooth Lipschitz problems. These cannot be carried over
to the smooth class of Theorem~\ref{th:zo}: both the function class
and the accuracy regime differ. For the smooth class, the answer is given in~\cite{DornEtAl2026}:
any method, including a randomized one, that uses only
exact function values requires, on the class of convex functions with an $L$-Lipschitz gradient (in
the Euclidean norm) and a minimizer in a ball of radius $R$,
\begin{equation}\label{eq:zolb}
   N(\eps)=\Om\Bigl(n\min\Bigl\{n,\ \sqrt{LR^2/\eps}\Bigr\}\Bigr)
\end{equation}
oracle calls to find an $\eps$-solution with probability~$2/3$. For
$LR^2/n^2\le\eps\le LR^2/32$, this matches the accelerated-method bound in
Theorem~\ref{th:zo}: in the smooth Euclidean case, the factor $n$
\emph{cannot be eliminated}
(for deterministic methods, a closely related result was independently obtained
in~\cite{WuEtAl2026}). At higher accuracy, the ``second
scale'' of \S\,\ref{sec:cutting} comes into play: for example, a quadratic function can
be reconstructed from $\Oh(n^2)$ values.

The dependence on the prox-structure is particularly instructive
(cf.\ \S\,\ref{sec:primaldualgrad}). Suppose that smoothness is still
measured in the Euclidean norm, while a minimizer lies in an $\ell_p$-ball
of radius $R_p$, $1\le p\le2$, $\frac1p+\frac1q=1$. Then, in the corresponding
accuracy range~\cite{DornEtAl2026},
\begin{equation}\label{eq:zolbp}
   N(\eps)=\widetilde\Theta\Bigl(n^{\frac12+\frac1q}
   \sqrt{LR_p^2/\eps}\Bigr).
\end{equation}
In the Euclidean case ($p=q=2$), the factor is $n$, whereas for an $\ell_1$-ball
($q=\infty$) it is only $\sqrt n$. Crucially,
the bound~\eqref{eq:zolbp} \emph{is attained}, precisely with
a non-Euclidean prox-structure: by an accelerated gradient-free method in
which the prox-function
is $1$-strongly convex with respect to $\nrm{\cdot}_p$, while the random directions
remain uniformly distributed on the Euclidean
sphere~\cite{GorbunovDvurGasn2022}. The source of the gain is concentration of measure:
the second moment of the estimator~\eqref{eq:twopoint}, measured in the dual norm
$\nrm{\cdot}_\infty$, has (for small $\mu$) order
$\ln n\cdot\nrm{\grad f(x)}_2^2$ instead of $n\nrm{\grad f(x)}_2^2$
in~\eqref{eq:zovar}, since the largest coordinate of a Gaussian
vector (of order $\sqrt{2\ln n}$) is much smaller than its Euclidean norm
(of order $\sqrt n$). The nonsmooth case is similar:
$\widetilde\Theta\bigl(n^{2/q}M^2R_p^2/\eps^2\bigr)$, and in
$\ell_1$-geometry the dimension remains only under
a logarithm~\cite{Shamir2017,DornEtAl2026}. In contrast to first-order
methods (\S\,\ref{sec:primaldualgrad}), the mismatch between norms
here does not prevent the prox-method from being optimal: the dimensional
factor plays the principal role.

Finally, the factor generally
depends on the details of the model: one-point and two-point oracles,
deterministic and stochastic values, shared or independent noise in
a pair of queries, the availability of a parallel batch, and the geometry of the set
give different dimension dependences. The phrase ``exactly a factor of $n$''
is meaningful only after all these details have been fixed.

\subsection{The full generalization: smoothing as a reduction to a smooth problem}
\label{sec:zoreduction}

The reasoning of~\eqref{eq:gsmooth}\,--\,\eqref{eq:twopoint} admits a far-reaching
generalization obtained in~\cite{GasnEtAlICML2022}: \emph{any}
accelerated method for smooth stochastic problems becomes a
gradient-free method for \emph{nonsmooth} problems. Let $f$ be convex and
$M$-Lipschitz. Smooth it by averaging over a Euclidean ball of radius $\mu$:
\begin{equation}\label{eq:ballsmooth}
   \hat f_\mu(x)=\Exp_u f(x+\mu u),
\end{equation}
where the vector $u$ is uniformly distributed in the unit Euclidean ball.
Then $f(x)\le\hat f_\mu(x)\le f(x)+\mu M$ for all $x$, and the gradient of $\hat f_\mu$ is Lipschitz with
constant $L_\mu=\Oh(\sqrt n\,M/\mu)$. An unbiased estimator of
$\grad\hat f_\mu$ from two values is given by the symmetric difference along
a direction $e$ uniformly distributed on the unit sphere:
\begin{equation}\label{eq:spheregrad}
   \hat g_\mu(x,e)=n\,\frac{f(x+\mu e)-f(x-\mu e)}{2\mu}\,e,\qquad
   \Exp_e\nrm{\hat g_\mu(x,e)}_2^2\le c\,nM^2,
\end{equation}
where $c$ is an absolute constant. Lipschitz continuity of $f$ alone would give
a second-moment bound of order $n^2M^2$ (hence the $n^2$ in
Theorem~\ref{th:zo} for differences of values); the improvement to $nM^2$
follows from concentration of measure on the sphere, and it is precisely this that yields
the bound $\Oh(nM^2R^2/\eps^2)$~\cite{Shamir2017}. Thus, with $\mu=\eps/(2M)$,
the nonsmooth problem has been replaced by a smooth one with constant
$L_\mu\simeq\sqrt n\,M^2/\eps$ and a stochastic gradient with variance
$\sigma^2\simeq nM^2$. It remains to apply an accelerated method with
batching~\eqref{eq:batch}\,--\,\eqref{eq:batch2}: it makes $N$
iterations and $T$ oracle calls, where
\[
   N=\Oh\Bigl(\sqrt{\frac{L_\mu R^2}{\eps}}\Bigr)
   =\Oh\Bigl(\frac{n^{1/4}MR}{\eps}\Bigr),\qquad
   T=\Oh\Bigl(\frac{nM^2R^2}{\eps^2}\Bigr).
\]
The total number of calls $T$ cannot be improved (up to a logarithmic factor;
see the discussion of lower bounds above), while the number of sequential iterations
$N$, when a batch is evaluated in parallel, is the same as for
\emph{first}-order methods using the same smoothing; whether it
is optimal is not yet known. The maximum allowable
noise level in function values for which the bounds
remain valid is also established. Thus the result~\cite{NesSpok2017} becomes
a prototype of the general scheme ``smoothing $+$ acceleration $+$ batching'', with
the inexact oracle of Definition~\ref{def:inexact} acting as the connecting
link.

The scheme readily extends to other formulations: stochastic,
saddle-point, and distributed problems~\cite{GasnEtAlICML2022}, and,
in~\cite{AbroninEtAl2026}, stochastic problems with functional
inequality constraints whose values are also accessible only through
an oracle. The number of sequential iterations remains
$\Oh(n^{1/4}MD/\eps)$, and the number of oracle calls remains
$\Oh(nM^2D^2/\eps^2)$; this bound cannot be improved up to
a logarithmic factor (here $D$ is the diameter of the feasible set, and the constant $M$
accounts for both the constraints and the size of the dual multipliers).

\subsection{What changed and the current state of the art}

Before~\cite{NesSpok2017}, gradient-free optimization was a collection of
recipes compared by experimental counts of function evaluations.
Afterwards, it became a branch of complexity theory with lower bounds, optimal
methods, and a clear ``price of not knowing the gradient''. The practical
consequences proved more far-reaching than one might have expected.

\practice \emph{Reinforcement learning.} The ``evolution strategies''
method proposed in 2017 as a scalable alternative to
reinforcement learning~\cite{Salimans2017} is literally
scheme~\eqref{eq:twopoint} with Gaussian smoothing: policy parameters
are perturbed by Gaussian noise, and the gradient is estimated from differences
in rewards.

\practice \emph{Attacks on neural networks.} Constructing adversarial
examples in a ``black-box'' setting, in which only a model's responses
and not its gradients are available, relies on the same two-point estimators.

\practice \emph{Fine-tuning large language models.} The MeZO
(memory-efficient zeroth-order optimization) method~\cite{Malladi2023} fine-tunes
language models using forward passes only, allowing memory
usage at the inference rather than training level (cf.\ the discussion of memory in
\S\,\ref{sec:ad}).

\hypoth Direct application of the worst-case bound in Theorem~\ref{th:zo}
predicts an unfavorable linear dependence on the full dimension for MeZO:
a factor of $n$ with $n\sim10^{9}$. One possible
explanation for its observed effectiveness is that the actual
dependence is governed not by the full dimension but by the \emph{effective}
dimension (the intrinsic dimension of the problem or the local rank of the Hessian). We emphasize that this
is a hypothesis, not a consequence of Theorem~\ref{th:zo}; refining bounds in
terms of effective dimension is an active topic today.

\smallskip
\emph{Heavy tails and an inexact oracle.} A separate line of research studies
gradient-free methods when function values are noisy:
with heavy-tailed random, bounded, or adversarial noise. Bounds
have been obtained that show how the allowable noise level is related to
attainable accuracy; for a systematic review of the current state of the art,
see~\cite{GasnZO2022}.

\smallskip
\emph{Adaptivity.} The ``smoothing $+$ acceleration $+$ batching'' scheme
requires knowledge of the problem constants: the smoothing radius, step size, and
batch size depend on them. \emph{Adaptive} variants have recently been developed.
The paper~\cite{DornOsmanovEtAl2026} proposes an accelerated two-point method for
stochastic convex problems on a compact set of known diameter~$D$; it selects all
three quantities from observed function values. Without knowing
$L$, $M$, or the second moment of stochastic gradients, it
recovers, up to logarithmic factors, the bounds of the nonadaptive
scheme in \S\,\ref{sec:zoreduction} (with $R$ replaced by $D$), in a single
run: the part of its bound that depends on the number $N$ of sequential iterations
is determined by the better of two rates---the ``smooth'' rate $LD^2/N^2$ and
the ``nonsmooth'' rate $n^{1/4}MD/N$ (for nonsmooth $f$, we set $L=+\infty$).
In essence, this is a gradient-free counterpart of the universal methods to be
discussed in \S\,\ref{sec:univ}, albeit so far only for the two extreme
cases of the smoothness scale.

\subsection{When the function-value oracle itself is inexact}

Let us return to the question posed in \S\,\ref{sec:adnograd}. If
the gradient is unavailable, then the \emph{exact}
function value is generally unavailable as well: the value comes from measurement,
simulation, or the solution of an auxiliary problem. Crucially, in
the last case the oracle accuracy is \emph{controllable}: one can spend
more time to obtain a more accurate value.

This leads to a formulation that adds a third parameter
to the usual ones. Suppose that the oracle returns the value $f(x)$ to accuracy $\delta$ in
time $c(\delta)$, with $c(\delta)\propto\delta^{-\gamma}$ for
some $\gamma>0$. One should then minimize not the number of calls
$N(\eps,\delta)$, but the total time
\[
   \mathcal T(\eps)=\min_{\delta>0}\ N(\eps,\delta)\cdot c(\delta),
\]
and the optimal trade-off between oracle accuracy and the number of calls
depends on $\gamma$, the dimension, and the function class. This
formulation is the subject of a separate line of research; see~\cite{SuvorikovaEtAl2026}.

Such a formulation arises outside the gradient-free setting as well. Suppose
that the \emph{gradient} is evaluated with controllable accuracy: an approximation with
accuracy $\delta$ costs on the order of $(c_0/\delta)^\gamma$. For smooth convex
problems, this model is studied in~\cite{PaingJacot2026}, and the answer seems
paradoxical: for $\gamma>2$, the optimal cost of solving the problem to
accuracy $\eps$ is $\Theta\bigl((c_0R/\eps)^\gamma\bigr)$, that is,
the same order as \emph{one} gradient evaluation
with accuracy $\eps/R$ (hence the paper's title: convex optimization
is ``free'' when accuracy is expensive). A naive scheme (gradient descent with
constant accuracy $\delta\simeq\eps/R$) would require
$\Oh(\eps^{-(\gamma+1)})$.
The gain comes from multilevel randomization: an unbiased estimator of the most
accurate gradient is assembled from corrections between approximations at neighboring
accuracy levels, with coarse levels queried frequently and accurate ones
rarely. Acceleration is beneficial in this model only for $\gamma<2$, which
agrees with the linear accumulation of error in accelerated methods
(\S\,\ref{sec:accbatch}). Finally, in \S\,\ref{sec:superfast} we shall see
that there, too, one must balance the accuracy of the inner solution against
the cost of the outer iteration.

%% file: en/r10_univ.tex

\section{Universal gradient methods (2013--2015)}
\label{sec:univ}

\subsection{Background: the problem of choosing a class}

By the 2010s, the complexity theory of convex optimization presented
a coherent but \emph{stratified} picture (Table~\ref{tab:map}): nonsmooth
problems---$\Oh(\eps^{-2})$, smooth problems---$\Oh(\eps^{-1/2})$, \emph{smooth}
strongly convex problems---logarithmic dependence on accuracy. An optimal method was known for each class; but
using it required \emph{knowing} which class
the problem belonged to and the value of the corresponding constant. In
practice, neither was known. Worse still, real functions often
do not fit entirely into any one class: for example,
$f(x)=\nrm{Ax-b}_p^{\,p}$ for
$1<p<2$ has an unbounded gradient Lipschitz constant, yet it seems a pity to
call it nonsmooth.

A natural scale of intermediate classes is given by \emph{H\"older continuity}
of the gradient (the letter $\nu$ here denotes the H\"older exponent and has
no connection with the barrier parameter in \S\,\ref{sec:ipm}):
\begin{equation}\label{eq:holder}
   \nrm{\grad f(x)-\grad f(y)}_{*}\ \le\ L_\nu\nrm{x-y}^{\nu},
   \qquad \nu\in[0,1],
\end{equation}
whence, by integration,
\begin{equation}\label{eq:holderub}
   f(y)\ \le\ f(x)+\dotp{\grad f(x)}{y-x}
   +\frac{L_\nu}{1+\nu}\nrm{y-x}^{1+\nu} .
\end{equation}
For $\nu=1$, this is ordinary smoothness; the case $\nu=0$ requires a separate
discussion (see \S\,\ref{sec:nu0} below).

It is important that this scale did not originate in 2015. Already in the
1985 paper by A.\,S.~Nemirovski and Yu.\,E.~Nesterov~\cite{NemNes1985},
intermediate (H\"older) smoothness was considered precisely from
the perspective of acceleration: optimal methods were obtained for the entire
family of classes~\eqref{eq:holder}, and the \emph{composite}
structure of the problem was also examined---the case $f(x)=\max\{f_1(x),\dots,f_m(x)\}$ and, more
generally, a convex monotone composition of a scalar outer function with
a smooth vector mapping. (The pointwise \emph{minimum} of convex
functions is generally nonconvex; it is specifically the maximum that is meant.)
This line of research subsequently developed into the notion of an \emph{inexact $(\delta,L)$-oracle} in~\cite{DevolderNesterovGlineur2014} (which we discuss in the next section), and was further extended to an \emph{inexact model} of a function (a model in place of the
function itself)~\cite{StonyakinEtAl2021}, in which
composite optimization, max-type problems, relative smoothness, and
an inexact oracle are special cases of a single construction.
The optimal methods of~\cite{NemNes1985}, however, required knowledge of
$\nu$ and $L_\nu$.
They solved the problem optimally for each fixed $\nu$, whereas universal methods solve it with the best $\nu$ over the entire
interval $\nu \in [0;1]$.

\subsection{The inexact $(\delta,L)$-oracle and the price of acceleration}
\label{sec:inexactoracle}

Before the universal method was developed, O.~Devolder,
F.~Glineur, and Yu.\,E.~Nesterov proposed a key tool in~\cite{DevolderNesterovGlineur2014}: a new model of inexact function and gradient evaluation, the $(\delta,L)$-oracle.
At first the result might have seemed purely technical, but it had an enormous influence on convex optimization as a whole;
in particular, it became a key tool and idea in the development of universal methods.

The key observation is that inexactness can be measured not separately for each object (say, an inexact gradient relative to the exact one, or an inexact function value relative to the exact one), but directly at the level of the linear model of the function. This yields the following definition.
\begin{definition}[$(\delta,L)$-oracle]\label{def:dLoracle}
Let $f$ be a convex function on a convex set $Q$. It is said
to be equipped with a \emph{$(\delta,L)$-oracle} if, for every point
$y\in Q$, one can compute a pair $\bigl(f_{\delta,L}(y),g_{\delta,L}(y)\bigr)$
such that, for all $x\in Q$,
\begin{equation}\label{eq:dLoracle}
   0\ \le\ f(x)-\bigl[f_{\delta,L}(y)+\dotp{g_{\delta,L}(y)}{x-y}\bigr]
   \ \le\ \frac{L}{2}\nrm{x-y}^2+\delta .
\end{equation}
\end{definition}

The left inequality in~\eqref{eq:dLoracle} says that the linearization returned by the oracle
remains a \emph{global lower model}---exactly what
convexity gives in the exact case. The right inequality says that it falls below the function by
at most a quadratic term plus $\delta$. For $\delta=0$,
condition~\eqref{eq:dLoracle} is exactly the definition of a convex function
with an $L$-Lipschitz gradient.

We now state the main result.

\begin{theorem}[\cite{DevolderNesterovGlineur2014}]\label{th:dng}
\emph{Class:} $f$ is convex on a convex compact set $Q$. \emph{Oracle:} a $(\delta,L)$-oracle. Then
\begin{enumerate}[label=\rm(\alph*),leftmargin=2.2em,itemsep=4pt]
\item \textbf{Ordinary gradient method.}
      \[
         f(x_N)-f^{*}\ \le\ \frac{LR^2}{2N}+\delta .
      \]
      Oracle error \emph{does not accumulate}: it enters the bound
      additively, with coefficient~1, for any number of iterations.
\item \textbf{Fast gradient method.}
      \[
         f(x_N)-f^{*}\ \le\ \Oh\Bigl(\frac{LR^2}{N^2}\Bigr)+\Oh(N\delta).
      \]
      The error accumulates \emph{linearly} with the number of iterations, and a lower bound shows that this linear accumulation
cannot be improved for an accelerated method.
\end{enumerate}
\end{theorem}

The conclusion is simple: acceleration is not free. An accelerated method outperforms
the ordinary one in rate ($N^{-2}$ versus $N^{-1}$), but is less tolerant of oracle
errors ($N\delta$ versus $\delta$). Between these extremes are intermediate methods with
bounds of the form $\Oh\bigl(LR^2/N^{1+p}+N^{p}\delta\bigr)$, $p\in[0,1]$,
where $p=0$ corresponds to the ordinary gradient method and $p=1$ to the fast one~\cite{DevolderGlineurNes2013}, and their universal version~\cite{KamzolovDvurGasn2019}.
Thus the choice is not simply ``yes, accelerate'' or ``no'': the degree of acceleration
can be adjusted to the inexactness of the oracle.

Here the seemingly technical paper becomes a change of paradigm.
If $\delta$ is viewed as interference, an error in evaluating the gradient or the function, then Theorem~\ref{th:dng} sounds like
a warning that acceleration may be unstable. But $\delta$
can instead be viewed as a parameter that we choose ourselves.
The inexact oracle then ceases to be a problem and becomes a tool
for extending the class of problems to which the machinery of smooth
optimization applies. This is precisely the observation used below:

\begin{itemize}[leftmargin=1.6em,itemsep=3pt]
\item \emph{Universality} (\S\,\ref{sec:idea}). Any function with
      a H\"older-continuous gradient is equipped with a $(\delta,L(\delta))$-oracle for
      \emph{every} $\delta>0$---with the exact value and the exact
      (sub)gradient, but an inflated $L$. Thus a nonsmooth problem
      \emph{is} a smooth problem with an inexact oracle, and no
      new theory is needed: it suffices to choose
      $\delta$ correctly. Theorem~\ref{th:dng} immediately tells us how:
      $\delta\sim\eps$ for a nonaccelerated scheme and
      $\delta\sim\eps/N$ for an accelerated one, so that the accumulated
      $N\delta$ does not exceed $\eps/2$.
\item \emph{Batching} (\S\,\ref{sec:accbatch}). A mini-batch of
      a stochastic oracle behaves like an inexact oracle with
      $\delta$ decreasing as the batch size grows. Part~(b) of
      Theorem~\ref{th:dng} explains a fact observed in practice:
      accelerated schemes need substantially larger batches than
      nonaccelerated ones, and this is not a heuristic but a consequence of linear
      accumulation.
\item \emph{An inexact model of a function}~\cite{StonyakinEtAl2021}. Replacing
      the quadratic term in~\eqref{eq:dLoracle} by an arbitrary
      Bregman divergence yields a construction whose special cases
      include composite optimization, max-type problems,
      relative smoothness~\cite{LuFreundNes2018,BBT2017}, and variational inequalities.
\end{itemize}

Thus, at first sight, the paper~\cite{DevolderNesterovGlineur2014}
answered a narrow question: what happens if the gradient is evaluated
approximately? In essence, however, it provided a language for constructing an entire
family of methods that work on a broad class of problems without knowing
which class a problem belongs to. The universal gradient
method, discussed in the next subsection, is the first and most
striking example of such an application.

\subsection{The idea: inexact smoothness}
\label{sec:idea}

The paper~\cite{Nes2015univ} (a 2013 CORE preprint) rests on
the following simple but highly fruitful observation: \emph{any
H\"older bound~\eqref{eq:holderub} implies a quadratic bound with
an arbitrarily small additive term}.

\begin{lemma}\label{lem:univ}
Suppose that~\eqref{eq:holderub} holds for some $\nu\in[0,1]$ and
$L_\nu>0$. Then, for every $\delta>0$ and all $x,y$,
\begin{equation}\label{eq:inexact}
   f(y)\ \le\ f(x)+\dotp{\grad f(x)}{y-x}
   +\frac{L(\delta)}{2}\nrm{y-x}^2+\delta,
\end{equation}
where
\begin{equation}\label{eq:Ldelta}
   L(\delta)=L_\nu^{\frac{2}{1+\nu}}
   \left[\frac{1-\nu}{2(1+\nu)\,\delta}\right]^{\frac{1-\nu}{1+\nu}}.
\end{equation}
\end{lemma}

\begin{proof}
It suffices to verify that $a\,r^{1+\nu}\le\frac{b}{2}r^2+\delta$ for all
$r\ge0$ with $a=\frac{L_\nu}{1+\nu}$ and $b=L(\delta)$.
First let $\nu<1$. The function
$g(r)=\frac{b}{2}r^2-a r^{1+\nu}$ attains its minimum at
$r_{*}=\bigl(a(1+\nu)/b\bigr)^{1/(1-\nu)}=(L_\nu/b)^{1/(1-\nu)}$, and
\[
   g(r_{*})=r_{*}^{1+\nu}\Bigl(\frac{b}{2}r_{*}^{1-\nu}-a\Bigr)
   =r_{*}^{1+\nu}L_\nu\Bigl(\frac12-\frac{1}{1+\nu}\Bigr)
   =-\,\frac{L_\nu(1-\nu)}{2(1+\nu)}
     \Bigl(\frac{L_\nu}{b}\Bigr)^{\frac{1+\nu}{1-\nu}} .
\]
The condition $-g(r_{*})\le\delta$ is equivalent to
$b\ge L_\nu^{2/(1+\nu)}
\bigl[\frac{1-\nu}{2(1+\nu)\delta}\bigr]^{(1-\nu)/(1+\nu)}$, which
gives~\eqref{eq:Ldelta}. For $\nu=1$, the exponent vanishes and
$L(\delta)=L_1$ for any $\delta$; for $\nu=0$, we obtain
$L(\delta)=L_0^2/(2\delta)$---both extremes agree with direct
verification.
\end{proof}

The meaning of Lemma~\ref{lem:univ} is this: from the viewpoint of a gradient method,
\emph{any} function in the class~\eqref{eq:holder} behaves like a smooth one,
provided we accept an additive error $\delta$.
In the language of \S\,\ref{sec:inexactoracle}, inequality~\eqref{eq:inexact}
together with convexity means that $f$ is equipped with
a $(\delta,L(\delta))$-oracle (Definition~\ref{def:dLoracle}) for every
$\delta>0$---an ``honest'' oracle returning exact $f(x)$ and
$\grad f(x)$; only the \emph{model} is inexact here, not the information.
How to choose $\delta$ is already known from Theorem~\ref{th:dng}:
$\delta\sim\eps$ for a nonaccelerated scheme (the error does not accumulate) and
$\delta\sim\eps/N$ for an accelerated one (it accumulates linearly). One step
remains: choose $L$ \emph{adaptively}, without knowing either $\nu$ or $L_\nu$.
This is exactly how the
universal method works: at each iteration, backtracking
adjusts the current estimate $L_k$ until
inequality~\eqref{eq:inexact} holds at a pair of trial points.
We emphasize that backtracking searches not for the ``true'' Lipschitz constant, but for a suitable $L$ in
inequality~\eqref{eq:inexact}, that is, for the \emph{locally best} $\nu$ and $L_\nu$ at the current point. This explains
the infimum over $\nu$ in Theorem~\ref{th:univ}.

\subsection{The extreme case $\nu=0$}
\label{sec:nu0}

Condition~\eqref{eq:holder} is written for the gradient and formally
assumes differentiability, whereas for $\nu=0$ we want
to cover arbitrary Lipschitz convex functions, which may be
nondifferentiable. This is resolved as follows: the gradient is replaced by
an arbitrary subgradient $g(x)\in\partial f(x)$, and one directly requires
inequality~\eqref{eq:inexact} with $\grad f(x)$ replaced
by $g(x)$. The extreme case then follows by direct calculation.

\begin{proposition}\label{pr:nu0}
Let $f$ be convex and \emph{globally} $M$-Lipschitz on $\R^n$. Then, for
every $x$, every $g(x)\in\partial f(x)$, and every $\delta>0$,
\[
   f(y)\ \le\ f(x)+\dotp{g(x)}{y-x}
   +\frac{L(\delta)}{2}\nrm{y-x}^2+\delta,
   \qquad L(\delta)=\frac{L_0^2}{2\delta},\quad L_0=2M .
\]
\end{proposition}

\begin{proof}
By convexity, $f(x)\ge f(y)+\dotp{g(y)}{x-y}$, whence
$f(y)-f(x)-\dotp{g(x)}{y-x}\le\dotp{g(y)-g(x)}{y-x}\le2M\nrm{y-x}$.
It remains to apply Young's inequality: $2Mr\le\frac b2r^2+\delta$ for all
$r\ge0$ is equivalent to $b\ge2M^2/\delta=L_0^2/(2\delta)$.
\end{proof}

Thus, for $\nu=0$, the quantity $L_0=2M$ is
an upper bound on the diameter
of the set of subgradients, and the resulting $L(\delta)$ agrees exactly with
formula~\eqref{eq:Ldelta} at $\nu=0$. The scale~\eqref{eq:holder}
is closed at both ends.

\subsection{Convergence rate bounds}

\begin{theorem}[\cite{Nes2015univ}]\label{th:univ}
\emph{Class:} $f$ is convex on a convex compact set $Q$ and satisfies
inequality~\eqref{eq:inexact} (in particular, any
condition~\eqref{eq:holder} for some $\nu\in[0,1]$).
\emph{Oracle:}
value and (sub)gradient. The \emph{method} uses no information about
$\nu$ or $L_\nu$, but uses the required accuracy $\eps$ and
an initial estimate $L_{-1}>0$. \emph{Guarantee:} deterministic. The infima below
are taken over those $\nu$ for which the constant $L_\nu$
in~\eqref{eq:holder} is finite (for all other values, we set
$L_\nu=+\infty$). Then, to attain accuracy $\eps$, it suffices to perform
\[
   N=\Oh\!\left(\inf_{\nu\in[0,1]}\left[
   \left(\frac{L_\nu R^{1+\nu}}{\eps}\right)^{\!\frac{2}{1+\nu}}\right]\right)
   \quad\text{iterations (nonaccelerated variant)},
\]
\[
   N=\Oh\!\left(\inf_{\nu\in[0,1]} \left[
   \left(\frac{L_\nu R^{1+\nu}}{\eps}\right)^{\!\frac{2}{1+3\nu}}\right]\right)
   \quad\text{iterations (accelerated variant)}.
\]
In particular, $\nu=1$ gives $\Oh(LR^2/\eps)$ and
$\Oh\bigl(\sqrt{LR^2/\eps}\bigr)$, while $\nu=0$ gives
$\Oh(M^2R^2/\eps^2)$ in both cases.
\end{theorem}

We emphasize four points.

First, there is \emph{one} method, but as many bounds as there are classes: it
automatically attains the best rate over all $\nu$ for which
the problem belongs to the class~\eqref{eq:holder}---hence the name
``universal''.
Second, the bound for the \emph{accelerated} method matches the known lower
bounds for every $\nu$, so universality comes almost
for free.
Third, as to the cost of backtracking, the standard
amortized analysis of line search gives a \emph{constant average}
number of additional oracle calls per iteration (at most two), plus
an \emph{additive} logarithmic term for the initial adjustment of
the estimate of $L$; no multiplicative logarithmic factor appears
in the bound. Fourth, there is a real price:
the methods require the target accuracy $\eps$ to be known in advance. As
a consequence, the methods are not ``anytime'': if one stops early or
wants to refine the answer, the computation must be restarted with a new $\eps$.
This limitation was removed only in 2024 by Yu.\,E.~Nesterov's student Anton Rodomanov and coauthors~\cite{RodomanovEtAl2024}.

\begin{example}\label{ex:univ}
The function $f(x)=\nrm{Ax-b}_2^{3/2}$ satisfies~\eqref{eq:holder} with
$\nu=1/2$, but not with $\nu=1$ (the constant is infinite), nor
reasonably with $\nu=0$ on an unbounded set. For this function, the universal accelerated
method gives $N=\Oh\bigl((L_{1/2}R^{3/2}/\eps)^{4/5}\bigr)$---
a rate intermediate between $\eps^{-1/2}$ and $\eps^{-2}$ that, without
Theorem~\ref{th:univ}, would have to be worked out separately.
\end{example}

\begin{remark}[on the transport problem]\label{rem:univtraffic}
It is often noted that the universal method behaves particularly well on the dual transport problem~\eqref{eq:sddual}. This needs to be stated carefully. The function $\sum_wd_wT_w(t)$ is piecewise linear, hence belongs to the class $\nu=0$, and \emph{no} theorem gives it a better worst-case bound than $\Oh(\eps^{-2})$. The observed acceleration is an empirical effect of local backtracking: over most of the trajectory the same combination of shortest paths is active, the function behaves locally as a smooth one, and the adaptive procedure chooses large steps. Thus this is a \emph{practical interpretation} of the method's behavior. Corresponding experiments can be found in~\cite{GasnTranspMono2020}.
\end{remark}

\subsection{Current state of the art}

The universal approach quickly became standard. In~\cite{GasnNes2018},
a universal similar triangles method was developed for stochastic
composite problems,
becoming a working tool for numerical solution of transport and
entropy-regularized problems; an accessible account of both
the universal method itself and its applications is given in
the textbook~\cite{GasnBook}.
This method combines backtracking with mini-batching and
therefore requires knowledge of the oracle variance to select the batch size;
it does not yet adapt to the noise level.

A key step was taken in~\cite{RodomanovEtAl2024}, which proposes
a line-search-free version of the universal method: instead of backtracking,
the coefficient $H_k$ is selected from an explicit ``balance equation'' between two
terms of the bound, making it computable by a single formula and, crucially,
introducing no statistical dependence between $H_k$ and $x_k$. Consequently,
the construction extends to a stochastic oracle. The resulting methods
adapt simultaneously to $\nu$, $L_\nu$, and the variance $\sigma$,
\emph{without prior
knowledge of $\eps$}: for the composite objective $F(x)=f(x)+\Psi(x)$
(see \S\,\ref{sec:composite}) and $F^{*}=\min_{x\in Q}F(x)$, the bounds
have the anytime form
\begin{equation*}
    \begin{gathered}
   \Exp[F(\bar x_k)]-F^{*}\le\inf_{\nu\in[0,1]}
   \left[\Oh\left(\frac{L_\nu D^{1+\nu}}{k^{(1+\nu)/2}}\right)\right]+\Oh\left(\frac{\sigma D}{\sqrt k}\right)\,\\
   \Exp[F(x_k)]-F^{*}\le\inf_{\nu\in[0,1]}
   \left[\Oh\left(\frac{L_\nu D^{1+\nu}}{k^{(1+3\nu)/2}}\right)\right]+\Oh\left(\frac{\sigma D}{\sqrt k}\right)
    \end{gathered}
\end{equation*}
for the nonaccelerated and accelerated schemes, respectively. The price is knowledge of
the diameter $D$ of the feasible set (instead of $R=\nrm{x_0-x^{*}}$) and
the boundedness of that set.

The idea of adaptively choosing a local
Lipschitz constant was developed independently, and in other forms, in machine learning:
AdaGrad-type adaptive steps~\cite{DuchiHazanSinger2011}, schemes without
knowledge of $L$, and methods with the Polyak step size.
The connection is more than conceptual: the step in~\cite{RodomanovEtAl2024}
is bounded above by the classical AdaGrad expression
$H_k'=D^{-1}\bigl(\sum_i\nrm{g_i-g_{i-1}}_{*}^2\bigr)^{1/2}$; thus
the universal method takes steps no shorter than AdaGrad, but with guarantees over
the entire scale~\eqref{eq:holder}. This connection was developed into a separate
result in~\cite{RodomanovJiangStich2024}: AdaGrad itself is shown to
be universal in the same
sense when expressed in terms of a $(\delta,L)$-oracle. Specifically,
gradient descent with classical AdaGrad steps $H_k'$
automatically adapts simultaneously to H\"older smoothness and to
the bias and variance of a stochastic oracle and, under additional
assumptions, provides implicit variance reduction
and acceleration---with the only required
parameter, as in~\cite{RodomanovEtAl2024}, being an estimate of the diameter
of the feasible set, $D$. In other words, the $(\delta,L)$-model idea does
not merely explain retrospectively why AdaGrad works on
H\"older-smooth problems; it also indicates which modifications of the step size extend the same
guarantee to accelerated and variance-reduced
schemes.

A separate
line of research extends universality to higher-order methods:
regularized Newton methods for functions with H\"older-continuous
Hessians~\cite{GrapigliaNes2017} and tensor methods~\cite{GrapigliaNes2020}
(\S\,\ref{sec:tensor}).

Lemma~\ref{lem:univ} also applies when smoothness differs across groups
of variables. The paper~\cite{TitovEtAl2021} considers saddle-point problems
$\min_{x\in Q_x}\max_{y\in Q_y}f(x,y)$ (here $Q_x$ and $Q_y$ are
convex compact sets) in which $f$ is $\mu_x$-strongly
convex in $x$ and $\mu_y$-strongly concave in $y$. The gradient $\grad_xf$ is assumed
to be H\"older continuous with exponent $\nu$ (with constants $L_{xx}$ in $x$ and
$L_{xy}$ in $y$), while $\grad_yf$ is Lipschitz continuous in $y$ with constant
$L_{yy}$. Then the function $\max_{y\in Q_y}f(x,y)$ has a H\"older-continuous gradient with
exponent $\nu/(2-\nu)$ and, accounting for inexact solution of the inner
problem, is equipped with an inexact oracle in the sense of \S\,\ref{sec:inexactoracle}.
An accelerated method is applied to it, while the inner problem in $y$ is solved
by a linearly convergent method, so that the scheme's parameters are simultaneously
adjusted to the constants of both variable groups: the total number of iterations
(including inner iterations) is
$\Oh\bigl(\sqrt{L(\eps)/\mu_x}\,\sqrt{L_{yy}/\mu_y}\bigr)$ up to
logarithmic factors, where $L(\eps)$ is constructed from
formula~\eqref{eq:Ldelta} for exponent $\nu/(2-\nu)$. The same paper also
considers, for the strongly monotone case, a version of the universal
Mirror Prox method with restarts,
adapting to an unknown level of smoothness. In~\cite{NguyenGasnikov2026},
the same reduction is combined with sliding (\S\,\ref{sec:sliding}): for
a problem of the form~\eqref{eq:spp} on bounded sets of diameter $D$, where
the gradients of $f$ and $g$ are H\"older continuous with exponents $\nu_x,\nu_y$ and
constants $L^f_{\nu_x},L^g_{\nu_y}$, each of the three oracles
is called at its own frequency. The gradient of $f$ is evaluated
$\Oh\bigl((L^f_{\nu_x}D^{1+\nu_x}/\eps)^{2/(1+3\nu_x)}\bigr)$ times---exactly
as in the accelerated variant of Theorem~\ref{th:univ}, and similarly for $g$, while
$\Oh(\nrm{B}D^2/\eps)$ multiplications by $B$ and $B^{\!\top}$ are required.
The H\"older exponents and constants are assumed known here.

Finally, the modern generalized smoothness condition
$(L_0,L_1)$ from \S\,\ref{sec:l0l1} is a direct relative of
the class~\eqref{eq:holder}: in both cases, the gradient's Lipschitz constant
is not constant but depends on the point, and the method must be able to
detect this.
The difference is in \emph{how} it varies: in~\eqref{eq:holder},
it degenerates as points approach each other (the scale is set by distance), while in
the $(L_0,L_1)$-model it grows with the gradient norm (the scale is set by
the function itself). Universal backtracking schemes work in both cases,
indirectly confirming that the method adapts not to a particular class but
to the locally valid inequality~\eqref{eq:inexact}.

This makes the contrast with the strong convexity constant~$\mu$ all the sharper: one can adapt
to it as well, but a fundamental difficulty arises. The analysis needs $L$
only through an inequality of the form~\eqref{eq:inexact} at points already
visited, and this can be checked: any
$H\ge L$ satisfies it, and backtracking finds a suitable value. By contrast, $\mu$ is
a \emph{global} lower bound on curvature, and no finite collection
of values and gradients guarantees that $\mu$ is at least a given
positive number: the function can coincide with a highly curved
quadratic at all visited points and be much flatter
far from them. Thus an overestimate of $\mu$ can be detected and
reduced, whereas any increase is only a guess. The preprint by
I.~Kuruzov and coauthors~\cite{KuruzovEtAl2026} shows that this difficulty
is real. A method which, with all other checks in place, doubles the estimate of
$\mu$ after each iteration that did not require reducing it converges, on
certain quadratic functions, only at the rate of
gradient descent, $1-\Theta(1/\kappa)$ per iteration, where
$\kappa=L/\mu$. The sequential subspace optimization method
SESOP (a development of A.\,S.~Nemirovski's orthogonal method), which uses
neither $L$ nor $\mu$ and in practice behaves like an accelerated method, also loses
acceleration in the worst case: for specially constructed functions, the residual,
from some step onward, is multiplied at each step by only $(1-3/\kappa)^2$.

The same work shows how to overcome this difficulty without periodic restarts
(cf.\ Remark~\ref{rem:nonmono}). If the estimate of $\mu$ is increased only under
control---within a fixed ``budget'' or below a ``ceiling''
that is lowered with each downward correction---then relative
accuracy $\eps$ is attained in
$2\sqrt\kappa\ln\frac1\eps+\Oh\bigl(\sqrt\kappa\ln^2\kappa\bigr)$
iterations. Subspace optimization methods can also be accelerated without knowing
$\mu$ and $L$ by replacing the slowly forgetting sums of SESOP with memory
that ``forgets'' the past sufficiently quickly: a pair of exponentially
weighted sums of gradients and steps, with a forgetting factor of order
$1-1/\sqrt\kappa$, stores the required horizon of length on the order of $\sqrt\kappa$ in just
two vectors, while the factor itself is chosen during the computation by
an observable progress check. But the difficulty has so far been overcome only
partly: increasing the estimate of $\mu$ does not reduce the leading term
$2\sqrt\kappa\ln\frac1\eps$ (it accelerates convergence along some
directions and slows it along others), and subspace methods
must choose between memory of $\Oh(\ln\kappa)$ vectors and infrequent
momentum resets, at most $\lceil\log_4\kappa\rceil$ of them; whether
bounded memory is possible with no resets at all is unknown.

%% file: en/r11_tensor_superfast.tex

\section{Implementable tensor methods (2018--2021)}
\label{sec:tensor}

\subsection{Background: the third-order frontier}

A natural continuation of the progression ``gradient method--Newton's method'' is methods using derivatives of order $p\ge3$. That such methods \emph{can be very fast} has long been known, and the simplest way to see this is through a classical example: computing the reciprocal $1/a$, $a>0$, without division. Newton's method for the equation $1/x-a=0$ leads to the Newton--Schulz iteration $x_{k+1}=x_k(2-ax_k)$, which uses only multiplication and subtraction. For the residual $u_k=1-ax_k$, this gives $u_{k+1}=u_k^2$, hence $u_k=u_0^{2^k}$: for $|u_0|<1$, the number of correct digits doubles at every step. The origin of the formula is transparent: for $|u_k|<1$,
\[
   \frac1a=\frac{x_k}{1-u_k}=x_k\bigl(1+u_k+u_k^2+\dots\bigr),
\]
and the Newton step truncates this series after two terms. Retain $r+1$ terms:
\[
   x_{k+1}=x_k\bigl(1+u_k+\dots+u_k^{r}\bigr)
   \quad\Longrightarrow\quad u_{k+1}=u_k^{r+1}.
\]
Indeed, $ax_k=1-u_k$, and the geometric sum formula gives
\[
   u_{k+1}=1-ax_{k+1}=1-(1-u_k)\bigl(1+u_k+\dots+u_k^{r}\bigr)
   =1-\bigl(1-u_k^{r+1}\bigr)=u_k^{r+1}.
\]
In other words, $x_{k+1}$ differs from $1/a$ by exactly the discarded tail of
the series, $x_k\bigl(u_k^{r+1}+u_k^{r+2}+\dots\bigr)=u_k^{r+1}/a$: the relative
error is raised to the power $r+1$, and the number of correct digits is multiplied
by $r+1$ at each step. For $r=1$, this is precisely the
Newton--Schulz iteration: $x_k(1+u_k)=x_k(2-ax_k)$. For $r=2$, one obtains a special case of P.\,L.~Chebyshev's method (1838), which uses the second derivative. For an arbitrary equation $f(x)=0$, this approach expands the inverse function $x=\varphi(y)$ in a Taylor series at $y_k=f(x_k)$ and truncates the series:
\[
   x_{k+1}=x_k-\frac{f}{f'}-\frac{f''}{2(f')^3}\,f^2
   -\frac{3(f'')^2-f'f'''}{6(f')^5}\,f^3-\dots
\]
(all derivatives are evaluated at $x_k$). Two terms give Newton's method (order~$2$), three give Chebyshev's method (order~$3$), and in general using derivatives through order $r$ gives a method of order $r+1$. For accuracy $\eps$ in the reciprocal example, one needs
\[
   k\ \ge\ \log_{r+1}\frac{\ln(1/\eps)}{\ln(1/|u_0|)}
\]
iterations: increasing the method's order changes only the \emph{base} of the outer logarithm.

In optimization, M.~Baes~\cite{Baes2009} and E.~Birgin and coauthors~\cite{BirginEtAl2017} obtained bounds for regularized Taylor models, while lower bounds for $p$th-order methods on the class of convex functions with an $L_p$-Lipschitz $p$th derivative take the form
\begin{equation}\label{eq:tenlb}
   N(\eps)\ =\ \Om\!\left(
   \left(\frac{L_pR^{p+1}}{\eps}\right)^{\!\frac{2}{3p+1}}\right),
   \qquad\text{that is,}\qquad
   f(x_k)-f^{*}\ =\ \Om\!\left(\frac{L_pR^{p+1}}{k^{(3p+1)/2}}\right)
\end{equation}
(N.~Agarwal and E.~Hazan~\cite{AgarwalHazan2018};
Y.~Arjevani, O.~Shamir, and R.~Shiff~\cite{ArjevaniShamirShiff2019} for
$p=2$; see also~\cite{GargEtAl2021}). For $\mu$-strongly convex functions in the same class ($p\ge2$), the lower bound is
\begin{equation}\label{eq:tenlbsc}
   N(\eps)\ =\ \Om\Bigl(\Bigl(\frac{L_pR^{p-1}}{\mu}\Bigr)^{\frac{2}{3p+1}}
   +\ln\ln\frac{\eps_0}{\eps}\Bigr),\qquad
   \eps_0=\mu^{\frac{p+1}{p-1}}L_p^{-\frac{2}{p-1}}
\end{equation}
(Y.~Arjevani, O.~Shamir, and R.~Shiff~\cite{ArjevaniShamirShiff2019} for $p=2$; G.~Kornowski and O.~Shamir~\cite{KornowskiShamir2020} for all $p$): the first term is an analogue of~\eqref{eq:tenlb}, and the second corresponds to the very same superlinear local convergence seen in the reciprocal example (for more detail, see \S\,\ref{sec:tensorrates}). Interestingly, lower bounds for tensor methods were first discussed in the book by A.\,S.~Nemirovski and D.\,B.~Yudin~\cite{NemYud1979} (as an exercise for the interested reader!).

However, all such methods were \emph{unimplementable}: the auxiliary problem---minimizing a polynomial of degree $p+1$---is itself a nonconvex global optimization problem, computationally no easier, and in general even harder, than
the original one. The barrier seemed fundamental.

The paper~\cite{Nes2021tensor} (a 2018 CORE preprint) overcame it with a single observation.

\subsection{Convexity of the regularized Taylor model}

\begin{theorem}[\cite{Nes2021tensor}]\label{th:tensorconv}
Let $f$ be convex, with a $p$th derivative that is Lipschitz continuous with constant $L_p$.
Denote by
\[
   \Phi_{x,p}(h)=\sum_{i=0}^{p}\frac{1}{i!}
   D^{i}f(x)[h]^{i}
\]
the Taylor polynomial of degree $p$. If $M_p\ge L_p$, the function
\begin{equation}\label{eq:tenmodel}
   \Omega_{M_p}(x;h)=\Phi_{x,p}(h)
   +\frac{M_p}{(p+1)!}\nrm{h}^{p+1}
\end{equation}
is a global upper model: $f(x+h)\le\Omega_{M_p}(x;h)$ for
all $h$. Moreover, if $M_p\ge p\,L_p$, then
$\Omega_{M_p}(x;\cdot)$ is \emph{convex} in $h$.
\end{theorem}

\begin{proof}[Proof of convexity]
All derivatives below are taken with respect to $h$; note that
$\grad^2\Phi_{x,p}(h)$ is exactly the Taylor polynomial of degree
$p-2$ for the mapping $z\mapsto\grad^2f(z)$, formed at $x$ and
evaluated at $x+h$. Thus the standard Taylor remainder bound gives
\[
   \nrm{\grad^2f(x+h)-\grad^2\Phi_{x,p}(h)}
   \ \le\ \frac{L_p}{(p-1)!}\nrm{h}^{p-1},
\]
and convexity of $f$ implies $\grad^2f(x+h)\succeq0$, whence
\[
   \grad^2\Phi_{x,p}(h)\ \succeq\ -\frac{L_p}{(p-1)!}\nrm{h}^{p-1}I .
\]

On the other hand, direct differentiation gives the following. The Hessian of the regularizer splits
into two parts:
\begin{equation}\label{eq:reghess}
   \grad^2\Bigl[\frac{M_p}{(p+1)!}\nrm{h}^{p+1}\Bigr]
   \ =\ \underbrace{\frac{M_p}{p!}\nrm{h}^{p-1}I}_{\text{isotropic}}
   \ +\ \underbrace{\frac{M_p(p-1)}{p!}\nrm{h}^{p-3}hh^{\!\top}}
        _{\text{rank-one}} .
\end{equation}
We have now reached the key point explaining why the convexity condition contains $pL_p$ rather than
$L_p$---this was precisely the barrier.
The rank-one term acts only along $h$ itself. In the direction
$h$, the curvature of the regularizer is $\frac{M_p}{(p-1)!}\nrm{h}^{p-1}$,
whereas in any orthogonal direction it is only
$\frac{M_p}{p!}\nrm{h}^{p-1}$,
that is, $p$ times smaller. What must be compensated, however, is the negative curvature of
the Taylor polynomial, and it is \emph{isotropic}: the bound
$\grad^2\Phi_{x,p}\succeq-\frac{L_p}{(p-1)!}\nrm{h}^{p-1}I$ applies
in all directions at once. Thus the rank-one part is simply discarded when deriving
the lower bound, and all the work is done by the isotropic
part with coefficient $M_p/p!$. This gives the requirement
$\frac{M_p}{p!}\ge\frac{L_p}{(p-1)!}$, that is, $M_p\ge pL_p$.
\end{proof}

It is now clear that the barrier arose from the properties of the multivariate derivative.
As a check, let us calculate in one dimension: for $n=1$ there are no transverse directions
at all, $\frac{d^2}{dh^2}\frac{M_p|h|^{p+1}}{(p+1)!}
=\frac{M_p}{(p-1)!}|h|^{p-1}$, and convexity requires only
$M_p\ge L_p$.
The factor $p$ is entirely a multidimensional effect.

The proof of Theorem~\ref{th:tensorconv}
takes only a few lines, but it was precisely what, throughout the period from the first bounds for
regularized Taylor models~\cite{Baes2009} to
2018, separated higher-order methods from practice.

Theorem~\ref{th:tensorconv} turns
the auxiliary problem from a nonconvex one into a convex one. More precisely, the problem
$\min_{h\in\R^n}\Omega_{M_p}(x;h)$ becomes a convex optimization
problem
admitting a controlled approximate solution by efficient convex
methods (which ones specifically is discussed in \S\,\ref{sec:superfast}).
Tensor methods became \emph{implementable}---
as reflected in the title of the paper.

\subsection{Convergence rate bounds}
\label{sec:tensorrates}

Four levels of convergence rates must be clearly distinguished here.

\begin{itemize}[leftmargin=1.6em,itemsep=4pt]
\item \emph{The basic method}
      $h_k=\argmin_{h\in\R^n}\Omega_{M_p}(x_k;h)$, $x_{k+1}=x_k+h_k$:
      $f(x_k)-f^{*}=\Oh\bigl(L_pR^{p+1}/k^{p}\bigr)$. This result,
      together with Theorem~\ref{th:tensorconv}, is contained
      in~\cite{Nes2021tensor}.
\item \emph{The accelerated method} (using the estimate sequence
      scheme): $\Oh\bigl(L_pR^{p+1}/k^{p+1}\bigr)$;
      also~\cite{Nes2021tensor}.
\item \emph{A nearly optimal method}, attaining the exponent
      $(3p+1)/2$ from~\eqref{eq:tenlb} up to a logarithmic
      factor (through a nested one-dimensional search in the
      Monteiro--Svaiter scheme): obtained simultaneously and independently
      in~\cite{GasnCOLT2019} and~\cite{Jiang2019,BubeckCOLT2019}; an improved search method was later proposed in~\cite{nesterov2021inexact}. For $p=2$,
      the corresponding construction had been developed earlier by
      R.~Monteiro and B.~Svaiter~\cite{MonteiroSvaiter2013}.
\item \emph{An optimal method}, attaining the exponent
      $(3p+1)/2$ from~\eqref{eq:tenlb} without a logarithmic
      factor (the nested one-dimensional search takes $\Oh(1)$): obtained independently
      in~\cite{KovalevGasnikov2022} and~\cite{CarmonEtAl2022}.
\end{itemize}

Checking consistency of the exponents: for $p=1$, the exponent
$(3p+1)/2=2$ gives exactly the fast gradient
method~\eqref{eq:fgmrate}; for $p=2$, we obtain $k^{-7/2}$, that is,
$N(\eps)=\Oh\bigl((L_2R^3/\eps)^{2/7}\bigr)$, matching the lower bound for second-order
methods in Table~\ref{tab:map}.

In the strongly convex case, a local phase is added to this picture, and
here the reciprocal example again provides an exact model. Near
a nondegenerate minimum, a $p$th-order tensor method converges with
order~$p$~\cite{Doikov2022local}---in agreement with the rule
``derivatives through order $r$ give order $r+1$'', since for the equation
$\grad f(x)=0$ it uses derivatives through order $p-1$---and at
the final stage it needs $\Oh\bigl(\log_p\ln\frac1\eps\bigr)$
iterations. No improvement is possible even in one dimension: already for strongly
convex functions on a half-line, the lower bound for methods of any
fixed order contains a term of order $\ln\ln\frac1\eps$ (this
problem is also treated in~\cite{NemYud1979}), and the method's order affects
only the base of the logarithm---exactly as in Chebyshev's example. In
the multidimensional case, the lower bound~\eqref{eq:tenlbsc} is attained:
an optimal tensor method with restarts in the global phase, combined
with a locally superlinearly convergent method in the final phase, yields
a complexity matching~\eqref{eq:tenlbsc} up to constants and
logarithmic factors independent of $\eps$~\cite{KornowskiShamir2020};
the first term in~\eqref{eq:tenlbsc} corresponds to the global phase, and the second to
the local one.

Thus, the deterministic theory of tensor methods for unconstrained convex
problems has essentially been completed: upper and lower bounds
match up to absolute constants (in the strongly convex
case, up to logarithmic factors independent of
$\eps$).
For a popular account of this story,
from Newton to the papers of 2022,
see~\cite{ZaNauku2023}.

\subsection{Iteration cost and applications}

The practical limit is $p=3$. The reason is the iteration cost, and here we must
return to \S\,\ref{sec:ad}. The third derivative as an object has $n^3$
components, and computing it in full is pointless. But
model~\eqref{eq:tenmodel} requires not the tensors themselves but their \emph{contractions} with
vectors, $D^pf(x)[h,\dots,h]$, and, for fixed
$p$, each such contraction can be computed by automatic differentiation at a cost
proportional to the cost of evaluating $f$ (see part~(ii) of
\S\,\ref{sec:ad}). The circle thus closes: the 1984 result on the low cost of
the gradient becomes the technical foundation of methods of the 2020s.
A systematic review of higher-order methods---their construction, bounds, and
implementation---is given in~\cite{KamzolovEtAl2023}.

A particular application showing where higher-order methods
really gain an advantage is \emph{distributed optimization with
statistical similarity}. If the data, totaling $\Ndata$ observations,
are distributed over $m$ nodes and
partitioned in a statistically homogeneous way, then local functions are close to
one another, and the auxiliary problem arising at each iteration of the
sliding scheme (\S\,\ref{sec:fgm}) contains only
$s\ll \Ndata$ terms. For
$s\sim n$, the cost of one iteration of a second- or third-order method
is determined not by evaluating the Hessian but by inverting it
($\Oh(sn^2)$ versus $\Oh(n^3)$); that is, higher-order
methods ``do not feel'' the finite-sum structure of the problem. This
observation underlies hyperfast local solvers for
statistically preconditioned distributed
optimization~\cite{DvurEtAl2022}.

\practice A limitation should also be noted. All the results listed
concern problems of moderate dimension: for large $n$, direct
inversion of a dense unstructured Hessian becomes
impractical, and tensor methods in this formulation are inapplicable
(sparse and matrix-free implementations constitute a separate regime in which
the range of applicability is determined by matrix structure, not dimension
alone). Their
natural domain is auxiliary problems inside other schemes
(as in~\cite{DvurEtAl2022}) and problems with $n\sim10^2$\,--\,$10^4$ requiring
high accuracy.

\subsection{Current state of the art}
\label{sec:tensornow}

For minimization in Euclidean geometry, as we have seen, the theory of tensor
methods is essentially complete. Beyond these limits---in non-Euclidean
geometry, for saddle-point problems, and with an inexact oracle---the picture
has become clearer only recently.

\smallskip
\emph{Non-Euclidean geometry.} Let the $p$th derivative be Lipschitz continuous
with respect to the $\ell_r$-norm, $1\le r\le\infty$:
$\nrm{D^pf(x)-D^pf(y)}_{*}\le L_p\nrm{x-y}_r$, where
$\nrm{T}_{*}=\max\{|T[h]^p|:\ \nrm{h}_r\le1\}$, and let
$R_r=\nrm{x_0-x^{*}}_r$. The paper~\cite{ContrerasGuzmanMartinezRubio2025}
constructs a method for which
\[
   N(\eps)=\Ot\Bigl(\Bigl(\frac{L_pR_r^{p+1}}{\eps}\Bigr)^{\frac{\rho}{\rho p+p+1}}\Bigr),
   \qquad \rho=\max\{2,r\},
\]
and proves matching lower bounds up to logarithmic factors,
including for randomized methods. For $r\le2$, this is
the exponent in~\eqref{eq:tenlb}; for $r=\infty$, it is $1/p$; and for
$p=1$, these are the classical bounds for accelerated methods in $\ell_r$-geometry.
One ingredient is a non-Euclidean counterpart of Theorem~\ref{th:tensorconv} on
convexity of the regularized Taylor model; in essence, this is the same
problem of choosing a prox-structure as in~\S\,\ref{sec:primaldualgrad}, but
for higher-order methods. The results also cover H\"older-continuous
derivatives.

\smallskip
\emph{Variational inequalities and saddle-point problems.} Let $F$ be
a monotone operator (for example, $F=(\grad_xf,-\grad_yf)$ for
convex-concave $f$), let a $p$th-order oracle return
$F,DF,\dots,D^{p-1}F$, and let $D^{p-1}F$ be Lipschitz continuous with constant $L_p$ (for
a saddle-point problem, these are again derivatives of $f$ through order $p$). Accuracy
is measured by the residual $\dist\bigl(0,F(z)+N_Z(z)\bigr)$, where $N_Z$ is
the normal cone to the feasible set $Z$; $R$ is the distance from
the initial point to a solution. The paper~\cite{JiangYoon2026} develops
an ``anchored'' extra-proximal method combining Halpern iteration
(each step is pulled slightly toward the initial point) with an inexact
proximal step in the spirit of Monteiro--Svaiter, and proves that
the complexity
\[
   N(\eps)=\Ot\Bigl(\Bigl(\frac{L_pR^p}{\eps}\Bigr)^{\frac{2}{3p-1}}\Bigr)
\]
is optimal up to logarithmic factors: a matching lower bound holds for any
deterministic method. For $p\ge2$ and compact convex $Z$, the earlier
work~\cite{ZhangChenPanEtAl2026} attains
$N(\eps)=\Oh((L_pR^p/\eps)^{2/(3p-1)})$ without logarithmic factors.
For saddle-point problems (with $F=(\grad_xf,-\grad_yf)$),
Chen et al.~\cite{ChenZhangLiuEtAl2026} had earlier established a lower
bound with the same exponent for a class of tensor methods.
The paper~\cite{LiEtAl2026} establishes this exponent for unrestricted
deterministic and randomized methods, for both the residual and the gap. The exponent $\frac{2}{3p-1}$ is worse than
the exponent~$\frac{2}{3p+1}$ in~\eqref{eq:tenlb}: for $p=1$, this is $\eps^{-1}$
versus $\eps^{-1/2}$ (the order $\eps^{-1}$ for the gap is achieved by Mirror
Prox and dual extrapolation from
\S\S\,\ref{sec:smooth},~\ref{sec:pd}); for $p=2$, it is $\eps^{-2/5}$ versus
$\eps^{-2/7}$.

\smallskip
\emph{Inexact and expensive second-order oracles.} Three recent papers
show that the optimal dependence $\eps^{-2/7}$ (the case $p=2$
in~\eqref{eq:tenlb}) is retained even when the Hessian is inexact or
expensive; its cost enters the bounds additively or through a factor
growing more slowly than under a naive approach.
\begin{itemize}[leftmargin=1.6em,itemsep=3pt]
\item \emph{A stochastic oracle.} The gradient and Hessian are available with
      errors whose second moments do not exceed $\sigma_1^2$ and
      $\sigma_2^2$. An accelerated variant of the stochastic cubic Newton method
      with a ``contracting'' domain~\cite{AgafonovEtAl2026} yields, after $T$
      iterations,
      \[
         \Exp f(x_T)-f^{*}=\Oh\Bigl(\frac{\sigma_1D}{\sqrt T}
         +\frac{\sigma_2D^2}{T^2}+\frac{L_2D^3}{T^{7/2}}\Bigr),
      \]
      where $D$ is the diameter of the feasible set; this matches the lower
      bound~\cite{Agafonov2023advancing} in all three terms.
\item \emph{Finite sums.} Let $f=\frac1\Ndata\sum_{i=1}^{\Ndata}f_i$,
      with convexity required only of $f$, with Hessians that are Lipschitz continuous in
      mean square,
      $\frac1\Ndata\sum_i\nrm{\grad^2f_i(x)-\grad^2f_i(y)}^2\le
      L_2^2\nrm{x-y}^2$, and with one oracle call returning the gradient and
      Hessian of one term. The paper~\cite{TsedenovKovalev2026} develops
      a method of complexity $\Ot\bigl(\Ndata+\Ndata^{6/7}(L_2R^3/\eps)^{2/7}\bigr)$
      and proves a matching lower bound for span methods; a method
      processing the entire sum at every iteration would require
      $\Oh\bigl(\Ndata\,(L_2R^3/\eps)^{2/7}\bigr)$ calls.
      It can be shown that the same scheme (a reference point at which the
      derivatives of all terms are computed, prediction using its Taylor model, and
      a correction from a randomly chosen term) also works for $p=3$
      if the oracle additionally returns the third derivatives of the term. When the $p$th
      derivatives of the terms are Lipschitz continuous in mean square with
      constant $L_p$, the complexity
      \[
         \Ot\Bigl(\Ndata+\Ndata^{\frac{3p}{3p+1}}
         \Bigl(\frac{L_pR^{p+1}}{\eps}\Bigr)^{\frac{2}{3p+1}}\Bigr),
         \qquad p=2,3,
      \]
      is attained and cannot be improved even for convex terms---at
      least for methods that, like
      span methods~\eqref{eq:span}, discover new coordinates only from
      oracle responses. It suffices to ``hide'' links of the chain
      (cf.~\eqref{eq:worst}) in randomly chosen terms: to discover
      the next coordinate, the method must examine half
      the terms on average. The gain over a full pass is modest---at most
      a factor of $\Ndata^{1/(3p+1)}$---and nonaccelerated variants, while improving over
      the accelerated one in the exponent of $\Ndata$, lose in the exponent of $\eps$ and
      come out ahead only where the term
      $\Ndata$ already dominates.

      But the number of oracle calls is not the best measure for second-order
      methods: the Hessian at the reference point need be
      factorized only once, while each step requires only matrix-vector
      products. Thus (for dense matrices), it is advantageous to
      update the reference point with probability of order $1/(\Ndata+n)$ rather than $1/\Ndata$, and
      for $p=2$ variance reduction is faster than its deterministic
      counterpart in arithmetic operations precisely when
      the heterogeneity of the terms (the ratio of $L_2$ to the Lipschitz constant
      of the Hessian of $f$ itself) does not exceed $\sqrt{\Ndata+n}$ up to
      a logarithmic factor. For saddle-point problems with finite-sum structure, the picture is
      incomplete: if $D$ is the diameter of the feasible set, the complexity
      (up to the additive term $\Ndata$ and logarithmic factors) lies between
      $\Ndata^{3/5}(L_2D^3/\eps)^{2/5}$ and
      $\min\bigl\{\Ndata^{2/3}(L_2D^3/\eps)^{2/3},\,
      \Ndata(L_2D^3/\eps)^{2/5}\bigr\}$; if only the sum is monotone, then for
      methods that discover coordinates only from oracle responses the lower
      bound is higher,
      $\Ndata^{4/5}(L_2D^3/\eps)^{2/5}$, and appears to be tight.
\item \emph{``Lazy'' Hessians} (cf.\ \S\,\ref{sec:cubic}). The Hessian
      is recomputed only once every $q$ iterations (its cost is comparable to $q$
      gradient evaluations), and only the gradient is available in the remaining
      iterations. The paper~\cite{ZhangChenLiuZhang2026} shows that the optimal
      total cost in this model (in gradient evaluations) is
      $\Ot\bigl(q+q^{1/7}(L_2R^3/\eps)^{2/7}\bigr)$, whereas recomputing
      the Hessian at every iteration would cost
      $\Oh\bigl(q\,(L_2R^3/\eps)^{2/7}\bigr)$. The lower bound is obtained
      by a block version of the chain construction~\eqref{eq:worst}: the gradient
      reveals one new coordinate, while the Hessian reveals an entire block of $q$
      coordinates at once.
\end{itemize}

\section{``Superfast'' second-order methods (2020--2021)}
\label{sec:superfast}

\subsection{The apparent paradox and its resolution}

For the class of convex functions with an $L_2$-Lipschitz Hessian, the lower bound
for \emph{second}-order methods is
$f(x_k)-f^{*}=\Om\bigl(k^{-7/2}\bigr)$~\cite{ArjevaniShamirShiff2019,
AgarwalHazan2018}. Nevertheless, Yu.\,E.~Nesterov constructs a method that
uses only gradients and Hessians yet attains the rate
$\Oh\bigl(k^{-4}\bigr)$: the idea itself already appears
in~\cite{Nes2021tensor}, and a detailed treatment is given
in~\cite{Nes2021superfast}.

There is no formal contradiction here, and this should be stated at the outset, not at
the end. The lower bound $\Om(k^{-7/2})$ is proved for the class of functions with
a Lipschitz-continuous \emph{second} derivative ($L_2<\infty$), whereas the method
of~\cite{Nes2021superfast} works on a different class, in which the constant for
the \emph{third} derivative is finite ($L_3<\infty$). The point is that the worst-case
second-order function is assembled from terms of the form $|t|^3$: its second
derivative $6|t|$ is Lipschitz continuous, but its third is $6\sgn(t)$ and has
a discontinuity at zero, so $L_3=\infty$. The worst-case second-order function simply
does not belong to the class ($L_3<\infty$), and the
lower bound built from it consequently ceases to apply. Meanwhile the oracle remains the same:
we still compute only the gradient and Hessian.

The substantive discovery is that second-order information is \emph{sufficient}
to exploit this additional smoothness: oracle order and
class smoothness are independent parameters and need not match.

Why had this not been noticed earlier? Because such a device is impossible for first-order
methods, and it was there that the habit was formed.
The worst-case first-order function is quadratic, and all derivatives of a quadratic
of order $\ge3$ are identically zero, so $L_p=0$
for every $p\ge2$. No additional assumption of higher
smoothness excludes it from the class, and the lower bound $\Om(k^{-2})$
remains valid.

The general picture is determined by parity. The worst-case function for
$p$th-order methods is built from terms $|t|^{p+1}$. For \emph{odd}
$p$, the exponent $p+1$ is even, and $|t|^{p+1}=t^{p+1}$ is an ordinary polynomial
whose constants $L_{p'}$ vanish for all $p'>p$; such a function
remains admissible in every higher-smoothness class, and the lower bound remains valid. For \emph{even} $p$, the exponent $p+1$ is odd, the function
$|t|^{p+1}$ loses smoothness at zero already in its $(p+1)$st derivative, and
the worst-case construction does not belong to a higher-smoothness class. The first two cases are
exactly $p=1$ (a quadratic function, nothing can be done) and $p=2$
(the function $|t|^3$, acceleration is possible).

This contrast nicely illustrates how subtle
the complexity theory of convex optimization is, and how carefully one must
specify the exact class on which a lower bound is proved.

\subsection{The mechanism}
We now turn to implementing a method that achieves the desired acceleration through higher smoothness.
The scheme consists of two ideas.

\smallskip
\emph{(1) The auxiliary problem via relative smoothness.}
The resulting convex problem---minimizing a fourth-degree polynomial
with regularizer $\nrm{h}^4$---is solved not by an ``ordinary'' gradient
method (its gradient Lipschitz constant is unbounded), but by a gradient
method \emph{relative} to the reference function $\eta(h) = \tfrac12\grad^{2} f(x)[h]^2 + \tfrac{L_3}{4}\nrm{h}^4$.
Recall the definition: a function $\phi(h)$ is called \emph{relatively $L$-smooth}
and \emph{relatively $\mu$-strongly convex} with respect to a convex reference function $\eta(h)$ if
$L\grad^2\eta(h) \succeq \grad^2\phi(h) \succeq \mu \grad^2\eta(h) $. The ordinary gradient step is then replaced by a step in
the Bregman geometry generated by $\eta(h)$, and the linear convergence rate is preserved.
The notion was introduced by H.~Bauschke, J.~Bolte, and
M.~Teboulle~\cite{BBT2017} and, independently, by H.~Lu, R.~Freund, and
Yu.\,E.~Nesterov~\cite{LuFreundNes2018}; the latter paper provides
exactly the tool needed here. Note the continuity of ideas: relative
smoothness directly develops the philosophy of self-concordance from
\S\,\ref{sec:ipm}, where the ``goodness'' of a function is also measured not by an absolute
constant but relative to an appropriate reference function.

\smallskip
\emph{(2) A third-order tensor step from a second-order oracle.}
Relative gradient descent steps require the Hessian $\grad^2 f(x)$ in the Bregman distance and only the gradient of the model~\eqref{eq:tenmodel} itself. Thus, only contractions of the third
derivative $D^{3} f(x)[h]^2$ are needed; they can be computed exactly by automatic differentiation at a cost comparable to that of evaluating the gradient, or approximated inexactly by finite differences. The
paper~\cite{Nes2021superfast} uses the \emph{central difference of
gradients}:
\begin{equation}\label{eq:centraldiff}
   D^3f(x)[h]^2\ \approx\
   \frac{\grad f(x+\theta h)+\grad f(x-\theta h)-2\grad f(x)}{\theta^2},
\end{equation}
for which, under a Lipschitz-continuous third derivative, the remainder is of order
$\Oh\bigl(L_3\theta\nrm{h}^3\bigr)$.

In~\cite{nesterov2021inexact}, Yu.\,E.~Nesterov proposed a proximal construction that does not even require evaluating or approximating the third directional derivative.

Thus, we need only one full Hessian evaluation for the Bregman distance and a logarithmic number of operations, each costing on the order of one gradient evaluation.
A further development is ``hyperfast'' methods attaining
$\Ot(k^{-5})$, that is, exponent $(3p+1)/2$ for $p=3$;
see~\cite{nesterov2021inexact, KamzolovGasnikov2020}.

%% file: en/r12_concl.tex
\section{Conclusion}
\label{sec:concl}

The results discussed above concern very different topics---from automatic differentiation to traffic equilibria, from semidefinite relaxations to tensor methods. Nevertheless, behind them lies the unified program formulated in \S\,\ref{sec:lookinside}: \emph{looking inside the black box}. Let us bring its outcomes together.

\subsection{Lower bounds and difficult problems}

Almost every work of Yu.\,E.~Nesterov begins with the question ``how fast is it possible to go \emph{in principle}?'' This standpoint is inherited from the complexity theory of optimization problems developed by A.\,S.~Nemirovski and D.\,B.~Yudin, and turns algorithm design from an art into a science: a method is considered understood only when accompanied by a proof that one cannot do better---in the specified problem class, with the specified oracle, and for the specified criterion. Each of these three qualifications is essential, and most misunderstandings in the literature arise from their omission.

A by-product of this approach is the \emph{technique} of constructing worst-case examples, that is, examples on which a method attains the lower bounds. The chain construction~\eqref{eq:worst} has proved to be a universal tool: today its variants yield lower bounds for high-order methods, coordinate descent and gradient-free schemes, distributed optimization, saddle-point problems, and variational inequalities (\S\,\ref{sec:complexity}).

\subsection{Looking into the black box}

The fast gradient method (\S\,\ref{sec:fgm}) attains the lower bound in the same information model, and there is nothing to improve on this class. This, however, is the only topic discussed here that has this structure. In all the others, the method converges \emph{faster} than the lower bound permits---and there is no contradiction, since the bound is proved for a broader problem class or for a less informative oracle. Technically, this means that the method is supplied with an \emph{objective function model} (\S\,\ref{sec:models}):

\smallskip
\begin{itemize}[leftmargin=1.6em,itemsep=2pt]
\item a \emph{geometric} model---an explicit upper model or an explicitly inscribed simple set: self-concordant barriers (\S\,\ref{sec:ipm}), cubic (\S\,\ref{sec:cubic}) and tensor (\S\,\ref{sec:tensor}) upper models, and an upper model for the residual norm (\S\,\ref{sec:gn});
\item an \emph{algebraic} model---a max representation of the function: the smoothing technique (\S\,\ref{sec:smooth}), primal-dual and barrier subgradient methods (\S\,\ref{sec:pd});
\item a \emph{structural} model---the interaction of simple and complex objects: composite functions (\S\,\ref{sec:composite}), relative smoothness (\S\,\ref{sec:superfast}), and generalized smoothness (\S\,\ref{sec:l0l1}).
\end{itemize}

\smallskip
To these are added four techniques:
\begin{itemize}[leftmargin=1.6em,itemsep=2pt]
\item taking the \emph{structure} of the problem into account when smoothing (\S\,\ref{sec:smooth});
\item moving to higher-order derivatives (\S\,\ref{sec:tensor});
\item \emph{narrowing the class} through higher smoothness while keeping the oracle order unchanged (\S\,\ref{sec:superfast});
\item separating different types of oracle calls and calling each at its own frequency (\S\S\,\ref{sec:fgm},~\ref{sec:cd}).
\end{itemize}
In each case, the key is to specify precisely \emph{which} part of the problem setting has changed; without this specification, the phrase ``the method is faster than the lower bound'' is meaningless.

Note that the program is not complete, and its author himself assesses the situation quite sharply: once the pessimistic bounds obtained in the black-box model had been circumvented in many directions at once, the theory of convex optimization was no longer covered by a unified system of lower bounds, and the development of new principles for classifying problems remains open~\cite{NesDiss2013}.

An important qualification concerning gradient-free methods (\S\,\ref{sec:zo}): they should \emph{not} be included in this list. Weakening the oracle does not accelerate but \emph{worsens} the minimax complexity---by a factor of $n$ in the smooth Euclidean case, and this factor is unavoidable in a moderate accuracy range. The correct general statement is that complexity theory allows us to measure the \emph{value} of different types of information: the gradient is worth a factor of $n$ compared with function values, the second derivative reduces the exponent from $1/2$ to $2/7$, and so on.

\subsection{The right notion of regularity}

A separate line of research concerns the search for structural properties and regularity conditions that guarantee the computational solvability of a problem. Self-concordance (\S\,\ref{sec:ipm}), H\"older continuity of the gradient (\S\,\ref{sec:univ}), gradient domination (\S\,\ref{sec:cubic}), relative smoothness (\S\,\ref{sec:superfast}), convexity of the regularized Taylor model (\S\,\ref{sec:tensor}), and generalized $(L_0,L_1)$-smoothness (\S\,\ref{sec:l0l1})---all these concepts have proved fruitful.

Here it is important not to succumb to the temptation to subsume them all under the term ``affine-invariant.'' That would be incorrect. It is precisely self-concordance that has \emph{strong} affine invariance: inequality~\eqref{eq:sc} transforms into itself under any nonsingular substitution $x=By$, and a method built on it is independent of the choice of coordinates. Relative smoothness is invariant under a \emph{consistent} transformation of the pair ``function $+$ reference function.'' The remaining conditions---H\"older continuity with constant $L_\nu$, strong convexity with parameter $\mu$, Lipschitz continuity of the $p$th derivative with constant $L_p$, and the Polyak--\L{}ojasiewicz condition with parameter $\tau$---\emph{depend on the chosen norm} and change under a linear transformation of coordinates. The choice of norm affects not only the constants: if smoothness and distance to the solution are measured in different norms, even the exponent in the optimal bound changes (\S\,\ref{sec:primaldualgrad}). The correct formulation of the general observation is as follows: successful notions of regularity relate the behavior of a function not to an absolute scale but to some \emph{reference}, be it the function itself (self-concordance), a chosen reference function (relative smoothness), the gradient norm ($(L_0,L_1)$-smoothness), or the optimality gap (gradient domination).

\subsection{Duality as a computational tool}

In the stable dynamics model (\S\,\ref{sec:traffic}), dual variables are queues and optimal tolls; in the smoothing technique (\S\,\ref{sec:smooth}), the dual problem provides both a method and an accuracy certificate; in primal-dual methods (\S\,\ref{sec:pd}), averaging the solutions of the inner problems delivers a primal solution together with a verifiable error bound at no additional cost. In Yu.\,E.~Nesterov's work, duality is always substantive rather than formal: it not only reformulates the problem but also operates within the method.

\subsection{From theory to everyday practice}

Let us list the influences in verifiable terms.

\begin{itemize}[leftmargin=1.6em,itemsep=3pt]
\item Nesterov acceleration and its variants (including what is called Nesterov momentum in neural network training libraries) are part of the standard optimization toolkit for training models; it is nevertheless incorrect to identify every ``momentum'' method with scheme~\eqref{eq:fgm}.
\item The dual averaging method~\eqref{eq:da} and adaptive regularized dual schemes (RDA, FTRL, AdaGrad) are closely related constructions with similar analyses; they provide the theoretical basis for the family of adaptive methods to which Adam also belongs.
\item Interior-point methods based on the theory of \S\,\ref{sec:ipm} are implemented in CVX, CVXPY, MOSEK, SeDuMi, and SDPT3, and determine which convex problems are considered solvable in practice.
\item Coordinate descent methods are implemented in \texttt{glmnet} and \texttt{LIBLINEAR} and are standard in applied statistics.
\item Gaussian smoothing \eqref{eq:gsmooth}\,--\,\eqref{eq:twopoint} underlies evolution strategies in reinforcement learning and gradient-free fine-tuning of language models.
\item Backpropagation is relation~\eqref{eq:adjoint}, and Theorem~\ref{th:ad} explains why training neural networks is possible at all.
\end{itemize}

It is precisely this density of mathematical constructions in everyday computational practice that accords with the citation of the Gauss Prize, which recognizes mathematical work with outstanding applications outside mathematics.

\subsection{The research school}

Finally, something that cannot be reduced to theorems. Through his postgraduate studies under B.\,T.~Polyak, his many years of collaboration with A.\,S.~Nemirovski, his work at CEMI, and subsequently through PreMoLab, MIPT, the Higher School of Economics, and Innopolis, Yu.\,E.~Nesterov became a link between the Soviet school of optimization and its present, largely international continuation. A substantial proportion of the results mentioned above as the ``current state of the art'' were obtained by his students, coauthors, and coauthors' students---one need only look at the bibliography. The award of the 2026 Gauss Prize can also be viewed as recognition of the scientific tradition that has formed around these works.

%% file: en/r13_bib.tex
\renewcommand{\refname}{References}